\documentclass[12pt,letterpaper]{amsart}

\usepackage{nomencl}
\usepackage[makeindex]{imakeidx}
\usepackage{amsmath,amsthm, amsbsy,amsfonts,amssymb, txfonts}
\usepackage[normalem]{ulem}
\usepackage{caption}
\usepackage{stmaryrd,mathrsfs,graphicx, supertabular, amscd, tikz-cd, tikz}
\usepackage{stackrel}
\usetikzlibrary{matrix,arrows,decorations.pathmorphing}

\DeclareFontEncoding{LS1}{}{}
\DeclareFontSubstitution{LS1}{stix}{m}{n}
\DeclareSymbolFont{symbols2}{LS1}{stixfrak} {m} {n}
\DeclareMathSymbol{\operp}{\mathbin}{symbols2}{"A8}

\usepackage[active]{srcltx}
\allowdisplaybreaks
\numberwithin{equation}{section}

\usepackage[
	hypertexnames=false,
	hyperindex,
	pagebackref,
	breaklinks=true,
	bookmarks=false,
	colorlinks,
	linkcolor=blue,
	citecolor=red,
	urlcolor=red,
]{hyperref}
\usepackage{hyperref}
\usepackage{nomencl}
\def\BB{{\mathbb B}}

\def\CC{{\mathbb C}}

\def\FF{{\mathbb F}} 
 
\def\HH{{\mathbb H}}
\def\II{{\mathbb I}}
\def\Hb{{\mathbb H}}
\def\KK{{\mathbb K}}

\def\LL{{\mathbb L}}
 
\def\MM{{\mathbb M}} 
\def\NN{{\mathbb N}} 
\def\PP{{\mathbb P}}

\def\QQ{{\mathbb Q}}

\def\ZZ{{\mathbb Z}}

\def\ssm{\smallsetminus}

\def\G{\Gamma}
\def\g{\gamma}

\def\odd{{\rm odd}} 
\def\even{{\rm even}} 
\def\orb{{\rm orb}} 
\def\irr{{\rm irr}} 
\def\spl{{\rm sp}} 
\def\lr{{\rm lr}} 
 
\def\reg{{\rm reg}} 
\def\refl{{\rm refl}} 
\def\sg{{\rm sg}} 
\def\ss{{\rm ss}} 
\def\std{{\rm std}}

\def\tbase{{\rm base}}
\def\tfiber{{\rm fiber}}

\def\A{{\rm A}}
\def\ct{{\rm ct}}
\def\DM{{\rm DM}}

\def\Hyp{{\mathcal H}\!y\!p}

\def\K{{\rm K}}

\def\int{{\rm int}}

\def\cox{{\rm Cox}} 
\def\st{{\rm st}}
\def\sing{{\rm sing}}

\def\bs{\backslash}
\newcommand{\eps}{\varepsilon}

\newcommand{\p}{\partial}

\def\Ccal{{\mathcal C}}
\def\Dcal{{\mathcal D}}
\def\Ecal{{\mathcal E}} 
 
\def\Gcal{{\mathcal G}} 
\def\Hcal{{\mathcal H}}

\def\Kcal{{\mathcal K}}
\def\Lcal{{\mathcal L}}
\def\Mcal{{\mathcal M}}

\def\Ocal{{\mathcal O}}
\def\Pcal{{\mathcal P}}

\def\Scal{{\mathcal S}}

\def\Qcal{{\mathcal Q}} 
\def\Ucal{{\mathcal U}}  
  
\def\Wcal{{\mathcal W}}  
\def\Xcal{{\mathcal X}}  
\def\Ycal{{\mathcal Y}}

\def\Mb{{\mathbf M}}

\def\Cscr{{\mathscr C}}

\def\Lscr{{\mathscr L}}

\def\Xscr{{\mathscr X}}

\def\la{\langle}
\def\ra{\rangle}
\def\half{{\tfrac{1}{2}}}

\def\disc{\mathrm{disc}}

\def\Hfrak{\mathfrak{H}}
\def\Sfrak{\mathfrak{S}}
\def\Ufrak{\mathfrak{U}}

\def\pt{{\scriptscriptstyle\bullet}}

\newcommand\arf{\operatorname{Arf}}
\newcommand\aut{\operatorname{Aut}}
\newcommand\out{\operatorname{Out}}

\newcommand\coker{\operatorname{coker}}

\newcommand\Hl{\operatorname{H}}
\newcommand\Hom{\operatorname{Hom}}

\newcommand\inn{\operatorname{inn}}

\newcommand\proj{\operatorname{Proj}}
\newcommand\pic{\operatorname{Pic}}
\newcommand\rk{\operatorname{rk}}
\newcommand\sym{\operatorname{Sym}}
\newcommand\GL{\operatorname{GL}}

\newcommand\PGL{\operatorname{PGL}}
\newcommand\PG{\operatorname{P\Gamma}}

\newcommand\Gg{\operatorname{G}}

\newcommand\PGtilde{\operatorname{P\tilde\Gamma}}

\newcommand\Lk{\operatorname{Lk}}

\newcommand\SL{\operatorname{SL}}

\newcommand\Tcal{\operatorname{{\mathcal T}}}

\newcommand\Orth{\operatorname{O}}

\newcommand\PSL{\operatorname{PSL}}
\newcommand\U{\operatorname{U}}

\newcommand\PU{\operatorname{PU}}
\newcommand\Sp{\operatorname{Sp}}

\newtheorem{theorem}{Theorem}[section]
\newtheorem{lemma}[theorem]{Lemma}
\newtheorem{proposition}[theorem]{Proposition}
\newtheorem{corollary}[theorem]{Corollary}

\newtheorem{conjecture}[theorem]{Conjecture}

\newtheorem{definition}[theorem]{Definition}

\theoremstyle{remark}
\newtheorem{convention}[theorem]{Convention}
\newtheorem{example}[theorem]{Example}
\newtheorem{examples}[theorem]{Examples}

\newtheorem{remark}[theorem]{Remark}
\newtheorem{remarks}[theorem]{Remarks}

\newtheorem{conventions}[theorem]{Conventions}
\newtheorem{conventionsnot}[theorem]{Notational Conventions}

\newcommand{\EL}[1]{\textcolor{red}{#1}}%can be used  for comments
\newcommand{\edit}[1]{\textcolor{blue}{#1}}  %correction
\title{Cubic threefolds moduli and the Monster group}
\begin{document}

\author{Eduard Looijenga}
\address{Mathematics Department, University of Chicago (USA) and Mathematisch Instituut, Universiteit Utrecht (Nederland)}
\email{e.j.n.looijenga@uu.nl}


\begin{abstract}
Allcock constructed a 13-dimensional complex ball quotient of which he conjectured that it admits a natural  covering with covering group isomorphic to the Bimonster (the wreath product of the Monster group and the group of order two). This ball quotient contains the moduli space of cubic threefolds as an open dense subset
of a 10-dimensional complex subball quotient. We prove that this subball quotient has a neighborhood in the Allcock ball quotient over which the conjectured cover exists.
\end{abstract}

\maketitle
\section*{Introduction}
 Allcock's  Monstrous proposal \cite{allcock:proposal}  concerns  a  $13$-dimensional ball quotient, 
he introduced. It takes the form of a conjecture, recalled below, which  connects  this ball quotient  with 
the largest sporadic simple  finite group, the Monster  $\Mb$ 
\nomenclature[A]{$\Mb$}{the Monster group}%
(also known as the  Friendly Giant), or rather with 
$\Mb\wr \mu_2 (:=\Mb^2\rtimes \mu_2)$, 
\nomenclature[A]{$\Mb\wr \mu_2$}{the Bimonster group}%
a group  Conway calls the \emph{Bimonster}. 
If true, then this would undoubtedly contribute to a better geometric understanding of both $\Mb$ and 
$\Mb\wr \mu_2$. For example, both groups come with a natural conjugacy class of involutions 
and this would realize these involutions as reflections.

This ball  quotient   has the structure of a Shimura variety and can be regarded as a moduli space for 
certain  polarized Hodge structures invariant under an automorphism of order 3.  
This suggests that if the conjecture holds, it might lead, as Allcock wonders, to  a modular interpretation 
of this ball quotient.  If true, this would make  fascinating  connections  between certain moduli spaces  and 
other sporadic finite simple  groups.  

This $13$-dimensional ball quotient is known to  contain a closed  subvariety which  is uniformized  by a subball  of dimension 10 and  which has an open-dense subset  the moduli space of cubic threefolds. The main result of this paper  may be summed up  by saying that it proves Allcock's  conjecture for  a regular neighborhood of this subvariety. This requires  us to establish properties of the moduli space of cubic threefolds that were not previously known or noted.
\\

In order to state the conjecture, we need to recall what  Eisenstein lattices are. Suppose $L$ is an even integral lattice and  $\rho\in \Orth(L)$ an orthogonal transformation
 satisfying $\rho^2+\rho +1=0$ (in other words, $\rho$ has order $3$ and its only fixed point is the origin). Then $L$ becomes a module over the Eisenstein ring  $\Ecal:=\ZZ[u]/(u^2+u +1)$, which we regard as
embedded in $\CC$ by sending $u$ to  $\omega:=e^{2\pi\sqrt{-1}/3}$. Since  $\Ecal$ is a Euclidean ring  for the norm, it is  a  principal ideal domain and hence the  $\Ecal$-module $L$ is free.  We convert the given form on $L$  into an $\Ecal$-valued hermitian form on this module by putting 
\[
\la x,y\ra =\theta\big((x\cdot \rho y)-\omega (x\cdot y)\big),
\]
where $\theta=\omega-\omega^2=\sqrt{-3}$
 (so this form is actually $\theta\Ecal$-valued; without the factor $\theta$ 
it would be skew-hermitian).  This is what we shall call an \emph{Eisenstein lattice}.  
When an even lattice is considered as an Eisenstein lattice,  we indicate this by  using  the corresponding blackboard font, in this case $\LL$. 
We regard $\LL$ as sitting in the complex vector space  $V(\LL):=\CC\otimes_\Ecal \LL$, 
to which we have extended $\la \; ,\, \ra$  as a hermitian form.  
The signature of this form is half that  of $L$ (which is indeed necessarily even). 
The case of interest here is when the signature of $L$ is of the form $(2n+2,2)$, for then  
$V(\LL)$ has \emph{Lorentzian} signature $(n,1)$. This implies that in the associated complex projective space is defined an $n$-dimensional  complex ball $\BB(\LL)$ 
\nomenclature[A]{$\BB$}{a complex ball associated with a Lorentzian complex vector space}%
whose points  represent the negative one-dimensional subspaces of $V(\LL)$. This ball is a hermitian symmetric domain  for the group $\PU(V(\LL))$ of projective unitary transformations. We may think of this  group as  the group of real points of an algebraic group over the Eisenstein field $\QQ(\omega)$. It contains the image of $\U(\LL)$ as an arithmetic subgroup. This subgroup is discrete and hence  acts properly discontinuously on $\BB(\LL)$.
This is of course still true for a finite index subgroup $\G$ of $\U(\LL)$.
The Baily-Borel theory tells us that the resulting  orbifold  $\Xcal:=\G\bs \BB(\LL)$
(a ball quotient) has a natural structure of a normal projective variety minus a finite subset,  the set of  \emph{cusps}.

In order to describe Allcock's ball quotient, we need the notions of triflection and mirror.
Any $r\in L$ with $r\cdot r =2$ defines an orthogonal  reflection $s_r$ in $L$ by the formula  $x\mapsto  x-(x\cdot r)r$. A small computation shows  that then $\rho(r)\cdot r=-1$ (so that $(r,\rho (r))$ is a root basis of type $A_2$) and that
\[
t_r(x):= s_{\rho (r)}s_r
\]
will preserve  $\LL$ and act on $\LL$ as a  unitary pseudo-reflection of order three, called by Allcock a \emph{triflection}.
The fixed point set of $t_r$ in $\BB(\LL)$ is a hyperplane section  (which is defined by the orthogonal complement of $r$ in $\LL$, which has signature $(n-1,1)$), called the \emph{mirror}   of $t_r$ with $t_r$ acting in the normal direction of hyperplane section   with eigenvalue $\omega$. One checks that $\la r,r\ra =3$, and because of this, we  refer to $r$ as a \emph{$3$-vector}. 
The triflections in the $3$-vectors make up finitely many $\U(\LL)$-conjugacy classes and their  mirrors make up a locally finite collection of hyperplane sections of $\BB(\LL)$. 

We get Allcock's ball quotient if we take  for $L$ a unimodular lattice of signature $(26,2)$. It is well-known that this lattice is unique. For example it can be obtained as the orthogonal direct sum of one of the 24  Niemeier lattices and 2 copies of the standard hyperbolic lattice of rank $2$. Allcock shows that $L$ admits an orthogonal transformation $\rho$ as above and that  all such $\rho$ make up a single conjugacy class in $\Orth(L)$. Hence the  resulting Eisenstein lattice $\LL$ is unique up to isomorphism; we will denote it here by $\LL^\A$.
  Allcock and Basak prove that 
$\U(\LL^\A)$ is generated by its  triflections, that  the triflections make up a single conjugacy class in $\U(\LL^\A)$ (\cite{basakI}, \cite{allcock:y555}), and that 
the sum  of all the triflection mirrors in $\BB(\LL^\A)$ is the zero divisor of a $\U(\LL^\A)$-automorphic form of weight 4 relative to some 
$\mu_3$-valued character (\cite{allcock:leech}, Thm.\ 7.1). 

The $\U(\LL^\A)$-action on $\BB(\LL^\A)$ defines  an orbifold $\Xcal$ and the sum $Z$ of the triflection mirrors  determines 
in $\Xcal$ a (discriminant) hypersurface $D_\Xcal$. This hypersurface is $2$-divisible in the Picard group of $\Xcal$, so that a double cover $\tilde \Xcal\to \Xcal$ with ramification locus $D_\Xcal$ exists. The latter comes with a natural (effective) orbifold structure and 
Allcock's conjecture is equivalent to the assertion  that the orbifold fundamental group of  $\tilde \Xcal$  is isomorphic to $\Mb\times\Mb$, or equivalently, that the covering  group $\hat G$ of $\hat \Xcal\to \tilde \Xcal\to \Xcal$ is isomorphic to the 
Bimonster $\Mb\wr \mu_2 $. 

All automorphisms of the Monster $\Mb$ are inner, and from this one derives that the involutions of $\Mb\wr \mu_2 $ over $\mu_2$ are the ones of the form  
\[
\iota_a: (x,y)\in\Mb^2\mapsto (aya^{-1} , a^{-1}xa)\in \Mb^2.
\]
 with  $a$ running over all elements of $\Mb$. The centralizer of such an involution is isomorphic to 
$\mu_2\times \Mb$, where the first factor is $\la\iota_a\ra $ and the second is the graph of $\inn(a)$ in $\Mb\times\Mb$.

Allcock and  Basak  very recently showed \cite{ab} that $\hat G$ is a quotient of $\Mb\wr \mu_2 $. Since $\Mb$ is simple, this implies  that either Allcock's conjecture is true or that $\tilde \Xcal$ is simply connected, which, needless to say,  would be disappointing. 
Our theorem (and the many remarkable `coincidences' that come along with it) is a strong indication that this is not the case.
\\

In order to explain the results in this paper, we must first discuss the case where the lattice $L$ we started with is positive definite.  
There are  four basic cases, which  can be characterized  as follows. Among the root lattices 
of type $A$, $D$ or $E$ (where all roots have squared norm $2$), there are four  with the property that some power of a Coxeter transformation   
has order  $3$ and turns that root lattice  into positive definite Eisenstein lattice: these  are   
$A_2$, $D_4$, $E_6$ and $E_8$.  We shall denote these  Eisenstein lattices $\LL_1$,  $\LL_2$, $\LL_3$ and $\LL_4$ respectively 
(so $\LL_i$ has rank $i$ over the Eisenstein ring).   
The associated triflections generate a subgroup $\Gg(\LL_i)$ of the finite group  $\U(\LL_i)$, which of 
course appears in the Shephard-Todd classification of finite pseudo-reflection groups.  

Allcock showed in \cite{allcock:y555} that every positive definite Eisenstein lattice $\MM$ 
generated by its $3$-vectors decomposes 
canonically into a direct sum of these lattices.  We can emulate the above  program for such an $\MM$ by  
replacing  $\BB(\LL)$ by  the complex vector space  
$V(\MM)=\CC\otimes_\Ecal\MM$ and   $\G$  by the (finite) subgroup $\Gg(\MM)\subset\U(\MM)$ 
generated by the 
triflections in $\MM$: We take for $Z$  the union of the triflection hyperplanes  (denoted here by $Z(\MM)$) and let the ${\Gg(\MM)}$-orbit space of the pair $(V(\MM),Z(\MM))$ take the role of the pair $(\Xcal, D_\Xcal)$. A classical theorem of Chevalley 
asserts that $\Gg(\MM)\bs V(\MM)$ is isomorphic to affine space
and that ${\Gg(\MM)}\bs  Z(\MM)$  is the discriminant  of the quotient  map $V(\MM)\to \Gg(\MM)\bs V(\MM)$. 

If we take $\MM=\LL_i$, then the  discriminant $\Gg(\LL_i)\bs Z(\LL_i)\subset \Gg(\LL_i)\bs V(\LL_i)$ turns out to be also the discriminant 
of the root system of type $A_i$.  This was first observed by Orlik-Solomon \cite{os}, but a natural proof  of this remarkable fact was later found by Couwenberg in his thesis \cite{couw}.  Since the Weyl group of that root system is a copy of the symmetric group 
$\Sfrak_{i+1}$, it follows that  $\Gg(\LL_i)\bs V(\LL_i)$ admits an $\Sfrak_{i+1}$-covering ramified over $\Gg(\LL_i)\bs Z(\LL_i)$ whose total space is again smooth.  For general $\MM$ we are simply dealing with  a product of these cases and so then the covering will also be smooth and the covering group, which we will denote by $\hat G(\MM)$, 
will be a product of  permutation groups of degree $\le 5$.
\\

We return to the Allcock lattice $\LL^\A$. Suppose $\MM\subset \LL^\A$ is a primitive  positive definite 
sublattice whose $3$-vectors generate a sublattice of finite index.  
Then $\BB(\MM^\perp)\subset \BB(\LL^\A)$ is an intersection of mirrors and any nonempty intersection of mirrors is  thus obtained. 
We expect that if we restrict the $\hat G $-covering $\hat \Xcal\to \Xcal$  to a transversal of  the image of 
$\BB(\MM^\perp)$ in $\Xcal$, then we obtain an embedding of $\hat G (\MM)$ 
(so a product of  permutation groups) in  $\hat G $ and,  as we shall now explain,  the  
$\hat G $-centralizer of this embedding to have an interpretation very much like that of  $\hat G $ itself.

We begin with noting that the  natural map $\Gg(\MM)\to \U(\LL^\A)$ is an embedding. The centralizer of its image, which we shall denote by $\G(\MM^\perp)$,  is a finite index subgroup of $\U(\MM^\perp)$. So the action of $\G(\MM^\perp)$ on $\BB(\MM^\perp)$ defines  a ball quotient  $\Xcal(\MM^\perp)$.
The  evident map of ball quotients  $\Xcal(\MM^\perp)\to \Xcal$ is finite.  Let $\hat \Xcal(\MM^\perp)\to \Xcal(\MM^\perp)$ be the finite Galois cover obtained by   pulling  back the 
$\hat G $-cover $\hat \Xcal/\Xcal$ along this map, normalize its total space and  take  a 
connected component. 

This also suggests a bijective correspondence between 
\begin{enumerate}
\item[($A$)] $\U(\LL^\A)$-orbits  of primitive  positive definite sublattices $\MM\subset \LL^\A$ whose $3$-vectors generate a sublattice of finite index and 
\item[($B$)] conjugacy classes of subgroups of the Bimonster that are products of good (in the sense of \cite{cns}) permutation groups of degree $\le 5$.
\end{enumerate}
This might be sufficiently concrete  as to be  verifiable with the aid of a computer. The correspondence  should however lift one level up: if we think of the collection ($A$) as labelling the set of strata of $\Xcal$, then the collection ($\hat A$) of  strata of the $\hat G $-cover $\hat \Xcal$  should be in bijective correspondence with the collection ($\hat B$) of members of  conjugacy classes appearing in ($B$). The  main result of this paper (which we shall state below) implies that this stronger version holds in some sense if in ($A$) we restrict to $\MM$ that contain a copy of $\LL_3$ and in ($B$)  to good subgroups that contain  a  good $\Sfrak_4$.
\\

Let us illustrate this with the  basic examples that are  pertinent here.
The first one is when $\MM\cong\LL_4$, a case that was already considered by Heckman \cite{heckman}.
Then $\MM^\perp\cong  2\LL_4\operp\HH$, where $\HH$ is the hyperbolic Eisenstein  lattice 
(it has by definition  isotropic generators $e, f$ with $\la e,f\ra=\theta$). The lattice  $\MM^\perp$ is well-known,  
as it appears in the Deligne-Mostow list. The associated ball quotient (of dimension 9) has the modular interpretation 
as the moduli space of Hilbert-Mumford stable positive degree 12 divisors on a smooth genus zero curve, 
or equivalently (and more to the point here) as a moduli space of hyperelliptic curves of genus $5$. 
As Heckman observed, the associated reflection cover is obtained by numbering the points of the divisor 
and hence the covering group  is of type $\Sfrak_{12}$.  This is compatible with the observation by  
Conway-Norton-Soicher \cite{cns} that the Bimonster $\Mb\wr\mu_2$ comes with a distinguished conjugacy class of 
subgroups isomorphic with  $\Sfrak_5$  whose centralizer in $\Mb\wr\mu_2$ is isomorphic to $\Sfrak_{12}$.  

We get two other examples leading to ball quotients with a known modular interpretation if we take  
$\MM\cong\LL_3\operp\LL_1$ resp.\ $\MM\cong\LL_3$ (so one contains the other). In these cases,  
$\MM^\perp$ is isomorphic to resp.\  $\LL_1\operp\LL_3\operp\LL_4\operp\HH$ resp.\  
$\LL_1\operp 2\LL_4\operp\HH$ and then the associated ball quotients of dimensions $9$ resp.\ $10$ 
contain Zariski open-dense subsets that have the interpretation of the moduli space of genus 4 curves (by a theorem of Kond\=o) 
resp.\ the moduli space of cubic threefolds (due independently to Allcock-Carlson-Toledo \cite{act} and 
Looijenga-Swierstra \cite{ls}). We prove (Corollary \ref{cor:ct2})  that the reflection 
covering groups of these ball quotients are given by the representations on the mod 2 homology of a 
genus four curve resp.\ a cubic threefold, which means that the groups in question  are isomorphic to 
$\Sp_8(\FF_2)$ resp.\  $\Orth^1_{10}(\FF_2)$ (the orthogonal group of a nonsingular quadratic form 
over $\FF_2$ in 10 variables of Arf invariant $1$). This, too is in agreement with what  Conway-Norton-Soicher 
\cite{cns} found:  they  proved that $\Mb\wr\mu_2$ comes with a distinguished conjugacy class of 
subgroups isomorphic with  $\Sfrak_4$  whose centralizer in $\Mb\wr\mu_2$ is isomorphic to 
$\Orth^1_{10}(\FF_2)$. This implies that the centralizer of  $\Sfrak_2\times\Sfrak_4$ in the Bimonster is isomorphic 
to $\Sp_8(\FF_2)$. Our main result, Theorem \ref{thm:main}, states that the $\Mb\wr\mu_2$-cover of $\Xcal$ 
predicted by Allcock exists on a  neighborhood of the image of this $10$-dimensional ball quotient in $\Xcal$.

Our main result  fits in a program  described in the short note \cite{looij:heegner} (that is in fact excised from an earlier version of this paper). 
We show there (by means of a weak Lefschetz theorem) that for the determination of the reflection covering groups 
such as $\hat G $ mentioned above, it suffices to do this for a regular neighborhood of what we call the 3-skeleton 
(a union of $3$-ball subquotients). Unfortunately, the closure of the cubic threefold locus  in the Allcock ball quotient does not 
contain the 3-skeleton (otherwise we would be done), but it is likely that the `missing' $3$-ball quotients 
also have a modular interpretation (it seems quite feasible to do this for the 3-skeleton) and this makes us cautiously optimistic that this will eventually lead to a proof of the conjecture. 

We note that the   main  result of the recent Allcock-Basak paper mentioned above implies  that the inclusion of a regular neighborhood of the closure of the cubic threefold locus in the Allcock ball quotient (with their modified  orbifold structures) is surjective on orbifold fundamental groups. 
Allcock's conjecture would of course follow if one can show that this induces an isomorphism on orbifold fundamental groups.
\\

Before we  review the individual sections, we  point  the reader to  some questions connected with this material.

The first is that another potential  example waiting to be probed  is when the role of  $\LL_3$ is taken by  $\LL_1\operp\LL_2$. This lattice embeds primitively in $\LL_4$ 
with orthogonal complement $\LL_1(2)$ (the lattice of rank one has a generator of square norm $6$).
We embed $\LL_1\operp\LL_2$ via this embedding in $\LL^\A$. Then $(\LL_1\operp\LL_2)^\perp\cong \LL_1(2)\operp 2\LL_4\operp\HH$ 
determines another  ball quotient of dimension 10. It has in common with the 10-dimensional ball parametrizing cubic 
threefolds that it contains the 
Deligne-Mostow ball quotient of dimension 9. Taking our cue from  \cite{cns}, we expect that the reflection 
group associated  to this ball quotient involves another sporadic  finite group, namely $\mathrm{Fi}_{23}$; the challenge to algebraic geometers is to find a modular interpretation of it.

The second question regards automorphic forms. We recall that a  $\U(\LL^\A)$-automorphic form of weight $w$ 
(with trivial character) is simply a $\G$-invariant holomorphic section of the automorphic line bundle
$\Ocal_{\BB(\LL)}(-w)$ (assuming, as we will,  that the rank $n+1$ of $\LL$ is $\ge 3$; otherwise a growth condition must be imposed). This means that the line bundle $\Ocal_{\BB(\LL)}(-1)$ descends to an ample bundle $\Lcal$ on $\Xcal$  in the orbifold sense so that the $\G$-automorphic forms of weight $w$ are the sections of $\Lcal^{\otimes w}$. If $\Gcal$ is finite, then this implies  that  $\hat\pi^*\Lcal$ is also ample and 
\[
\oplus_{w=0}^\infty \Hl^0(\hat \Xcal, \hat\pi^*\Lcal^{\otimes w})
\]
is a finitely generated graded $\CC$-algebra endowed with a $\Gcal$-representation. This may well be worth studying in its own right, even for the
ball quotients with a known modular interpretation. If one can come up with a  decent proposal (in a Moonshine spirit, say) for what to expect for the Allcock ball quotient (or for  one of its mirrors modulo its stabilizer) together  with  a candidate of a graded  $\CC$-algebra $A^\pt$ as above with a faithful $\Mb\wr\mu_2$ resp.\ $\Mb$-action, then  its subalgebra of invariants must  reproduce  $\oplus_{w=0}^\infty \Hl^0(\Xcal, \Lcal^{\otimes w})$ (whose proj defines the Baily-Borel compactification $\Xcal^*$), so that  $\proj(A^\pt)$  is the sought-for  $\hat \Xcal$ (or rather the normalization of  $\Xcal^*$ in $\hat \Xcal$).
\\

We have divided the paper into two parts (but  numbered the sections  consecutively).
Part 1 already contains the proof of the main theorem, but depends on some properties of the moduli space of cubic threefolds, which are  
needed at the very end (Theorems \ref{thm:ctcover} and  \ref{thm:match0}). By  thus giving  algebraic geometry a bit of a black box treatment we hope  to make an essential  part  of this paper more accessible for (and perhaps appealing to) researchers mostly  working in (finite) group theory,  lattices or  ball quotients.  Part 2 is entirely about properties of the moduli space of cubic threefolds. It reviews previous work  and also proves many new properties about it, which we believe have an interest in their own right. Its sole goal,  however, is  to establish  Theorems \ref{thm:ctcover} and  \ref{thm:match0} of Part 1. 

The first half of Section \ref{sect:quadratic forms} collects known facts regarding quadratic forms on $\FF_2$-vector spaces and the second half focuses on a particular case,  relevant for this paper.

Section \ref{sect:Bimonstercomplex} gives a characterization of the Bimonster in a way that is useful for our purpose. It is based on the work of a group of 
Conway, Norton 
and many others in the 1990s. We use the work in Section \ref{sect:quadratic forms}, to classify the conjugacy classes of  subgroups  of the Bimonster  (here denoted by $\Gcal$) which we have called \emph{geometric}: these are direct products of groups whose factors are isomorphic with $\Sfrak_4$ or $\Sfrak_5$ (\emph{good}  in the sense of Conway \emph{et al.}). Subgroups of $\Mb\wr\mu_2$ of this type form a partially ordered set $\Cscr$ for the inclusion relation and  hence define a simplicial complex $\Sigma(\Cscr)$ which comes with a $\Mb\wr\mu_2$-action. 
One of the key results here is that this  simplicial complex is simply connected.
 
Section \ref{sect:universalcovers} introduces  the notion of a reflection cover. Its appearance  is so ubiquitous that it is a bit surprising  that it  had not been formalized before. We give several examples; many of these (or special cases thereof) show up later in this  paper.

Section \ref{sect:eisenstein}  reviews the general theory of Eisenstein lattices that we shall need (and adds to this a bit as well), states 
Allcock's conjecture and introduces the Eisenstein lattices that are central to what follows.

Section \ref{sect:appC}  classifies the sublattices of the Allcock lattice that are isomorphic to $p\LL_3\operp q\LL_4$. We call such lattices \emph{geometric}  for  a reason.
Such sublattices form a partially ordered set $\Lscr$ and hence also gives rise to a simplicial complex $\Sigma(\Lscr)$. This complex  comes with  an action of the unitary group of the Allcock lattice (here denoted by $\G$). We find that the quotients $(\Mb\wr\mu_2)\bs \Sigma(\Cscr)$ and 
$\G\bs \Sigma(\Lscr)$ (which are still simplicial complexes) are isomorphic. The proof  involves no algebraic geometry at all (it is essentially a proof by inspection) and deprived of its context, appears to be a miraculous coincidence. This  is however a consequence  of the validity of Allcock's conjecture, and so  this may be regarded as evidence for it.

In the final section \ref{sect:neighborhood} of Part 1, we state what we need to know about the moduli space of cubic threefolds (Theorems \ref{thm:ctcover} and  \ref{thm:match0}) and subsequently  derive from this our main result.

The first section of Part 2 (Section \ref{sect:cubic3fold} ) is about the 10-dimensional  ball quotient of which an open subset parametrizes  cubic threefolds. 
Although there is quite a bit of literature on this topic, we have not found there the results we need in a 
directly quotable  form. %In the case where the moduli space of genus 4 curves is concerned we  also  somewhat deviated from the  way this is presented in the work of Kond\=o \cite{kondo:g=4} and  Casalaina-Jensen-Laza \cite{cjl}.  

In Section \ref{sect:level2}  we show that by imposing on a smooth cubic threefold $X$  a full level 2 structure on $\Hl_3(X)$,  we obtain a reflection cover of its moduli space. We then prove  that this cover  becomes simply connected if we extend it as a ramified cover over the moduli space of stable  smooth cubic threefolds. This property has an important consequence for the associated 10-dimensional ball quotient, for it gives us  Theorem \ref{thm:ctcover}, which we do not know how to obtain without this modular interpretation.
We then also have almost all what is needed  for the proof of our main Theorem  \ref{thm:match0}.  

The proof itself occupies Section \ref{sect:coarsestrata} and  involves the following geometric result: if a cubic threefold $X$ has  only $A_3$-singularities (it can have at most 4 such), then  we prove that the property for two such  singularities to  be coplanar, i.e., lie in a plane  contained in $X$,  is an equivalence relation on its singular set.  We find that this equivalence relation separates the connected components of the loci in the moduli space of cubic threefolds with a fixed number of $A_3$-singularities. Although this involves small sets, here again there is a remarkable match: this set of  connected components is 
 in bijection  with  the set of  $\G$-orbits of geometric sublattices of the Allcock lattice isomorphic with a direct sum of copies of $\LL_3$
(in a way that links coplanarity with imprimitivity). This leads us  to the proof of Theorem  \ref{thm:match0}. 
\\
% We found it striking how well the complex geometry here dovetails with the characterization of the Bimonster due to  Conway-Norton-Soicher \cite{cns} and Conway-Pritchard \cite{cp} (proved by Ivanov \cite{ivanov} and Norton \cite{norton}) and we hope the reader concurs.

\emph{Acknowledgements.} A paper of this length builds on the work of many, many  others.  This is especially so  for the inspiring work  by the sporadic group community  on the Monster and allied groups. Although this was   done several decades ago (see the references), part of it turns out to be tailor made for the kind of applications that appear here.

Gert Heckman told me about the Allcock conjecture about ten years ago. Although my interest  remained passive, we had some interesting exchanges, which I am sure helped me when things changed and  I developed a more active interest  in this at the end of  2022. This first led to a proof that the Deligne-Mostow 9-ball quotient has a regular neighborhood in the Allcock ball quotient which admits a reflection cover with the Bimonster as covering group. At the time I even believed that this would imply the Allcock conjecture, but both Gert and  Daniel Allcock quickly corrected me on this.
Subsequent correspondence with Daniel and  Tathagata Basak and our getting together for a few days at Iowa State University helped me to further shape ideas. It was thus that the suggestion came up to extend the (considerably  easier) argument for the Deligne-Mostow ball quotient to the cubic threefold 3-ball quotient. I am grateful  to all three colleagues.

I am especially grateful to Leonard Soicher for  feedback regarding  Subsection \ref{subsect:bimonster}. He alerted me to an incomplete argument in the original version and numerous other comments of his helped me to improve the presentation. I thank Stephen Humphries for consultation in connection with Section \ref{sect:level2}.

The input of  referees  of an earlier version  not only led to an improved exposition, but also prompted me to rearrange the material and  to alter (and correct)  some of the proofs. After I had completed this version, I used Claude to check  for mathematical correctness (a test the paper passed) and for presentation (which led to minor changes).
%\tableofcontents

\begin{conventions}
In this paper most of our geometric objects are  orbifolds in the conventional sense (of Thurston, say), 
but taken  in  the complex-analytic or complex-algebraic category. This notion still allows for  nontrivial 
global inertia groups, so that if  $\G$ is a group acting properly discontinuously on a complex manifold $M$, 
then the orbit space $\G\bs M$ is in  a natural manner a complex-analytic orbifold: if $p\in M$ and $U\ni p$ is 
a neighborhood of $p$ with the property that 
each $\G$-orbit in $M$  which meets $U$, meets $U$ in a $\G_p$-orbit, then the map 
$U_p\subset M\to  \G\bs M$ is an orbifold chart. We shall denote this orbifold by $M_\G$.
We will encounter certain  Deligne-Mumford stacks over $\CC$ and we indicate this in general by an underlined symbol. Such 
stacks have  an underlying orbifold and we will sometimes not distinguish the two when there is no risk of confusion.

We adhere to the convention (used in \cite{atlas} for example) to denote a possibly nonsplit extension of a group $G$ by a group $N$ as $N.G$ (a split extension is a semi-direct product and usually  denoted $N\rtimes G$).
\end{conventions}
\tableofcontents
\newpage 

 \begin{center}
 \textbf{PART I: THE BIMONSTER REALIZED AS A REFLECTION GROUP}
 \end{center}
 
\section{Quadratic forms over $\FF_2$}\label{sect:quadratic forms}

In this section all vector spaces are  finite dimensional and over $\FF_2$. 
The first subsection reviews standard properties of quadratic forms  over $\FF_2$ and in doing so establishes some notation. In the next subsection  we focus on the cases that are particularly relevant for this paper. The results therein may be  less standard and will be  used in an essential way in Section \ref{sect:Bimonstercomplex}. 

\subsection{Basics}Let $V$ be a finite dimensional vector space over $\FF_2$. 
We consider bilinear forms $(x,y)\in V\times V\mapsto (x\cdot y)\in \FF_2$ that are \emph{alternating},  
i.e., which vanish on the diagonal (this implies that $x\cdot y=y\cdot x$). 
The  \emph{radical} of such a form is the subspace $K_V\subset V$ of $v\in V$ for which $v\cdot x=0$ for all $x\in V$.
If that kernel is trivial, then the form is called \emph{symplectic} and we denote the group of linear transformations of 
$V$ that preserve it by $\Sp(V)$. This forces $V$ to have even dimension, $2m$ say,  and it then admits a dual basis 
$(x_{\pm 1} ,\dots, x_{\pm m})$ on which it takes the standard form $x\cdot y=\sum_{|i|=1}^g x_i y_{-i}$ (no need for signs here, as $-1=1$ in $\FF_2$).

So in general, $V/K_V$ inherits from $V$ a symplectic form.
We then  denote the group of linear transformations of $V$ that preserve the form and act as the identity on
$K_V$ by  $\Sp_o(V)$. The obvious map $\Sp_o(V)\to \Sp(V/K_V)$ is onto and its kernel can be identified with the vector group
$V/K_V\otimes K_V$ which acts unipotently by assigning  to $a\otimes  n$ the transformation
\[
T(a\otimes n): x\mapsto x+(x\cdot a)n
\]
Every $a\in V$ defines a \emph{symplectic transvection} 
$\tau\in\Sp_o(V)$ by
\[
\tau_a(x)= x+ (x\cdot a)a
\]
and these generate $\Sp_o(V)$ (notice  that if  $n\in K_V$, then $\tau_a\tau_{a+n}=T(a\otimes n)$).

If $T\subset V$ is a subspace, then we write $T^\perp$ for the set if $v\in V$ for which $v\cdot t=0$ for all $t\in T$. So $K_V=V^\perp$.
In general,  $K_T=T\cap T^\perp$ and so we have a decomposition $V=T\operp T^\perp$ precisely when $T$ is symplectic.

A \emph{quadratic form} on $V$ is a function   $q:V\to \FF_2$ for which $q(0)=0$ and 
\[
(x,y)\in V\times V\mapsto q(x+y)-q(x)-q(y)\in \FF_2
\] 
is bilinear; we then call $V$ or rather the pair $(V,q)$ a \emph{quadratic space}.  Its  associated bilinear form is then alternating and the restriction  $q|K_V$ is linear. If that restriction  is zero, then $V/V_K$ inherits a quadratic form  from $V$; otherwise
it only inherits a symplectic form. 

We say that $q$  is \emph{nonsingular}  if the associated alternating form is symplectic, i.e., if $K_V=0$. 
There are two isomorphism types of nonsingular quadratic forms in 2 variables:  
one is represented by the hyperbolic form $x_1x_2$  and the other by the form
$x_1^2+x_1x_2+x_2^2$ (which we will denote by $A_2$). 

A linear map $V\to V'$  between two quadratic spaces is called a \emph{morphism of quadratic spaces} if its composite with the quadratic form on $V'$ yields the quadratic form on $V$. Such a map is injective if $V$ is nonsingular.
 
 We will call $a\in V$ a  \emph{root} (in the belief that this will not cause confusion with 
 other uses of that term) when $q(a)=1$. In the hyperbolic case, only $(1,1)$ is a root, whereas in the $A_2$-case, all three nonzero elements are roots.
For a root $a$ the 
associated symplectic transvection preserves $q$ and is therefore called an 
\emph{orthogonal reflection}; we then denote it $\sigma_a$. 
So $\sigma_a(x)=x+(x\cdot a)a$ and $\sigma_a$ lies in the group $\Orth_o(V)$ of orthogonal 
transformations that act as the identity on the radical. If that radical is zero, then by a theorem of Dieudonn\'e (see \cite{chevalley}, p.\ 20), reflections in the roots  
generate $\Orth(V)$, except when $V$ is isomorphic to the orthogonal sum of two copies of $A_2$; the 
orthogonal reflections then generate a copy of  $\Sfrak_3\times \Sfrak_3$ and hence do not contain an 
orthogonal transformation that exchanges the two $A_2$ copies.

There are two isomorphism types of   nonsingular quadratic spaces of  a given (necessarily even) dimension. 
These can be distinguished by their  \emph{Arf invariant}.  This is an element of $\FF_2$ that takes the
value  $0$ or $1$  according to whether $\#( q^{-1}(0))>\#( q^{-1}(1))$ or $\#( q^{-1}(0))<\#( q^{-1}(1))$.
We denote the isomorphism type of the orthogonal group of a nonsingular quadratic form
of dimension $n$ and of Arf invariant $\eps\in\FF_2$ by $\Orth^\eps_n(\FF_2)$  
(\footnote{This should not be confused with  the notation used in the 
$\mathbb{ATLAS}$  \cite{atlas},  where   $\Orth_n^+(\FF_2)$ resp.\  $\Orth_n^-(\FF_2)$ 
stands for the (simple) subgroup of $\Orth^0_n(\FF_2)$ resp.\   
$\Orth^1_n(\FF_2)$ of transformations that can be written as a  word in an even number of  reflections.}).

Both $\Sp(V)$  and  $\Orth(V)$ have the structure of  a Chevalley group: each is the group of $\FF_2$-valued points of an algebraic group defined over $\FF_2$.

\begin{lemma}\label{lemma:embeddingcriterion}
Let  $V$ be a nonsingular quadratic space and $T$ a quadratic space whose quadratic form vanishes of its radical. 
Then an orthogonal embedding  $T\hookrightarrow V$ exists if and only if $\dim V\ge \dim T+\dim K_{T}$ and, in case of equality, in addition $\arf(T/K_{T})=\arf(V)$.  These  embeddings then lie  in a single $\Orth(W)$-orbit. 
\end{lemma} 
\begin{proof}
We prove the existence part. Write $I$ for $K_{T}$. Choose a supplement $T_1$ of $I$ in $T$ so that $T_1\cong T/I$. 
If we endow $I\oplus I^*$ with the (hyperbolic) quadratic form $(x, \xi)\mapsto \xi(x)$, then this space is nonsingular and has Arf invariant $0$. 
So $T_2:=T\oplus I^*=T_1\operp (I\oplus I^*)$ is nonsingular and has the same Arf invariant as $T/I$. So when
 $\dim V= \dim T+\dim I$ and $\arf(T/I)=\arf(V)$, then $V$ must be isomorphic with $T_2$. When 
 $\dim V>\dim T+\dim I$, then choose a nonsingular quadratic space $T_3$ of dimension  $\dim V-\dim T-\dim I$ with Arf invariant $\arf(V)-\arf(T_1)$
 and then $T_2\oplus T_3\cong V$. In either case we get a quadratic embedding of $T$ in $V$.

The proof that all embeddings so arise,  walks through this argument  in the opposite direction and is left to the reader.
The last clause then also follows.
\end{proof}

Let $T\subset V$ be as in the conclusion of Lemma \ref{lemma:embeddingcriterion}: $V$ is nonsingular and $T\subset V$ a subspace for  which the quadratic form is identically zero on  $I:=T\cap T^\perp(=K_{T})$.  
We want to describe  the structure of the pointwise stabilizer $\Orth(V)_{T}$ in $\Orth(V)$ and its relation to the subset stabilizer  $\Orth(V)_{[{T}]}$. The latter is a simple type of parabolic subgroup of $\Orth(V)$  and what follows is a special case  of  the theory of  such subgroups of orthogonal Chevalley groups.

We first note that ${T}/I$ and  $I^\perp/{T}$ are nonsingular quadratic spaces.  The natural  homomorphism $\Orth(V)_{T}\to \Orth (I^\perp/{T})$ is onto (formation of  the Levi quotient) and its kernel, the unipotent radical $R_u(\Orth(V)_{T})$,   consists of \emph{Eichler transformations}: these associate with  a pair $(e,f)\in I\times {T}^\perp$  
the orthogonal  transformation
\[
\exp(e\wedge f): x\mapsto x+(x\cdot e)f+(x\cdot f)e + q(f)(x,e)e, 
\]
(thus denoted, because this transformation indeed only depends on $e\wedge f\in \wedge^2V$ and $\exp$ functions here as an $\FF_2$-analogue of the  exponential map defined on  a Lie algebra).  These transformations make up a unipotent group with  center 
the vector group $\wedge^2 I$ and with  quotient the vector group $I\otimes ({T}^\perp/I)$. The commutator map 
factors through the  alternating map $I\otimes ({T}^\perp/I)\times  I\otimes ({T}^\perp/I)\to \wedge^2I$ which 
is obtained  by contraction with the help of the  quadratic form on  ${T}^\perp/I$. When  $\dim I=2$, this is a Heisenberg group (which we shall denote 
$\Hfrak({T}^\perp/I)$,  as we can take for its argument  any nonsingular quadratic space). We record this as follows.

\begin{corollary}\label{cor:parabolic}
Let ${T}\subset V $ be a subspace on the nonsingular quadratic space $V$ such that the quadratic form is zero on its radical $I=T\cap T^\perp$. Then the map $\exp$ as defined above identifies  $\wedge^2 I$ with a central subgroup of  $\Orth(V)_{T}$. The  quotient  is  a semidirect product of a reductive group by a vector group, 
\[
\Orth(V)_{T}/\exp(\wedge^2 I)\cong  (I\otimes ({T}^\perp/I))\rtimes \Orth ({T}^\perp/I), 
\]
where the action of $\Orth ({T}^\perp/I)$ on $I\otimes ({T}^\perp/I)$ is the obvious one.
Similarly, 
\[
\Orth({T})\cong (T/I\otimes I)\rtimes \big(\Orth({T}/I)\times \GL(I)\big).
\]
The  restriction map  $\Orth(V)_{[{T}]}\to \Orth({T})$ is onto with kernel $\Orth(V)_{T}$ and 
the  outer action of $\Orth({T})$ on $\Orth(V)_{T}$ is via the  Levi factor $\GL(I)$ of the former on the 
unipotent radical $R_u(\Orth(V)_{T}\cong \wedge^2 I. (I\otimes ({T}^\perp/I))$ of the latter in the obvious manner.
\end{corollary}
\begin{proof}
The  first paragraph is a rewording of the discussion which preceded the corollary. The  proof of the second paragraph is similar, but easier.
The proof of the last assertion is an exercise.
\end{proof}

\begin{remark}\label{rem:parabolic}
We can also state this more in terms of abstract groups, for $\Orth(V)_{T}$ is the $\Orth(V)$-centralizer of $\Orth ({T})$ and $\Orth(V)_{[{T}]}$ is the $\Orth(V)$-normalizer of $\Orth ({T})$. Then  $\Orth ({T})$ appears as the quotient of these two. If in the last clause of Corollary \ref{cor:parabolic} we exchange  the roles of ${T}$ and ${T}^\perp$, we  find that  the outer action of $\Orth(V)_{[{T}]}$ on $\Orth ({T})$ is through its Levi factor $\GL(I)$.
\end{remark}

\subsection{The quadratic spaces of type $A_n$}\label{subsect:stabilizer} Let $E$ be a nonempty finite set. We regard $\FF_2^E$ as the space of characteristic functions of subsets of $E$ and denote by $H_E\subset \FF_2^E$ the subspace  spanned by characteristic functions of even sized subsets of $E$. An $\FF_2$-valued quadratic form $q_E:H_E\to \FF_2$ is induced by the 
$\frac{1}{2}\ZZ$-valued  form on $\ZZ^E$ given by $\frac{1}{2}\sum_{e\in E}\chi_e^2$. In other words,  it takes on 
the characteristic function $\chi_A$ of a subset $A\subset E$ the value $0$ or $1$ according to the 
parity  of half its size $\#( A)$ and the associated bilinear form assigns  to a pair of even sized subsets the parity of 
the size of their intersection.  In particular, a root in $H_E$ is given by a subset of $E$ whose size is $2\pmod{4}$.

The alternating form on $H_E$ is nondegenerate (and so the quadratic form nonsingular) when $\#( E)$ is odd.
When $E$ is nonempty and of even size, then the radical of the alternating form is spanned by $\chi_E$ and so 
$H_E/\la \chi_E\ra$ comes with a nondegenerate alternating form. We denote that symplectic space by $H'_E$.  
It is a quadratic space precisely when   $q(\chi_E)=0$ and this is equivalent with $\#( E)\equiv 0\pmod{4}$. 
Observe  that this construction comes with an embedding of the permutation  group $\Sfrak_E$ in $\Orth(H'_E)$.

When $E$ has  even size and   $e\in E$, then the  natural map $H_{E\ssm \{e\}}\to H'_E$ is an isomorphism of symplectic spaces.

When $E=\{0,1,\dots, n\}$  we also write  $A_n$ for  $H_E$ and $\Orth (A_n)$ for $\Orth(q_E)$.
This choice of notation  is motivated by the fact  that the $2$-element subsets of  $E$ define roots and that these are  also the roots (in the more conventional sense of the word) of a subsystem of type $A_n$, where opposite roots  are identified. 
If we number the elements of $E$: $e_0, \dots ,e_n$, then 
\[
\{\alpha_i:=\chi_{e_i}-\chi_{e_{i-1}}\}_{i=1}^n
\]
is a root basis of that type. 
The orthogonal reflections in such a basis elements generate the permutation group $\Sfrak_{n+1}$,  here realized as a subgroup of $ \Orth(A_n)$.  This is in general a proper subgroup, for  a root  need not be given by $2$-element subset of $E$.  When  $E$ has size $\le 6$,  this is however the case and hence   
 $\Orth(A_n)\cong \Sfrak_{n+1}$ when $n\le 5$. 
 
For $n\ge 5$, the orthogonal complement of the obvious inclusion $A_4\subset A_n$ has as basis the roots  $e_0+e_1+e_2+e_3+e_4+e_5$ and 
$e_5-e_6, \dots , e_n-e_{n-1}$. These make up a copy of $A_{n-4}$. So if we write $n=4m+i$ with $i\in \{0,1,2,3\}$, then  $A_n\cong A_i\operp m A_4$.
Let us therefore focus on the first four cases.
\begin{itemize}
\item[$A_1$:] $q_E$ takes the value 1 on the generator and represents the form $x^2$.
\item[$A_2$:] $q_E$ takes the value $1$ on its nonzero elements (it represents the form $x_1^2+x_1x_2+x_2^2$); so it  is nonsingular  of Arf invariant $1$.
\item[$A_3$:] $q_E$ is singular with $\chi_E$ as its unique element in the radical. If we divide out by the radical we get a form of type $A_2$. Note that Corollary \ref{cor:parabolic} gives us the Chevalley group description of 
 $\Orth(A_3)=\Sfrak_4$ as the semidirect product $A_2\rtimes \Orth(A_2)\cong\FF_2^2\rtimes \Sfrak_3$
 (the affine transformation group of an affine plane).
\item[$A_4$:] The partial basis  $(\alpha_1, \alpha_2)$ splits off an $A_2$-type as a direct summand.
Its orthogonal complement has basis $(\alpha_1+\alpha_3, \alpha_4)$ which represents the form $x_3x_4+x_4^2$. The latter  is nonsingular  of Arf invariant $0$
and hence, by the  additivity of the Arf invariant,  $q_E$ is nonsingular  of   Arf invariant 1. 
%It represents the form $x_1^2+x_1x_2+x_2^2+x_3x_4+x_4^2$. 
\end{itemize}
It follows that $A_{4m}$ and  $A_{4m-2}$ are nonsingular with Arf invariant $m\pmod{2}$. We also note $A_{4m+1}$ is the direct sum of its  radical (spanned by $\chi_E$) and $A_{4m}$.

\begin{definition}\label{def:geomsubgroup1}
Given a nonsingular quadratic space $V$, we say that a subgroup of $\Orth(V)$ isomorphic  with $\Sfrak_{k+1}$ 
is \emph{of reflection type}  if it arises from an embedding of $A_{k}$  in $V$.
We say that it is \emph{strictly of reflection type} if $k\le 4$ and that  it is \emph{geometric} if $k\in \{3,4\}$. 
More generally, a subgroup of  $\Orth(V)$  isomorphic with a product of permutation groups is said to be 
(strictly) of reflection type or  geometric if every factor is.

We write  $\Cscr_{\Tcal}$ for the geometric conjugacy classes of $\Orth(\Tcal)$.
\end{definition}

The reason that we here single out the $A_k$ with $k\le 4$ has ultimately to do with the fact that the singularities of stable 
cubic threefolds  are of that type (but now taken in the sense of singularity theory). 
This can also be justified in intrinsic terms, as these are the cases for which the inclusion  
$\Sfrak_{k+1}\subset \Orth(A_k)$ is an equality.  The special role we attribute to the cases 
$k\in \{3,4\}$ will  become  clear later.

So a subgroup of reflection type isomorphic with $\Sfrak_{i_1+1}\times\cdots\times  \Sfrak_{i_r+1}$ comes from a 
morphism of quadratic spaces $j: A_{i_1}\operp \cdots \operp A_{i_r}\to V$. 
This map need not be an embedding: the radical  $K_{A_n}$ of $A_n$  is of dimension 1 when $n$ is 
odd and trivial otherwise and  the kernel of this map can be any subspace 
$K\subset K_{A_{i_1}}\oplus \cdots \oplus K_{A_{i_r}}$ which not contain any nonzero summand.  %If $i_t\le 4$ for all $i$, then  $\Orth(A_{i_t})\cong \Sfrak_{i_{i+1}}$ and so the reflection subgroup  is then completely determined by the map $j$. 

Our main interest will be in the situation when $E$ has size $12$ (so that would be $A_{11}$). The  quadratic space 
\begin{equation}\label{eqn:defT}
\Tcal:=H'_E:=H_E/\la\chi_E\ra
\end{equation}
 is then nonsingular of rank 10 and has 
Arf invariant $1$. Indeed, it is isomorphic with $A_{10}$, but  the construction  endows it with an embedding of 
$\Sfrak_{12}$ in $\Orth(H'_E)\cong \Orth^1_{10}(\FF_2)$.
We encounter such $\Tcal$  when we are dealing with a hyperelliptic curve $C$ of genus $5$ and $E$ is its  set  of Weierstra\ss\ points: Then we have a natural  isomorphism $\Tcal\cong H^1(C; \FF_2)$ of symplectic vector spaces
(note that a difference of two  Weierstra\ss\ points can be understood as an element of order two of the Jacobian of $C$). Via  this isomorphism the symplectic form on $H^1(C; \FF_2)$  lifts  to a quadratic form. This quadratic form coincides with the one  given  by a natural (odd) theta characteristic on $C$, namely the divisor class defined by $4$ times a Weierstra\ss\ point (which is independent of the Weierstra\ss\ point).

Given a   $W\in \Cscr_{\Tcal}$,  we will write $W^\perp$ for its centralizer in $\Orth(\Tcal)$. So its $\Orth(\Tcal )$-normalizer can be written as
\[
N_{\Orth(\Tcal)}(W)=(W\times W^\perp).\out_{\Orth(\Tcal)}(W),
\]
where $\out_{\Orth(\Tcal)}(W)$  is the image of $N_{\Orth(\Tcal)}(W)$ in $\out(W)$. We shall see that the  action of $\out_{\Orth(\Tcal)}(W)$ on
$W^\perp$ may be nontrivial, but need not be faithful.

We classify such $W$ up to conjugacy  and determine $N_{\Orth(\Tcal)}(W)$  and its constituents.

\begin{proposition}\label{prop:a3classification}
A geometric reflection subgroup $W\in \Cscr_{\Tcal}$  is given by an orthogonal map 
$rA_3 +sA_4\to  \Tcal $ whose kernel $K(W)\subset K_{A_3}^r$ does not contain a summand. Its $\Orth(\Tcal)$-conjugacy class  is completely determined by the position of $K(W)$ in the direct  sum $K_{A_3}^r$ and we have the following cases:
\begin{enumerate}
\item[$A_3$:] an $\Orth(\Tcal)$-orbit  of embeddings $A_3\hookrightarrow \Tcal$ ($K(W)=0$), 
\item[$2A_3$:] an $\Orth(\Tcal)$-orbit  of embeddings $2A_3\hookrightarrow \Tcal$ ($K(W)=0$), 
\item[$(2A_3)'$:] an $\Orth(\Tcal)$-orbit  of the image of an orthogonal map $2A_3\to \Tcal$  with  $K(W)$  the diagonal in $2K_{A_3}$, 
\item[$3A_3$:] an $\Orth(\Tcal)$-orbit  of the image of an orthogonal map $3A_3\to \Tcal$  with  $K(W)$  the main diagonal in 
$3K_{A_3}$, 
\item[$(2A_3)'+A_3$:] an $\Orth(\Tcal)$-orbit  of the image of an orthogonal map  $3A_3\to \Tcal$ with  $K(W)$ the diagonal of the first two summands of $3K_{A_3}$,  
\item[$(3A_3)'$:] an $\Orth(\Tcal)$-orbit  of the image of an orthogonal map  $3A_3\to \Tcal$ whose  kernel $K(W)\subset 3K_{A_3}$ is the span of the diagonals of each pair of summands (so is of dimension  2). 
\item[$A_4$:] an $\Orth(\Tcal)$-orbit  of embeddings $A_4\hookrightarrow \Tcal$,
\item[$A_3+A_4$:] an $\Orth(\Tcal)$-orbit  of embeddings $A_3+A_4\hookrightarrow \Tcal$,
\item[$2A_4$:] an $\Orth(\Tcal)$-orbit  of embeddings $2A_4\hookrightarrow \Tcal$,
\item[$(2A_3)'+A_4$:] an $\Orth(\Tcal)$-orbit  of embeddings $(2A_3)'+A_4\hookrightarrow \Tcal$ whose restriction to 
$(2A_3)'$ is as above: $K(W)$  the diagonal in $2K_{A_3}$.
\end{enumerate}

In all these cases, the group $\out_{\Orth(\Tcal)}(W)$ is the product of the permutation groups of the summands of $rA_3 +sA_4$ which preserve type and the above  bundling (so its outer action on $W$ is faithful). Its outer action on $W^\perp$  has the same kernel as its action 
on the radical $K_{\Tcal (W)}$ of the $W$-fixed point space $\Tcal^W$ (on which $W\times W^\perp$ acts trivially and which is also the image of $
K_{A_3}^r\to \Tcal^W$).

Furthermore, there is only one  conjugacy class $\Cscr_{\Tcal}(A_{11})$ of reflection subgroups of $\Orth({\Tcal})$ isomorphic with $\Sfrak_{12}$; it is represented by $\Sfrak_E$.
Two geometric subgroups contained in $\Sfrak_E$ that are conjugate in $\Orth(\Tcal)$  are conjugate in $\Sfrak_E$. An $\Sfrak_E$-conjugacy class  consisting of geometric subgroups isomorphic  with $\Sfrak_4^r\times \Sfrak_5^s$  is given by a pair $(r,s)$ for which $4r+5s\le 12$ (and is represented by specifying $r$ 4-element subsets and $s$ 5-element subsets of  $E$, all 
pairwise disjoint).
\end{proposition}
\begin{proof}
We first assume that $s=0$. So $W$ is given by an orthogonal map $rA_3\to \Tcal$. Its image will be  a quotient $\Tcal(W)$ of $rA_3$ whose radical  $K_{\Tcal(W)}$  is  
$K_{A_3}^r/K(W)$.
Since the quadratic form on $A_3$ is zero on its radical, the same is true for $\Tcal(W)$ and  hence $\Tcal(W)\cong rA_2 \operp K_{\Tcal(W)}$.  The classification of
embeddings of $rA_3$ in $\Tcal$ is provided by Lemma \ref{lemma:embeddingcriterion}: it is necessary and sufficient that  $2r+2\dim K_{\Tcal(W)}\le 10$ and if we have equality, $r$ must be odd. This implies $r\le 3$ and that when  $r=3$, $K(W)$ is nonzero. Lemma \ref{lemma:embeddingcriterion} also tells us that in these cases  the embedding is unique up to an orthogonal transformation of $\Tcal$. 
The determination of the $\Orth(\Tcal)$-orbits of orthogonal maps $rA_3\to \Tcal$ is now straightforward and gives the listed cases.

When $s>0$,  the orthogonal complement of an embedding $A^s_4\hookrightarrow \Tcal$ is type $A_{12-4s}$   and hence nonsingular.  The classification of  conjugacy classes of subgroups of reflection type isomorphic with $\Sfrak_4^r\times \Sfrak^s_5$   with $s>0$ boils
down to conjugacy class of subgroups of reflection type isomorphic with $\Sfrak_4^r$ in  $\Orth(A_{12-4s})$. The proof is then similar.

This argument  shows that $\Tcal^W=\Tcal(W)^\perp$ and that the action of $N_{\Orth(\Tcal)}(W)$ on $\Tcal^W$ is through  $W^\perp.\out_{\Orth(\Tcal)}(W)$, with  $W^\perp$ appearing as the subgroup of $W^\perp.\out_{\Orth(\Tcal)}(W)$ acting trivially on $K_{\Tcal(W)}$. This implies the assertion of the second paragraph.

For the last assertion, we note that $A_{11}\cong H_E$ has its  radical spanned by $\chi_E$ and so any orthogonal map  $A_{11}\to \Tcal$ must kill it and factor through an isomorphism $H'_E\cong \Tcal$.  This implies  that all such copies of $\Sfrak_{12}$ are in the same $\Orth(\Tcal)$-conjugacy class. The rest is left to the reader. 
\end{proof}

We put most of this information in Table \ref{table:a10}.
Here  $\Hfrak(A_2)$  stands for the Heisenberg group obtained from $K_{\Tcal(W)}\otimes A_2$, where $K_{\Tcal(W)}=K_1+ K_2$. 
The action of $\Orth(A_2)\cong\Orth^1_2(\FF_2)\cong  \Sfrak_3$ is on the  $A_2$ tensor factor. The outer group $\Sfrak_2$ acts on $\Hfrak(A_2)$ through $K_{\Tcal(W)}$ 
and exchanges $K_1$ and $K_2$.

\begin{small}
\begin{table}
\begin{centering}
\begin{tabular}{|c|c|c|c|}
\hline 
conj.\ cl. of $W$ & $K_{\Tcal(W)}$ &  $W^\perp=$ & $\out_{\Orth(\Tcal)}(W)\twoheadrightarrow$ \\
 in $\Orth(\Tcal)$ && $C_{\Orth(W)}(W)$ & $\out_{\Orth(\Tcal), W}(W^\perp)$\\
  \hline\hline
$A_3$ &  $K_1$  & $\FF_2^6\rtimes \Sfrak_8$ & $\{1\}=\{1\}$\\
\hline
$A_4$  &  $0$   & $\Sfrak_8$ &  $\{1\}=\{1\}$\\
\hline
$2A_3$   &  $K_1\oplus K_2$ &   $\Hfrak(A_2)\rtimes \Sfrak_3$ & $\Sfrak_2=\Sfrak_2$\\
%\hline
$(2A_3)'$ &  $K_1\equiv K_2$ &  $\FF_2^4\rtimes \Sfrak_5$ & $\Sfrak_2\twoheadrightarrow\{1\}$\\
\hline
$3A_3$ & $\frac{K_1\oplus K_2\oplus K_3}{\text{main diag.}}$ & $\wedge^2K_{\Tcal(W)}\cong \FF_2$ & $\Sfrak_3=\Sfrak_3$\\
%\hline
$(2A_3)'+A_3$ &  $K_1(\equiv K_2)\oplus K_3$  & $\wedge^2K_{\Tcal(W)}\cong \FF_2$ & $\Sfrak_2\twoheadrightarrow\{1\}$ \\
%\hline
$(3A_3)'$ &  $K_1\equiv K_2\equiv K_3$    & $\FF_2^2\rtimes\Sfrak_2$ & $\Sfrak_3\twoheadrightarrow\{1\}$\\
\hline
$A_3+A_4$  &  $K_1$ &    $\Sfrak_4$ & $\{1\}=\{1\}$\\
\hline
$(2A_3)'+A_4$  &  $K_1\equiv K_2$   &$\{1\}$ & $\Sfrak_2\twoheadrightarrow\{1\}$\\
\hline
$2A_4$  &  $0$ &    $\Sfrak_3$ & $\Sfrak_2\twoheadrightarrow\{1\}$\\
\hline
\end{tabular}
\caption{The conjugacy classes of geometric subgroups of $\Orth(\Tcal)$. Here  $K_i$ stands for the radical of the $i$th $A_3$-summand and $K_{\Tcal(W)}$ and the second column describes $K_{\Tcal(W)}$  as the quotient of their direct sum. The last column gives $\out_{\Orth(\Tcal)}(W)$ and its image in the outer automorphism group of $W^\perp$.}\label{table:a10}
\end{centering}
\end{table}
\end{small}

\section{The Bimonster and its geometric subgroups}\label{sect:Bimonstercomplex}

The structure which we discuss below  was essentially exhibited more than three decades  
ago by Conway, Norton and others. As we mentioned in the introduction, their work seems to be tailor made for the purposes of this paper. 
The characterization of the Bimonster given by Proposition \ref{prop:bimonsterchar} is based on their work.

\subsection{Geometric subgroups of the  Bimonster}\label{subsect:bimonster}
It is well-known that an automorphism group of the permutation group $\Sfrak_n$ is inner when
$n\not=6$. So for such $n$ and $k\le n$, the conjugacy class of the standard subgroup   
$\Sfrak_k\le \Sfrak_n$ is intrinsically defined; we call this the \emph{reflection conjugacy class} of 
$\Sfrak_k$ in $\Sfrak_n$ or of any group isomorphic to $\Sfrak_n$ ($n\not=6$).
\\

Consider  a group $\Gcal$  which  comes with a distinguished conjugacy class  
$\Cscr(A_4)$ of subgroups isomorphic to $\Sfrak_5$ whose union  generates $\Gcal$ .  
In view of the preceding, then is defined for every $k\le 5$ a reflection conjugacy class 
$\Cscr(A_k)$ of $\Gcal$ of subgroups isomorphic to $\Sfrak_k$. 
Assume the following properties hold:
\begin{enumerate}
\item[(G0)] The graph whose vertex set is $\Cscr(A_3)\sqcup\Cscr(A_4)$ and whose edges 
are defined by the pairs $(W',W)\in \Cscr(A_3)\times \Cscr(A_4)$ with $W'\le W$ is connected.
\item[(G1)] If  $W\in \Cscr(A_4)$, then the $\Gcal$-centralizer of $W$  is isomorphic with $\Sfrak_{12}$ and the members of $\Cscr(A_3)$ resp.\ $\Cscr(A_4)$ contained in that centralizer make up its reflection conjugacy class of subgroups of that type.
\item[(G2)] If $W\in \Cscr(A_3)$, then the $\Gcal$-centralizer of $W$ is isomorphic with $\Orth^1_{10}(\FF_2)$ and the members of 
$\Cscr(A_3)$ resp.\  $\Cscr(A_4)$ contained in 
that centralizer  make up  its (standard) conjugacy class of subgroups of that type.
 \item[(G3)] Let $(W',W)\in \Cscr(A_3)\times \Cscr(A_4)$ be such that $W'\le W$. Then 
the $\Gcal$-centralizer of $W$  is in the reflection conjugacy class of the $\Gcal$-centralizer of $W'$.
\end{enumerate}

\begin{proposition}\label{prop:bimonsterchar}
Any group ${\Gcal}$ satisfying the four properties (G0), (G1), (G2) and (G3) is isomorphic to the Bimonster.
\end{proposition}
\begin{proof}
Conway and Pritchard  \cite{cp} consider groups $\Gcal$ not isomorphic
to $\Sfrak_{17}$ that satisfy a weaker form of (G1) and are minimal for that
(weaker) property (it is their `Axiom'). One such group is the Bimonster
(see \cite{norton}, Lemma 13 and the paragraph following its proof). Conway and
Pritchard then show in Theorem 1 of \cite{cp} and the discussion at the end
of their paper that $\Gcal$ must be a quotient of the group $Y_{555}$
as defined by Norton in \cite{norton} (but beware that $Y_{555}$ is defined in
various different (but related) ways in \cite{atlas}, \cite{cp}, \cite{cns}, and \cite{norton}).
Norton's $Y_{555}$ is the quotient by a single particular relation of a
simply laced Coxeter group whose diagram is $Y$-shaped: a central
vertex with three legs, each containing five more vertices. That this
$Y_{555}$ is indeed the Bimonster was subsequently shown by the combined
work of Ivanov \cite{ivanov:monstergeom} and Norton \cite{norton}. This was later simplified by Ivanov in \cite{ivanov}. (I thank Leonard Soicher for pointing me to these references and explaining to me how they relate to each other.)

Now let $\Gcal$ be any group satisfying the four properties (G0), (G1), (G2), (G3). Let
$(W,W')\in \Cscr(A_4)\times \Cscr(A_3)$ with  $W'\le W$.
Since (G1) holds, from Conway and Pritchard \cite{cp} (Theorem 1 and its
proof and Theorem 7), we see that either
\begin{enumerate}
\item[(i)]  $W$ is contained in a copy
of $\Sfrak_{17}$ in which the centralizer of $W'$ is a copy of $\Sfrak_{13}$, or 
\item[(ii)] $W$ is contained in a subgroup $\Gcal_W$ of $\Gcal$ that is generated
by $W$ and a copy of $\Orth_{10}^1(\FF_2)$ centralizing $W'$ and this $\Gcal_W$ is a
quotient of Norton's $Y_{555}$. 
\end{enumerate}
The first case cannot hold, since by
(G2), the $\Gcal$-centralizer of $W'$ is a copy of $\Orth_{10}^1(\FF_2)$, which
contains no element of order $13$.  In the second case, which must hold,
the Ivanov-Norton result shows that the subgroup $\Gcal_W\le \Gcal$ introduced there is the Bimonster.
So it remains to show that $\Gcal_W =\Gcal$.

Let us first observe that Property (G2) implies  that the copy of $\Orth_{10}^1(\FF_2)$ mentioned in case (ii)
must be the full $\Gcal$-centralizer  $C_{\Gcal}(W')$ of $W'$.

Since the $\Sfrak_4$-subgroups of $W$ are pairwise $W$-conjugate, the group $\Gcal_W$
then contains the $\Gcal$-centralizer (isomorphic to $\Orth_{10}^1(\FF_2))$ of
every $\Sfrak_4$-subgroup of $W$, so $\Gcal_W$ does not depend on the choice $W'$
of $\Sfrak_4$-subgroup of $W$. We further claim that if $W_1\in \Cscr(A_4)$ with
$W'\le W_1$, then $\Gcal_W=\Gcal_{W_1}$.  By assumption there exists
a $g\in \Gcal$ such that $gW_1 g^{-1} = W$. Conjugation by $g$ takes $W'$
to a possibly different $\Sfrak_4$-subgroup of $W$, but since all such subgroups
are conjugate inside $W$, we may (and will) assume that $gW'g^{-1} = W'$.
Then $g\in N_{\Gcal}(W') = W' \times C_{\Gcal}(W')$.  By the observation above, the latter is 
a subgroup of $\Gcal_W$. This implies that $\Gcal_{W_1} = g^{-1}\Gcal_W g
= \Gcal_W$.

The connectivity of the graph defined in (G0) then implies that 
$\Gcal_W$ is independent of the choice of $W\in \Cscr(A_4)$.  Since we assumed
that $\Cscr(A_4)$ generates $\Gcal$, it follows that $\Gcal = \Gcal_W$,
a copy of the Bimonster.
\end{proof}

\begin{definition}\label{def:geometric subgroup}
For  ${\Gcal}$ as above, a subgroup of $W< {\Gcal}$ will be called \emph{geometric} of type $(p,q)$ if it is  isomorphic with $\Sfrak_4^p\times \Sfrak_5^q$ for certain $p$ and $q$ with  each factor in $\Cscr(A_3)$ or $\Cscr(A_4)$.  We denote the collection of all geometric subgroups of $\Gcal$ (partially ordered by the inclusion relation) by $\Cscr$.
\end{definition}

If $W\subset \Gcal$ is a geometric subgroup, then we write  $\out_{\Gcal}(W)$ for the group of  permutations of the factors of $W$ induced by the $\Gcal$-normalizer  $N_{\Gcal}(W)$.
Note that then   
\begin{equation}\label{eqn:outW}
N_{\Gcal}(W)=(W\times C_{\Gcal}(W)).\out_{\Gcal}(W).
\end{equation}

If we choose one of the factors of $W$, then the product of the remaining factors is in the centralizer of that factor and so in a copy of $\Orth^1_{10}(\FF_2)$  or $\Sfrak_{12}$  and is therein a geometric subgroup in the previous sense.

 With the help of  the Table  \ref{table:a10}  we can enumerate their conjugacy classes.

\begin{proposition}\label{prop:geomclassinbimonster}
There are precisely 13  conjugacy classes of geometric subgroups $W\le \Gcal$, as listed in Table \ref{table:ghat}, where (as explained below) the notation  indicates  the incidence relations and how   $\out_{\Gcal}(W)$  permutes its  factors. 

\begin{small}
\begin{table}
\centering
\begin{tabular}{|c|c|c|c|}
\hline 
conj. class of $W$& $K_W^{\Gcal}$ & $C_{\Gcal}(W)$ & $\out_{\Gcal}(W)\twoheadrightarrow $\\
in $\Gcal$ &&&$\out_{\Gcal, W}(C_\Gcal(W))$\\
\hline\hline
 $\Cscr(A_3)$  & $K_1$ & $\Orth^1_{10}(\FF_2)$  & $\{1\}=\{1\}$\\ 
 \hline
$\Cscr(A_4)$  & $0$ & $\Sfrak_{12}$   & $\{1\}=\{1\}$ \\ 
\hline
$\Cscr(2A_3)$  & $K_1\oplus K_2$  & $\FF_2^6\rtimes \Sfrak_8$ &$\Sfrak_2=\Sfrak_2$  \\ 
\hline
$\Cscr(A_3+A_4)$ & $K_1$  & $\Sfrak_8$  & $\{1\}=\{1\}$\\ 
\hline
$\Cscr(2A_4)$  & $0$ & $\Sfrak_7$ &  $\Sfrak_2\twoheadrightarrow\{1\}$\\ 
\hline
$\Cscr(3A_3)$    & $K_1\oplus K_2\oplus K_3$ & $\Hfrak(A_2)\rtimes \Sfrak_3$ & $\Sfrak_3=\Sfrak_3$\\ 
\hline
$\Cscr(3A_3)'$  & $K_1(\equiv K_2\equiv K_3)$ &$\FF_2^4\rtimes \Sfrak_5$ & $\Sfrak_3\twoheadrightarrow \{1\}$\\ 
\hline
$\Cscr(2A_3+A_4)$ &  $K_1\oplus K_2$   &$\Sfrak_4$ & $\Sfrak_2=\Sfrak_2$ \\ 
\hline
$\Cscr(A_3+2A_4)$ & $K_1$ & $\Sfrak_3$ & $\Sfrak_2\twoheadrightarrow \{1\}$\\ 
\hline
$\Cscr(3A_4)$ & $0$ &  $\Sfrak_2$ & $\Sfrak_3\twoheadrightarrow \{1\}$\\ 
\hline
$\Cscr((3A_3)'+A_3)$ &$K_1(\equiv K_2\equiv K_3)\oplus K_4$  & $\FF_2$ & $\Sfrak_3\twoheadrightarrow \{1\}$\\ 
\hline
$\Cscr(4A_3)'$ & $\frac{K_1\oplus K_2\oplus K_3\oplus K_4}{\text{main diag.}}$ & $\FF_2^2\rtimes\Sfrak_2$ & $\Sfrak_4= \Sfrak_4$\\ 
\hline
$\Cscr((3A_3)'+A_4)$  & $K_1(\equiv K_2\equiv K_3)$ & $\{1\}$  & $\Sfrak_3\twoheadrightarrow \{1\}$\\ 
\hline
\end{tabular}
\caption{Conjugacy classes of geometric subgroups  of $\Gcal$. For the meaning of the second column and its relation with the fourth column, which gives the image of $\out_{\Gcal}(W)$ in $\out(C_\Gcal(W))$,  see Lemma \ref{lemma:virtualradical}. }
\label{table:ghat}
\end{table}
\end{small}
\smallskip
\end{proposition}

The notation is chosen such that the number of $A_3$'s resp.\ the number of $A_4$'s give the number of 
$\Sfrak_4$ resp.\  $\Sfrak_5$ factors. A permutation of the factors of $W$ is in $\out_{\Gcal}(W)$  if and only if it preserves the  bundled factors (in parentheses and primed).
Furthermore, the types of the  geometric subgroups contained in  a given one are determined by their labels, by which we mean 
the following.  There is only one conjugacy class of geometric subgroups of type $(4,0)$. It is denoted $\Cscr(4A_3)'$, because  a geometric subgroup of that type gives after removal of a  $\Sfrak_4$-factor always a member of $\Cscr(3A_3)'$.  On the other hand, if we do this for  a geometric subgroup of  type $\Cscr(3A_3)'$, then we get  a member of $\Cscr(2A_3)$, as there is only one of that type.  Hence a member of $\Cscr((3\A_3)'+A_3)$ contains three members of  type $\Cscr(3A_3)$ and one member of type $\Cscr(3A_3)'$. Likewise for $\Cscr((3\A_3)'+A_4)$. 

This also shows that these incidence properties suffice to determine the label  of any  given geometric subgroup.  

Proposition \ref{prop:geomclassinbimonster} follows from the three lemmas below.  We first treat  the case when $q>0$, as this is easier and helps us to do the remaining cases.
 
\begin{lemma}\label{lemma:$q>0$}
The geometric conjugacy classes of type $(p,q)$ with $4p+5q\le 17$ and $q>0$, make up a single $\Gcal$-conjugacy class. If $W\le \Gcal$ is a member of this conjugacy class, then any permutation of its factors is induced by an inner automorphism.
\end{lemma} 
\begin{proof}
If $W\in \Cscr(A_4)$, then  $W^\perp\cong \Sfrak_{12}$.  There is a  geometric $\Sfrak_{12}$-conjugacy class of type $(p,q-1)$ in $\Sfrak_{12}$ if and only if $4p+5(q-1)\le 12$ or equivalently $4p+5q\le 17$ and it is then unique. So there is then  exactly one geometric 
${\Gcal}$-conjugacy class of type $(p,q)$. It also tells us that if $W\le \Gcal$ is a member of this conjugacy class, then every permutation of factors which fixes one $\Sfrak_5$-factor is induced by its $\Gcal$-normalizer. This implies that 
every permutation is obtained when $q\not=2$. However, the $q=2$ implies that  leaves room for either an additional  $\Sfrak_5$-factor  or a single $\Sfrak_4$-factor (for they must be subgroups of $\Sfrak_7$). In the first case, we  
 get a geometric subgroup $\tilde W$ of type  $(0,3)$) and then we see that the two $\Sfrak_5$-copies can be   exchanged by an  element of the  $\Gcal$-normalizer of $\tilde W$. This element is in the $\Gcal$-normalizer of $W$. 
 
 The remaining case is when $(p,q)=(1,2)$. Write  $W=W_0\times W_1$ with $W_0\in \Cscr(A_3)$ and  
 $W_1\in \Cscr(2A_4)$. Then $W_1$ embeds in $W_0^\perp\cong \Orth(\Tcal)$.  By Lemma \ref{lemma:embeddingcriterion} there is only one orbit  of embeddings of the nonsingular quadratic space $2A_4$ in $\Tcal$. In particular, the exchange of the two summands is induced by an element of $\Orth(\Tcal)$ and hence by an element of $N_\Gcal(W)$.
\end{proof}

We denote these conjugacy classes as in Table \ref{table:ghat} (these are the items with $q>0$).
This leaves us to deal with the case when $q=0$.

\begin{lemma}\label{lemma:$q>0$}
A  geometric conjugacy class of type $(p,0)$ exists   if and only if $p\le 4$.  It is unique when $p\le 2$; we denote them  $\Cscr(A_3)$ and  $\Cscr(2A_3)$.

For $p=3$, there are precisely two, denoted   $\Cscr(3A_3)$  and  $\Cscr(3A_3)'$.  They can be distinguished by the fact that the centralizer of a member of the former contains no  member of $\Cscr(A_4)$, whereas the latter does. 
\end{lemma} 
\begin{proof}
Let  $W\in \Cscr(A_3)$. Then  $W^\perp\cong \Orth^1_{10}(\FF_2)$. Proposition \ref{prop:a3classification} 
classifies the $\Orth^1_{10}(\FF_2)$-conjugacy classes of geometric subgroups of type $(p-1, q)$ in $\Orth^1_{10}(\FF_2)$. We thus get all the geometric conjugacy classes in $\Gcal$ of type $(p,q)$ with $p>0$ and a $\Sfrak_4$-copy singled out.
This implies that  type $(p,0)$ with $p>0$ occurs if and only if $p-1\le 3$.
Some cases occur in this manner more than once: this happens when $p-1=2$  and   $p-1=3$.

We now analyze the case $p=3$. The case $(p, q)=(3,1)$ was already covered by Lemma \ref{lemma:$q>0$}: that lemma tells us that this case is represented by the  $(3A_3)'$-conjugacy class in $\Scal_{12}$. So this must appear in Table \ref{table:a10} as the case 
$(2A_3)'+A_4$. In view of the inclusion relation in $\Orth(\Tcal)$,  it follows that there is a unique  geometric  conjugacy  class  in   $\Gcal$ of type $(3,0)$ whose members contain in their centralizer a member of $\Cscr(A_4)$; we denote that conjugacy class  $(3A_3)'$. For a member $W$ of  this conjugacy class all the permutations of factors are induced by its normalizer.
In particular, if we  single out any  $\Sfrak_4$-factor $W_0$ of $W$ and write $W=W_0\times W_1$, then $W_1$ appears
in $W_o^\perp\cong\Orth(\Tcal)$ as the geometric conjugacy class labeled $(2A_3)'$. This also proves that the other 
geometric conjugacy class in $\Orth(\Tcal)$  labeled $2A_3$ cannot arise  in this manner.  We denote the  associated geometric conjugacy class in $\Gcal$ by $\Cscr(3A_3)$.
\end{proof}

\begin{lemma}\label{lemma:}
There are in $\Gcal$  precisely  two  geometric conjugacy classes of type $(4,0)$, namely of type $\Cscr((A_3)'+A_1)$ and $\Cscr(A_4)'$ (a member of which contains 1 resp.  4   geometric subgroups  of type $(3A_3)'$). The $\Gcal$-normalizer of such a geometric subgroup permutes its factors according to the notation and is maximal for that property (giving a copy of $\Sfrak_3$ resp.\  $\Sfrak_4$). 
\end{lemma}
\begin{proof}
A geometric subgroup of $\Gcal$ of type $(4,0)$ comes (up to conjugacy) from geometric subgroup of 
$\Orth(\Tcal)$ of type $(3,0)$ and Table \ref{table:a10} shows that there at most three   such cases denoted there 
$3A_3$ and $(3A_3)'$ and $(2A_3)'+A_3$. To these correspond  at most three geometric conjugacy classes of  type $(4,0)$, which we denote  $\Cscr(4A_3)$,  $\Cscr((3A_3)'+A_3)$ and $\Cscr(4A_3)'$. It suffices to  show that $\Cscr(4A_3)$ does not exist.

Suppose the contrary, and let $W=\prod_{i\in I}W_i$,  where $I$ is a 4-element set with  $W_i\in \Cscr(A_3)$,   such that  for every factor $W_i$, the product of the remaining factors (denoted $W'_i$) is   (in the notation of Table \ref{table:a10}) of type $(3A_3)$, when considered as a subgroup of  $W_i^\perp$.  According to Table \ref{table:a10}, the  
group $N_\Gcal(W_i)\cap N_\Gcal(W)$  induces the full permutation group of the  factors of $W'_i$.  So $N_\Gcal(W)$ induces the full permutation $\Sfrak_I$ group of all  four factors of $W$. 

We also see from Table \ref{table:a10} that $C:=\cap_{i\in I} W_i^\perp$  is of order $2$ and that for every $i\in I$, it 
 admits the natural description as $C\cong \wedge^2 K_{W'_i}$. We number the elements of $I$, so that we can assume that 
 $I=\{1,2,3,4\}$. We have in $W_4^\perp$  the equalities  $e_1+e_2+e_3=0$, but $e_1\wedge e_2$, $e_2\wedge e_3$,
 and $e_3\wedge e_1$ all give the same nonzero generator of $C$. An element such as $e_1\wedge e_2$ lies  in fact in 
$(W_3\times W_4)^\perp=W_3^\perp\cap W_4^\perp$ and  has  the same meaning when considered as an element of $W_3^\perp$. So each $e_i\wedge e_j$ with $i\not=j$ produces the nonzero element of $C$. This however  contradicts  the identity  $e_1\wedge e_4+e_2\wedge e_4+e_3\wedge e_4=0$.
\end{proof}

%\begin{remark}\label{rem:virtualradical}

Let $W\in \Cscr$. We will  need to know how 
$N_\Gcal(W)/W$ acts on $C_{\Gcal}(W)$ and the following lemma addresses this. 
It gives   a $\Gcal$-analogue of the space $K_{\Tcal(W)}$ we assigned to a geometric subgroup $W$ of $\Orth(\Tcal)$. 

\begin{lemma}\label{lemma:virtualradical}
Let the finite set $P$ index the $\Sfrak_4$-factors of $W$. Then there  exists a $\out_\Gcal(W)$-equivariant quotient $K_W^{\Gcal}$ of the $\FF_2$-vector space $\FF_2^P$ (given as in the second column in Table \ref{table:ghat},  where $\FF_2^P$ is written as $K_1\oplus K_2\oplus \cdots$) such that the  outer action of $\out_\Gcal(W)$  on
$C_{\Gcal}(W)$ is via this quotient and  has the property that if we single out a factor:  $W=W_1\times W'$ with $W_1$ of type $\Sfrak_4$ or $\Sfrak_5$ and with the supplementary product $W'$ identified with a geometric (or trivial) subgroup of $\Orth(\Tcal)$, then via this identification $K_{\Tcal(W')}$ embeds in $K_W^{\Gcal}$ in a manner  that is equivariant for the natural homomorphism 
$N_{\Orth(\Tcal)}(W')/W'\to N_\Gcal(W)/W$.
\end{lemma}
\begin{proof}
An element of $N_\Gcal(W)$ which does not permute the $\Sfrak_4$ factors of $W$ (so only its 
$\Sfrak_5$-factors) will act on $C_\Gcal(W)$ as an inner automorphism. In other words, 
if  we write $W=W'\times W''$ with $W'$ collecting the $\Sfrak_4$-factors and $W''$ collecting the $\Sfrak_5$-factors, then  $K_W^{\Gcal}$ will only depend on how $N_\Gcal(W)$ acts on $P$. 
We shall therefore assume that $W$ is a product of  $\Sfrak_4$-factors. 
If $W$ can be obtained by suppressing a $\Sfrak_5$-factor (so that  $W_1\times W \in \Cscr$ for some 
$W_1\in\Cscr(A_4)$), then $W$ embeds in $\Sfrak_E\subset \Orth(\Tcal)$ and no genuine quotient is taken: we put  $K_W^{\Gcal}:= \FF_2^P$.

When $W$ is not of this type, then $P$ has size $3$ or $4$. We can then split off an $A_3$-factor:  
$W=W_o\times W'$ and identify $C_{\Gcal}(W_o)$ with $ \Orth(\Tcal)$. 
So $W'$ is  a geometric subgroup of $ \Orth(\Tcal)$ of type $2A_3$, $(2A_3)'$ or  $(3A_3)$,  
$(3A_3)'$ or $(2A_3)'+A_3$. By a smart choice of $W_o$, we can avoid the last case (it then gets replaced by
the type $(3A_3)'$). Then $W=W_o\times W'$ lies in $\Cscr(3A_3)$, $\Cscr(3A'_3)$, $\Cscr(4A_3)'$ or $\Cscr((3A_3)'+A_3)$ respectively.  

We know that the outer action of $\out_{\Tcal}(W')$  on $W'{}^\perp=C_\Gcal(W)$ is via 
$K_{\Tcal(W')}$. This suggests that we take for 
$K_W^{\Gcal}$ the full $\FF_2^P$ when $W\in \Cscr(3A_3)$,  
 the quotient of $\FF_2^P\to \FF_2$ given by the sum when $W\in \Cscr(3A'_3)$,
the   quotient of $\FF_2^P$ by its main diagonal when $W\in \Cscr(4A'_3)$, 
the direct sum of  the one-dimensional  $K_\Gcal^{W_o}$ and the 
 one-dimensional $K_{\Gcal}^{W'}$ when $W'\in \Cscr (3A_3)'$.
 
By construction it then  has the property stated in the last clause. 
\end{proof}

The following corollary will be used at the very end of the paper.

\begin{corollary}\label{cor:extendedaction}
In the situation of Lemma \ref{lemma:virtualradical}, the action of 
$N_{\Orth(\Tcal)}(W')/W'$ on $\Tcal^{W'}$ extends in a unique manner to $N_{\Gcal}(W)/W$.
\end{corollary}
\begin{proof}
When  $W_1$ is of type $A_4$, the map  $N_{\Orth(\Tcal)}(W')/W'\to N_{\Gcal}(W)/W$ is an isomorphism and so there is nothing to show. We shall therefore assume that $W_1$  of type $A_3$.
In that case, we need to rely on the classification given by Table \ref{table:ghat}. 
Fortunately there are then only a few cases to deal this, for the last column of that Table shows that often  
(namely when the  surjection is to the trivial group), the group $N_{\Gcal}(W)$ acts on $C_{\Gcal}(W)$ by inner automorphisms.  In these cases, $C_{\Gcal}(W)$ has trivial center, so that we get a retraction $N_{\Gcal}(W)\to C_{\Gcal}(W)$ and thereby an action $N_{\Gcal}(W)\to C_{\Gcal}(W)$ on $\Tcal^{W'}$.
  
The remaining  cases  are when $W$ is of type $2A_3$, $2A_3+A_4$, $3A_3$ or $(4A_3)'$. 
Since $2A_3+A_4$ is subsumed by the case $2A_3$, there are in fact three cases to consider:
$W$ of type   $2A_3$, $3A_3$, $(4A_3)'$. Then $K^\Gcal_W=\FF_2^k$ ($k=2,3,4$)  and $\out_{\Gcal}(W)$ is the full permutation group of the factors, so isomorphic with $\Sfrak_k$.
Furthermore, $W'$ is of $\Orth(\Tcal)$-type $(k-1)A_3$, the group  $\out_{\Orth(\Tcal)}(W')$ is the full permutation group of the factors, so isomorphic with $\Sfrak_{k-1}$ and acts as such on the radical $K_{\Tcal (W')}\cong \FF_2^{k-1}$  (which is also the radical of $\Tcal^{W'}$). 

Table \ref{table:ghat} shows that $N_{\Orth(\Tcal}(W')/W'=W'{}^{\perp}. \Sfrak_{k-1}$ is realized as the 
subgroup of $\Orth(\Tcal^{W'})$ which acts on the radical $K_{\Tcal (W')}\cong \FF_2^{k-1}$ as $\Sfrak_{k-1}$. The group $\Orth(\Tcal^{W'})$ acts on this radical as 
its  full linear group, which for $k=1,2$ a copy of $\Sfrak_2$ resp.\ $\Sfrak_3$. So this identifies 
$\Orth(\Tcal^{W'})$ with $W'{}^{\perp}. \Sfrak_{k}\cong N_{\Gcal}(W)/W$ and thus defines an action 
of $N_{\Gcal}(W)/W$ extending the action of  $N_{\Orth(\Tcal}(W')/W'$.

For $k=3$, we have in fact an equality:  $W'{}^{\perp}. \Sfrak_3=\Orth(\Tcal^{W'})$, but in that case we have a surjection $N_{\Gcal}(W)/W\cong W'{}^{\perp}. \Sfrak_4\to W'{}^{\perp}. \Sfrak_3$ (its  kernel is the extension of the Kleinian 4-group $V_4\subset \Sfrak_4$ by $W'{}^{\perp}$) and hence also an action 
of $N_{\Gcal}(W)/W$ extending the action of  $N_{\Orth(\Tcal}(W')/W'$.

It is clear that in all cases this extension is unique.
\end{proof}

Here is a partial   summary of the preceding. If we fix some $W_o\in \Cscr(A_3)$ and 
identify  $C_\Gcal(W_o)$ with $\Orth(\Tcal)$, then  we have a map
$W'\in \Cscr_{\Tcal}\mapsto W_o\times W'\in \Cscr$. If we also include $W_o$ itself (so take for $W'$ the trivial group),   then this produces representatives of all the  geometric conjugacy classes in $\Gcal$ of type $(p,q)$, with $p>0$. This need not be an embedding: $W'$ and $W''$  may be  nonconjugate subgroups of $\Orth(\Tcal)$ while 
$W_o\times W'$  and $W_o\times W''$ are conjugate in $\Gcal$. 

Next fix a  $W_1\in \Cscr(A_4)$ with $W_1\supset W_o$ and identify 
$C_{\Gcal}(W_1)$ with $\Sfrak_{12}$. If we denote by $\Cscr_{W_1^\perp}$ the collection  of geometric subgroups of $W_1^\perp$, then we get in the same manner a map 
$W\in \Cscr_{W_1^\perp}\mapsto W_1\times W\in \Cscr$. If we include $W_1$ itself, then this produces representatives of all the  geometric conjugacy classes in $\Gcal$ of type $(p,q)$, with $q>0$. This is in fact an embedding.

We thus hit  in terms of a  fixed pair $W_o<W_1$ all the  geometric conjugacy classes in $\Gcal$.

\subsection{Simplicial complexes associated with the Bimonster}
The set $\Cscr$ of geometric subgroups of $\Gcal$ is  partially ordered by inclusion and hence  determines a simplicial complex  $\Sigma(\Cscr)$:
a $d$-simplex of $\Sigma(\Cscr)$ is  given by a strictly increasing  chain $W_\pt=(W_0<W_1<\cdots <W_d)$ of such subgroups and will be denoted  
$\la  W_\pt\ra$. Its  $\Gcal$-stabilizer is   $\Gcal_{[\la W_\pt\ra]}=\Gcal_{\la W_\pt\ra}$ and this is just the group $N_{W_\pt}(\Gcal)=\cap_k  N_{W_k}(\Gcal)$. 

No two $0$-simplices which span a $1$-simplex lie in the same $\Gcal$-orbit and hence  $\Gcal\bs \Sigma(\Cscr)$  is still a simplicial complex.  We denote this simplicial complex by $\check\Sigma(\Cscr)$. By Proposition \ref{prop:geomclassinbimonster}  it has  thirteen vertices (hence is finite).

\begin{corollary}\label{cor:scnerve}
The complex $\Sigma(\Cscr)$ is simply connected\end{corollary}
\begin{proof}
Let $\tilde\Sigma \to \Sigma(\Cscr)$ be  a universal covering with Galois group $\Pi$. 
Then  $\tilde\Sigma $ comes endowed with an action of a group $\tilde \Gcal$ that is an extension of 
$\Gcal$ by $\Pi$. Since the $\tilde \Gcal$-stabilizer of a simplex of  $\tilde\Sigma$ maps isomorphically to 
the $\Gcal$-stabilizer its image in  $\Sigma(\Cscr)$, we see that $\tilde \Gcal$ comes with a conjugacy class of subgroups isomorphic with $\Sfrak_5$  which satisfies  the four properties (G0), (G1), (G2), (G3) (and denoted there by $\Cscr(A_4)$).  But these properties characterize the 
Bimonster and so  $\tilde\Gcal=\Gcal$ and $\Sigma(\Cscr)$ is simply connected. 
\end{proof}

Given a  $W\in \Cscr$, then the full subcomplex of $\Sigma(\Cscr)$  
consisting of simplices defined by chains  whose members are $\ge W$ will be  denoted 
$W\bs\Sigma(\Cscr)\subset \Sigma(\Cscr)$: it is the `upper star' of  $\la W\ra$ (we follow Quillen's notation here).
This subcomplex is contractible, because it has a minimal member (namely  $W$).
An intersection  $W\bs\Sigma(\Cscr)\cap  W'\bs\Sigma(\Cscr)$ is empty unless
$W\cup W'$ is contained in a geometric subgroup and is then equal  to  $W''\bs\Sigma(\Cscr)$, where 
$W''$ is the   geometric subgroup generated by  $W\cup W'$.
This implies:

\begin{lemma}\label{lemma:coarsenerve}
The collection   $\{W\bs\Sigma(\Cscr)\}_{W\in \Cscr(A_3)}$ is a Leray  covering of $\Sigma(\Cscr)$  and  so its nerve,  denoted $\Sigma(\Cscr(A_3))$,  has the same homotopy type as  $\Sigma(\Cscr)$. In particular, it is simply connected. \hfill $\square$
 \end{lemma}

We shall encounter $\Sigma(\Cscr)$ in  a specific geometric setting as the  nerve of an open covering of a space $\hat U$. 
This may motivate the following discussion. We begin with the following elementary result.

\begin{lemma}[Generalized van Kampen theorem]\label{lemma:generalizedvK}
Let  $T$ be a path connected space and $\Ufrak$ a covering of $T$ by simply connected subspaces such that every nonempty intersection of two of its members  is path connected. Then the fundamental groups of $T$ and of the nerve of $\Ufrak$
are isomorphic.
\end{lemma}
\begin{proof}
This is surely known, but as we did not find a reference, we give a proof. We discuss the cases when $\Ufrak$ has 2 resp.\ 3 elements first and do the general case. If  $\Ufrak$ has only two elements, then  its nerve is a $1$-simplex and 
the classical  van Kampen theorem implies that $T$ is simply connected.  
Now consider the case when $\Ufrak$  has only three elements,  $\{U_i\}_{i\in \ZZ/3}$.  

If some twofold intersection is empty, say $U_{12}=U_1\cap U_2$ is,  then the nerve is a string, and hence simply connected. 
On the other hand,  $T$ is then covered by the simply connected unions $U^{01}:= U_0\cup U_1$ and  $U^{02}:= U_0\cup U_2$, 
whose intersection  is the path connected $U_0$ and so it follows from the previous case that $T$  is also simply connected. 
We therefore assume that each intersection $U_{i-1,i}$ is nonempty. 

We first do the case when $U_{012}:=U_0\cap U_1\cap U_2$ is empty.   Since $U^{12}$ is simply connected  and  meets $U_0$ 
(which is also simply connected)  in two path components (namely $U_{01}$ and $U_{02}$), it easily follows that the fundamental group of $T$  is infinite cyclic 
and generated  by a loop which visits successively $U_0, U_1, U_2, U_0$. This also represents a loop of the nerve 
(which  is the boundary of a $2$-simplex)  and defines a generator of its  (infinite cyclic) fundamental group. 

Next we do the case when  $U_{012}\not=\emptyset$.  Then the nerve is a $2$-simplex and hence simply connected. Since each $U^{i,i+1}$ is path connected, any  loop based in $U_{012}$, is  homotopic to a composite of loops, each of which stays in some $U^{i,i+1}$. Since  $U^{i,i+1}$ is simply connected, such a  loop is null homotopic. 

The general case follows formally from these cases: if a loop visits successively members $U_0, \dots, U_N=U_0$ of  $\Ufrak$, then $U_{i-1,i}$ must be nonempty for all $i$ and hence define an edge loop of the nerve of $\Ufrak$. A homotopy between edge paths is simplicial: it is a composition of elementary homotopies
(involving $2$-simplices). The preceding shows that a homotopy of edge paths can be lifted to $T$. That conversely, a homotopy of loops in $T$ gives rise  to a  simplicial homotopy of edge paths, is clear.
\end{proof}

Lemmas \ref{lemma:coarsenerve} and \ref{lemma:generalizedvK} give:

\begin{corollary}\label{cor:simply connected cover2}
Let $\hat U$ be a $\Gcal$-space and $\{\hat U\la W\ra\}_{W\in \Cscr(A_3)}$ a $\Gcal$-equivariant  open covering of $\hat U$
whose nerve is $\Sigma(\Cscr(A_3))$. If each $\hat U$ is simply connected and each nonempty intersection 
$U\la W\ra\cap U\la W'\ra$ is path connected, then  $\hat U$ is simply connected.\hfill $\square$
\end{corollary}

\section{Universal Reflection covers}\label{sect:universalcovers}

\subsection{Reflection covers}
In this section we mostly work in the setting of complex-analytic geometry.

Let $X$ be a normal complex-analytic variety and $D\subset X$ a (possibly empty) subvariety of codimension one 
such  that $X\ssm D$ has the structure of an orbifold (let us then call $(X,D)$ an \emph{orbifold pair}).  
The following notion will be  central to this paper.

\begin{definition}\label{def:bifcover} 
A \emph{reflection cover of an orbifold pair $(X,D)$} is a surjective morphism $f: \hat X\to X$ with the property that  
$X$ can be covered by open subsets $U$ such that each connected component of  $f^{-1}U$ is finite over $U$, 
is an orbifold cover  over $U\ssm D$ and has order 2 ramification over $U\cap D$ at worst.

We say that a  reflection cover $f: \hat X\to X$  is \emph{universal} if  $\hat X$ is connected and every reflection cover $\hat X'\to X$ of $(X,D)$ factors through $f$ by a morphism $f': \hat X\to \hat X'$. 

We say that $f$ is \emph{universal at $p\in X$}  (resp.\  is \emph{locally universal}) 
if $p$ admits a basis of  neighborhoods (resp.\ $X$ admits a basis of open subsets) $U$ such that the restriction of $f$ to every connected component  $f^{-1}U$ is universal for the pair $(U, D\cap U)$.
\end{definition}

This notion generalizes  (in an obvious manner, which we shall not spell out) to the \emph{germ} 
of such a pair $(X, D)$ at a closed subvariety of $X$.

So a reflection cover is an ordinary orbifold cover over $X\ssm D$. 

\begin{remark}\label{rem:doublecover}
For many (but certainly not all!) orbifold pairs $(X,D)$ of interest here, there exists a  connected double cover 
$\pi: \tilde X\to X$ of degree 2  whose ramification locus is $D$. This of course  makes  $(\tilde X,\pi^{-1} D)\to (X,D)$ a reflection cover. 
These examples often have in addition the  property that the (induced) orbifold structure on $\tilde X\ssm \pi^{-1} D$ extends to an orbifold structure  on $\tilde X$ in such a manner that  $\tilde X$ has no orbifold points in codimension one. Such an extension will then be unique and we will say that $\pi$  is an \emph{orbifold double cover} of $(X,D)$. In that case,  precomposition of $\pi$ with an 
ordinary  orbifold covering $\tilde\pi: \hat X\to \tilde X$ produces a reflection cover of $(X,D)$. It is not hard to see that if we take for $\tilde\pi$ a universal orbifold covering, 
then $\hat\pi:=\tilde\pi\pi$ is a universal reflection cover.
\end{remark}

The following example,  mentioned to us by Daniel Allcock, shows  that  a universal reflection cover need not  always exist. 

\begin{example}[Allcock]\label{example:allcock}
Take for $X$ the complex plane $\CC^2$ and for $D$ the union of four distinct lines through the origin. 
Then the fundamental  group of $X\ssm D$ has four generators  $g_1,g_2, g_3, g_4$, 
represented by simple loops around the lines, such that their product $c:=g_1g_2g_3g_4$ is central.
 This is also a presentation of that group.  
The quotient $G$ obtained by killing $g_1^2, \dots, g_4^2$ has the property that each of its finite quotients defines a reflection covering of $(X,D)$. One may check that  the profinite completion of $G$ is infinite. This precludes the existence of  a universal  finite quotient of $G$ and hence also of a universal  reflection cover of $(X,D)$. 
(One can also see this by observing that $(X, D)$ has  a simple-elliptic  singularity of degree two as an orbifold double cover; the local fundamental group of such singularity  is  a central  extension of $\ZZ^2$ by $\ZZ$.  Such a group is residually finite and hence has no universal finite quotient.)
\end{example}

The following  construction shows however  that universal reflection covers exist in codimension $1$. We assume here without loss of generality that $X$ is connected. Choose a base point $o\in X\ssm D$ and let $N\subset\pi^\orb_1(X\ssm D,o)$ be the subgroup normally generated by  the square of every simple loop encircling a  smooth point of $D$ (it suffices to choose one for every connected component of the smooth part of $D$) and put
\begin{equation}\label{eqn:reflgroup}
W(X,D):= \pi^\orb_1(X\ssm D,o)/N.
\end{equation}
Let $D_\infty\subset D$
be the set of $s\in D$ which possess a basis of simply connected neighborhoods $U$ for which there exists  a path 
$\g$ in $X\ssm D$  from  some  $o'\in U$ to $o$ such that the composite $\pi^\orb_1(U\ssm D, o')\xrightarrow{\g}  \pi_1^\orb(X\ssm D,o)\to 
W(X,D)$ has infinite image. By means of a topologically locally trivial (e.g.,  Whitney)  stratification of the pair $(X,D)$ one can see that 
$D_\infty$ is a closed complex subvariety of the singular part of  $D$ and hence  of codimension $\ge 2$ in $X$. 
The $W(X,D)$-cover of $X\ssm D$ extends to a universal reflection cover over $X\ssm D_\infty$ and $X\ssm D_\infty$ is maximal for that property (and  
$W(X,D)=W(X\ssm D_\infty,D\ssm D_\infty)$). 

If a universal reflection cover exists for the germ of $(X,D)$ at  $p\in X$, then  clearly $p\notin D_\infty$ and so:

\begin{corollary}\label{cor:local-global}
A universal reflection cover of $(X,D)$ exists if it has one at every point of $X$.
\end{corollary}

\begin{lemma}\label{lemma:bifcover}
If the pair $(X,D)$ admits a  universal reflection cover, then this cover is unique up to a covering transformation and its total space is simply connected as an orbifold. Conversely, a reflection cover $(\hat X, \hat D)\to (X,D)$ with $\hat X$ simply connected as an orbifold which ramifies along $\hat D$ is universal.

A  universal reflection cover has the property that every automorphism of the pair $(X,D)$ lifts to an automorphism 
of the cover. That lift is unique up to covering transformation. 
\end{lemma}
\begin{proof}
Let $f: \hat X\to X$ and $f': \hat X'\to X$ be universal reflection covers. Then each factors through the other: $f: \hat X\to \hat X'\xrightarrow{f'}X$ and $f': \hat X'\to \hat X\xrightarrow{f}X$. 
The resulting composite of covers  $\hat X\to \hat X'\to \hat X$ is then necessarily a covering transformation. 
Hence $\hat X\to \hat X'$ is an isomorphism.

If a reflection cover $(\hat X, \hat D)\to (X,D)$  ramifies along $\hat D$, then every reflection cover dominating 
$(\hat X, \hat D)$ must be unramified in codimension one and hence if $\hat X$ is simply connected as an orbifold,
be an isomorphism.

If $f: \hat X\to X$ is a universal reflection cover, then its precomposition with an (ordinary) universal cover of $\hat X$  is also a reflection cover. It is clearly universal (since $f$ is) and so the universal cover of $\hat X$  must be trivial. In other words, $\hat X$ is simply connected.

The above uniqueness property implies the last clause of the lemma. 
\end{proof}

\begin{remark}\label{rem:}
It is not true that a reflection cover with simply connected total space is universal. Probably the simplest example is $(\CC, \{0\})$, which  has in fact no universal reflection cover at all.
\end{remark}

\subsection{The case of a real reflection group}\label{subsect:Coxeter}
Let $V_o$ be a real finite dimensional vector space,  denote by $V$ its complexification and regard   $V_o$ as a subspace of $V$.
Let $W\subset \GL(V_o)$ be a finite subgroup generated by reflections. Then $W$ is  a  finite Coxeter group which (via the inclusion $\GL(V_o)\subset \GL(V)$)  acts freely on $V\ssm Z(V)$, 
where $Z(V)\subset V$ denotes the union of its (complexified) reflection hyperplanes.  
According to Chevalley, the underlying algebraic variety of the orbifold $V_W$ is   isomorphic to an affine space. The image $Z(V)_W\subset V_W$ 
\nomenclature[A]{$Z$}{refers to a locally finite union of hyperplanes or hyperplane sections}%
of $Z(V)$ is the discriminant of the orbit space formation $V\to V_W$ and 
\[
f:(V,Z(V))\to (V,Z(V))_W:=(V_W,Z(V)_W)
\] 
is a universal reflection cover.  The obvious homomorphism $\det: W\to \mu_2$ gives a factorization through an  orbifold double cover and hence this is of the type described in Remark \ref{rem:doublecover}. 

Every $p\in V_W$ has a neighborhood basis of open subsets $U$ for which the connected components of its preimage in $V$ are simply connected 
and so this   reflection cover $f$ is locally universal. There is also an hereditary  property of  a somewhat different  nature. In order to explain,  
we assume (without loss of generality) that the origin is the only fixed point of $W$. Let   $C\subset V_o$ be a fundamental chamber of $W$. 
Its  closure $ \overline C$ is a strict fundamental domain for the $W$-action on $V_o$ in the sense that every 
$W$-orbit in $V_o$ meets  $\overline C$ in a single point. The resulting  map $r: V_o\to \overline C$ is continuous  piecewise-linear in that it maps 
each face of $V_o$ isomorphically onto one of $\overline C$.

Let $V'_o\subset V_o$ be an intersection of (real) reflection hyperplanes and $V'\subset V$ its complexification. We denote by $W_{[V']}\subset W$  the $W$-stabilizer of $V'$ (here $[V']$ stands for the point of a Grassmannian defined by $V'$---this explains  our notation), by $W_{V'}\subset W_{[V']}$ its pointwise stabilizer and  by $W(V')$ the quotient $W_{[V']}/W_{V'}$,  viewed as a subgroup of $\GL(V'_o)$ or $\GL(V')$, and by $W(V')^{\refl}$  the subgroup of $W(V')$  generated by its reflections.

We will see that the image of $V'$ under the orbifold cover $V\to V_W$ is conveniently described in terms of 
$r|V'_o$. We first note that $r(V'_o)$ is a  closed union of faces of $\overline C$. 

Let $F_1$ and $F_2$ be two distinct faces  of $V_o$ which are contained and \emph{open} in $V'_o$ such that $\overline F_1\cap \overline F_2$ is a \emph{wall}: the linear span $L_o$ of $\overline F_1\cap \overline F_2$ is of codimension one in  $V'_o$.
\\

Case ($a$): $r(F_2)=r(F_1)$. \label{sss:i} 

Then there exists a $w\in W$ which takes $F_1$ to $F_2$. It is clear that $w$ then preserves $V'_o$ and hence defines an element of $W(V')$.
We also note that $w$ must fix $\overline F_1\cap \overline F_2$ (and hence $L_o$) pointwise. So $w$  acts in  $V'_o$ as a reflection and hence lies in $W(V')^{\refl}$.

\begin{remark}\label{rem:elaborate_case_a}
The pointwise $W$-stabilizer of $L_o$ is a Coxeter  subgroup  $W_{L_o}$ of $W$ whose natural representation space is $V_o/L_o$. The images of $F_1$ and $F_2$ in $V_o/L_o$ are then opposite rays (=one-dimensional faces). These rays define  the same vertex of the Coxeter diagram $\cox(W_{L_o})$. The fact that they are opposite means that this vertex is a fixed point of the opposition involution of the diagram $\cox(W_{L_o})$.
\end{remark}

Case ($b$): $r(F_2)\not=r(F_1)$.  

After applying an element of  $W$, we may assume that $F_1\subset\overline C$. Then $V'_o\cap \overline C=\overline F_1$ and hence  $r(F_2)$ is not contained in $V'_o$. 
Let $w\in W$ be such that $w(F_2)=r(F_2)$. Then $w$ must fix $\overline F_1\cap \overline F_2$ (and hence $L_o$) pointwise, but  $w(F_2)$ is not contained in $V'_o$.  

This can be understood geometrically as that the natural map $V'_{W(V')^{\refl}}\to  V_W$ self-intersects along the image of $L_o$ (see also Corollary \ref{cor:inherit} below).
Here are two examples  illustrating each case.

\begin{example}\label{example:Coxeter}
A  Coxeter group $W$ of type $E_8$ contains a Coxeter group $W'$ of type 
$E_6$. The centralizer of $W'$ in $W$ is  a product of a Coxeter group of type $A_2$ (acting as such in the plane 
$V^{W'}$) and  its center $\{\pm 1\}$ and  the normalizer of $W'$  in $W$ is $W'$ times this centralizer. 
Either acts in $V^{W'}$ as a Coxeter group of type $G_2$.  The plane $V^{W'}$ has 12 open faces 
(the $G_2$-chambers), all defining the same face of $\overline C$. We are in case (a).
\end{example}

\begin{example}\label{example:Coxeter-b}
We assume $W$ is of type $E_6$ and $W'$ of type $A_2\times A_2$: the centralizer of $W'$ in $W$ is a Coxeter group of type $A_2$. Its normalizer in $W$ contains  an element which induces in  the plane $V^{W'}$ minus the identity  and  exchanges the two factors of $A_2\times A_2$. It acts in $V^{W'}$ as the Coxeter group $G_2$ (as above), the  12 $G_2$-chambers  are $W$-faces, but lie in two different $W$-orbits: two adjacent  $G_2$-chambers in the same $A_2$-chamber  are represented by  two distinct 2-dimensional faces of $\overline C$. We are then in case (b).
\end{example}

The following is then clear. 

\begin{corollary}\label{cor:inherit}
Suppose we are the situation above: $V'\subset V$ is an intersection of  reflection hyperplanes of $W$, $W_{[V']}$ is the stabilizer of $V'$, $W_{V'}$ the pointwise stabilizer of $V'$, $W(V'):=W_{[V']}/W_{V'}$ and  $W(V')^{\refl}\subset W(V')$ the  normal (Coxeter) subgroup generated by its reflections. Then in the factorization
\[
V'\to V'_{W(V')^{\refl}}\to V'_{W(V')}\to V_W, 
\]
the first map is the formation of an orbit space for the action of a finite Coxeter group (and so  $V'_{W(V')^{\refl}}$ is isomorphic to affine space),  the second a universal orbifold covering with covering group $W(V')/W(V')^{\refl}$ and the third is a normalization of its image. 
In particular,  if we let $Z(V')\subset V'$ stand  the union of  the  reflection hyperplanes of  $W(V')$,  then 
\[
(V',Z(V'))\to (V',Z(V'))_{W(V')}
\] 
is a universal reflection cover.
\hfill $\square$
\end{corollary}
 
\begin{remark}\label{rem:discrchar}
Since  the discriminant of $V'_{W(V')^{\refl}}\to V'_{W(V')}$ is everywhere of  codimension $>1$, we can characterize 
$Z(V')_{W(V')}\subset V'_{W(V')}$ as the closure of the part of the discriminant  $V'\to V'_{W(V')}$ where it is of codimension one.  \end{remark}

All the cases of interest to us are locally of this form  and so we make the following definition.

\begin{definition}\label{def:coxtype}
We say that an orbifold pair  $(X,D)$ is \emph{locally of Coxeter type} if each point of $X$ has an orbifold chart for which the preimage of $D$ is isomorphic to the discriminant of a finite Coxeter group.
\end{definition}

Recall that a stratification (of a complex-analytic variety $X$, say) is a partition of $X$ into \emph{connected} smooth subvarieties, called \emph{strata}, such that the closure of every stratum is a subvariety and a union of strata. There is an  orbifold version of this notion for which the strata are orbifolds.  

The discriminant of a finite Coxeter group determines the type of that Coxeter group. So if  $(X,D)$ is locally of Coxeter type, then $X$ is naturally stratified, a stratum being a connected component of the locus where the discriminant is locally of a fixed Coxeter type.

Corollaries \ref{cor:local-global} and \ref{cor:inherit}  give:

\begin{corollary}[Existence of a universal reflection cover]\label{cor:inheritance}
If $(X,D)$ is locally of Coxeter type, then $(X,D)$ admits a universal reflection cover. This pair also comes with a natural stratification with the following property: if $Y^\circ$ is stratum and $Y$ is the normalization of its closure, then $Y$ has the structure of an orbifold  and if  $D_Y\subset  Y$ 
stands for  the locus that in the above model  corresponds to $Z(V')_{W(V')}$, then the pair $(Y,D_Y)$ is again locally of Coxeter type. 
\hfill $\square$
\end{corollary}

\subsection{Reflection covers of hyperelliptic moduli}

We begin with discussing a basic example.

\begin{example}[Stable $\PGL_2(\CC)$-orbit of divisors on $\PP^1$]\label{example:DM}
Let  $H^0_d$  
be the  space of effective divisors on $\PP^1$ of degree $d$ and let $D^0_d\subset H^0_d$ be the set of nonreduced divisors. A reflection cover of $(H^0_d,D^0_d)$ is given 
by numbering the points of such a divisor, i.e., by the obvious map $(\PP^1)^d\to H^0_d$. Since $(\PP^1)^d$ is simply connected, this covering is universal. 

This remains true if  we consider such divisors modulo projective equivalence, but restrict to stable ones (where we then take $d\ge 3$). This is the open subset $H^{0,\st}_d$
\nomenclature[A]{$H^{0,\st}_d$}{The space of stable  effective divisors of degree $d$ on  $\PP^1$\nomrefpage}%
  of positive degree $d$ divisors on $\PP^1$ for which each point has multiplicity $<d/2$. The group $\PGL_2(\CC)$ acts properly on this set and yields an orbifold $\Hcal^{0,\st}_d$ with  the  nonreduced divisors defining a hypersurface $\Dcal^{0,\st}_d\subset \Hcal^{0,\st}_d$.
The preimage of $H^{0,\st}_d$ in $(\PP^1)^d$ is an open subset $U$ of  $(\PP^1)^d$ on which $\PGL_2(\CC)$ acts freely and the $\PGL_2(\CC)$-orbit space of $U$ then  defines a 
reflection cover 
\[
(\hat \Hcal^{0,\st}_d,\hat \Dcal^{0,\st}_d)\to (\Hcal^{0,\st}_d,\Dcal^{0,\st}_d).
\]
This cover is universal: $U$ is simply connected because the complement of $U$ in $(\PP^1)^d$ is of complex codimension $>1$ and since $\PGL_2(\CC)$ is connected, its  $\PGL_2(\CC)$-orbit space $\Hcal^{0,\st}_d$ is then
also  simply connected. 

Note that the sign homomorphism $\Sfrak_d\to \mu_2$ defines an intermediate orbifold double cover.
\end{example}

This example is very much related to the moduli space of hyperelliptic curves of genus 
$g\ge 2$.  A smooth complete hyperelliptic curve curve of genus $g$ comes with an involution $\iota$ whose orbit space is of genus zero and whose fixed point set (its set of 
Weierstra\ss\ points) has size $2g+2$. This identifies the coarse moduli space $\Hyp_g^\st$ of smooth hyperelliptic curves with $\Hcal^{0}_{2g+2}$. So the projective  compactification $H^{0,\st}_{2g+2}\supset H^{0,\st}_{2g+2}$ defines a projective  compactification $\Hyp_g^\st\supset \Hyp_g$.

A pair $z,z'$ of   Weierstra\ss\ points of a smooth hyperelliptic curve $C$ defines an element  of $H_1(C; \FF_2)$: just choose an arc $\g$ in $C$ from $z'$ to $[z]$  and take the class $[z]-[z']$ of $\g-\iota_*\g$ in $H_1(C; \FF_2)$ 
(in these $\FF_2$-vector spaces, the  signs can of course be ignored). This is indeed independent of $\g$. A numbering of the Weierstra\ss\ points: $z_0, \dots , z_{2g+1}$ yields a basis $\{[z_{i}]-[z_{i-1}]\}_{i=1}^{2g}$ of $H_1(C; \FF_2)$ and thus puts  a full level 2 structure on $C$. Example \ref{example:DM} implies:

\begin{proposition}\label{prop:hypkernel}
Let  $\widehat \Hyp_g\to \Hyp_g$ be  the 
$\Scal_{2g+2}$-covering defined by a numbering of its Weierstra\ss\ points. Then the normalization  of $\Hyp^\st_g$ in this covering yields a reflection covering $\widehat \Hyp^\st_g\to \Hyp^\st_g$ which is universal.
\end{proposition}

\subsection{Monodromy covers as reflection covers} We first do this for  versal deformations of  hypersurface singularities.
\label{subsect:hypersurfacesing}
Let $(X,x_o)$  be a germ of an $n$-dimensional hypersurface that has $x_o$ as an isolated singular point. 
Such a singularity has a semi-universal deformation $F: (\Xcal,x_o)\to (S,o)$ where $(\Xcal,x_o)$ and $(S,o)$ are smooth 
complex-analytic germs, $F$ is flat and the central fiber (over $o$) has been identified with $(X,x_o)$. 
Let us choose a `good representative' for $F$ in the sense of \cite{looij:icis} and still denote it $F:\Xcal\to S$. 
The critical set of $F$ is a connected complex submanifold  $\Ccal\subset \Xcal$ of codimension  $n+1$.
The restriction $F|\Ccal$ is finite and so $\Dcal:=F(\Ccal)$ is an irreducible  hypersurface in $S$, the \emph{discriminant} of $F$.  The points $x\in \Xcal$ with the property that the $F$-fiber through $x$ has at $x$ a \emph{nonsimple} singularity (i.e., singular, but not of type $A$, $D$ or $E$) define a closed  subvariety $\Dcal_\infty\subset\Dcal$ of codimension $\ge 7$ in $S$ (even $\ge 8$ if $n=1$). Singularity theory tells us that $(S\ssm \Dcal_\infty, \Dcal\ssm \Dcal_\infty)$ is locally like the discriminant of a Coxeter group whose irreducible factors are of type $A$, $D$ or $E$ and that $\Dcal_\infty$ can indeed be characterized as the `bad locus' $D_\infty$ in the above sense. So $(S\ssm \Dcal_\infty, \Dcal\ssm \Dcal_\infty)$ admits a universal reflection cover that cannot be extended to a larger open subset of $S$. 

The  isomorphism type of
$(S, \Dcal)$ does not change if we replace $(X,x_o)$ by its suspension (if the singularity is defined by $f\in\CC\{z_0, \dots , z_n\}$, then its \emph{suspension} is defined by $f+z_{n+1}^2$; somewhat more intrinsically, if we regard  $(X,x_o)$ as a  hypersurface germ in a complex manifold germ $(M,x_o)$, then this is the double cover of  $(M,x_o)$ ramified along $(X,x_o)$). So when it comes to reflection covers, there is no loss of generality in assuming that $n$ is even. In that case the intersection pairing on the reduced cohomology in degree $n$ of a smooth (Milnor) fiber of $F$ over $S\ssm \Dcal$ is symmetric and the monodromy representation 
on this is generated by Picard-Lefschetz reflections  defined by vanishing cycles that are orthogonal with respect to that pairing (more about this below).
An associated monodromy cover of $S\ssm \Dcal_\infty$ extends across 
$S\ssm \Dcal_\infty$ as a reflection cover. This   reflection cover is locally universal: the total space of the cover is smooth and $\Gg(S,\Dcal)$ acts as a reflection group whose point stabilizers are Coxeter groups  of the type mentioned. As was observed in \cite{looij:icm}, certain period maps that can be defined in this setting extend across  this  reflection cover. We do not know whether this cover is simply connected, so in order to get the universal one, we must pass to the universal cover of this ramified monodromy cover.
\\

A variant of this example, which also a generalizes Example \ref{example:DM},  is the moduli space of projective hypersurfaces in $\PP^{n+1}$ of fixed degree $d\ge 3$. That is, we  take for $S$ the projective space 
\begin{equation}\label{eqn:}
H^n_d=\PP(\CC[Z_0, \dots, Z_{n+1}]_d)
\end{equation}
of effective degree $d$-divisors on $\PP^{n+1}$ (this is a complete linear system) and for $D$ the discriminant divisor $D^n_d$ in $H^n_d$ which parametrizes the singular divisors. Then $D_\infty$ is the locus $D^n_{d,\infty}$ which parametrizes the hypersurfaces  that
have a singularity not of type $A$, $D$ or $E$.  Since $d\ge 3$, $D_\infty$  has codimension $\ge 7$ in $H^n_d$.  
The monodromy representation gives rise to a  reflection cover of the pair 
\[
(S\ssm D_\infty,D\ssm D_\infty):=(H^n_d\ssm D^{n}_{d, \infty}, D^n_d\ssm D^n_{d,\infty})
\]
as follows.  
Recall that if $X\subset \PP^{n+1}$ is a smooth hypersurface (so with $s\in S\ssm D$), then its primitive homology is 
(by definition) the kernel of $\Hl_\pt(X)\to \Hl_\pt (\PP^{n+1})$. It is concentrated in the middle dimension $n$  
and if we denote it $\Hl^\circ_n(X)\subset \Hl_n(X)$, then the intersection pairing 
$\Hl^\circ_n(X)\times \Hl^\circ_n(X)\to \ZZ$ is nondegenerate. For $n$ even, $\Hl^\circ_n(X)$ is of corank 1 in 
$\Hl_n(X)$ and the pairing is 
 symmetric, whereas  for $n$ odd,  $\Hl^\circ_n(X)=\Hl_n(X)$ and the pairing is  skew and unimodular. 
 This primitive homology comes with a distinguished subset $\Delta_X\subset \Hl^\circ_n(X)$
  whose elements are 
 called  \emph{vanishing cycles}: these are $\delta\in \Hl^\circ_n(X)$ with the property that 
 $\delta\cdot \delta=(-1)^{n/2}2$ in case $n$ even (for $n$ odd, we have $\delta\cdot \delta=0$, of course). 
 The (\emph{Picard-Lefschetz})  transformation
\[
T_\delta : a\in \Hl_n(X)\mapsto a+(-1)^{n(n-1)/2}(\delta\cdot a)\delta\in \Hl_n(X)
\]
describes the monodromy of a simple loop around the discriminant. For $n$ odd, this is a symplectic transvection 
 and  for $n$ even  an orthogonal reflection. 
The set $\Delta_X$ generates $\Hl^\circ_n(X)$ and  the associated Picard-Lefschetz  transformations generate the 
monodromy group in $\GL(\Hl^\circ_n(X))$ of the universal family over $S\ssm D$. Moreover,  $\Delta_X$  is a single orbit of this group. 

When is $n$ is even, this monodromy group  is of finite index in the orthogonal group of $\Hl^\circ_n(X)$ (see Beauville \cite{beauville}). It produces a Galois cover of $S\ssm D$, which, since each $T_\delta $ is of order two, extends to a `monodromy reflection cover' 
\[
(\hat S\ssm D_\infty, \hat D\ssm D_\infty)\to (S\ssm D_\infty,D\ssm D_\infty).
\] 
It is locally universal when for  every hypersurface with at worst  $ADE$-singularities it restricts to a versal  deformation of its singular multi-germ. This is often the case. (A weaker property suffices, however: all we need is that the map to a versal deformation is surjective on the fundamental groups of the discriminant complements.) 

In case  $H^n_d\ssm D^n_{d, \infty}$ contains the $\SL_{n+2}(\CC)$-stable locus   $H^{n,\st}_d\subset H^n_d$, 
then it is more natural to proceed as in  Example \ref{example:DM} and 
take for $S$ the orbifold $\Hcal^{n,\st}_d:=\PGL_{n+2}(\CC)\bs H^{n,\st}_d$
and for $D$ the image $\Dcal^{n,\st}_d$ of the $\PGL_{n+2}(\CC)$-stable part of $D^n_d$. The reflection covers defined above are then still defined. 
  
\begin{example}[The moduli space of cubic surfaces]\label{example:}
We take for base $S$ the  GIT quotient $\Hcal^{2,\st}_3=\PGL_{4}(\CC)\bs H^{2,\st}_d $ of stable cubic surfaces. This parametrizes all the cubic surfaces with $A_1$-singularities at worst. 
So here $D=\Dcal^{2,\st}_3 \subset S$ is the locus which parametrizes all the singular semistable cubic  surfaces (having $A_1$-singularities only).  It is a classical fact that  the monodromy group of this family  is a Weyl group of type $E_6$. It is perhaps slightly less classical (but  known)   
that the corresponding reflection  cover of $(S,D)$ is simply connected and hence universal. 
(The proof uses the fact that blowing up 6 numbered points in $\PP^2$ determines a smooth marked cubic surface up to projective equivalence. This leads to an  identification of the cover in question with the $\PSL_3(\CC)$-orbit space of the complement  of a closed subvariety in $(\PP^2)^6$ of codimension $\ge 2$ on which $\PSL_3(\CC)$ acts properly, from which one deduces that it is simply connected, but as this is  not needed in what follows, we will not spell this out.)
\end{example}

When $n$ is odd, we look at the monodromy action on $\Hl_n(X; \FF_2)$ instead, for then  each $T_\delta $ acts in
$\Hl_n(X; \FF_2)$ with order two. In that case we let our cover be given by the  monodromy representation of  $\Hl_n(X; \FF_2)$. According to Beauville \cite{beauville} (see also Randal-Williams \cite{rw}) the image of this representation is for $d$ even  the full symplectic group of $\Hl_n(X; \FF_2)$ and for $d$ odd the full orthogonal group of the quadratic form $q:\Hl_n(X; \FF_2)\to \FF_2$  that takes the value one on the vanishing cycles; it has the intersection pairing as its associated bilinear form: $q(a+b)=q(a)+ q(b)+ (a\cdot b)$.

We do not know whether  for the reflection covers just constructed, the total space $\hat S$ is simply connected 
(so that the covering  in question is universal). Let us observe that  the mod 2 monodromy cover  need not be \emph{locally} universal as a reflection cover, the reason being that for example a Weyl group of type $E_8$ contains $-1$ and hence is not faithfully represented on the mod 2 reduction of its root lattice. So if  $s\in S$ is such that $X_s$ has a singular point of that type, then  the family cannot be locally universal at $s$. One of our main results is that when $n=d=3$ (the case of cubic threefolds), the resulting reflection cover is indeed universal.

\section{Eisenstein lattices and ball arrangements}\label{sect:eisenstein}

\subsection{Eisenstein lattices}\label{subsect:eislattices}
 A general reference for this subsection is  Allcock \cite{allcock:leech}, \cite{allcock:y555}. Our notation  differs slightly though and our 
sign conventions are those of \cite{allcock:y555}, rather than  of \cite{allcock:leech}.

Let $L$ be an \emph{even} $\ZZ$-lattice endowed with an orthogonal transformation
$\rho\in \Orth(L)$ satisfying $\rho^2+\rho +1=0$ (in other words, $\rho$ has order $3$ and $1$ is not an eigenvalue of $\rho$). Then 
$L$ becomes a module (denoted  $\LL$) over the Eisenstein ring 
$\Ecal:=\ZZ[\omega]/(\omega^2+\omega +1)$. A hermitian form $\la \; ,\; \ra: \LL \times \LL\to \Ecal$ is then characterized by the property that its real part equals 
$\frac{3}{2}$ times the given bilinear form on $L$:
\begin{multline*}\label{eqn:hermform}
\la x,y\ra =(x\cdot y)+(x\cdot \rho y)\omega  +(x\cdot  \rho^2y)\big) \omega ^2=\\
=(x\cdot y)+(x\cdot \rho y)\omega+ \big(-(x\cdot y)-(x\cdot  \rho y)\big)\bar\omega
=\\=(x\cdot y)(1-\bar\omega)+ (x\cdot \rho y)(\omega-\bar\omega) =\big((x\cdot \rho y)-\omega (x\cdot y)\big)\theta,
\end{multline*}
where $\theta:=\omega -\bar\omega$ (a square root of $-3$).  We thus end up with an \emph{Eisenstein lattice}: a free $\Ecal$-module of finite rank endowed with an $\Ecal$-valued hermitian form. 
Here the hermitian form even takes its values  in $\Ecal \theta$. The hermitian form is also characterized by the property that $\la x,y\ra +\la y,x\ra =3(x\cdot y)$.

We embed  $\Ecal$ in $\CC$ by sending $\omega$ to the primitive third root of unity with positive imaginary part. Then $V(\LL):=\CC\otimes_\Ecal \LL$ is a complex vector space which contains $\LL$ as a lattice. We extend  $\la \; ,\; \ra$ to a hermitian form $V(\LL)\times V(\LL)\to \CC$ and continue to denote this extension by $\la \; ,\; \ra$.  Its signature is half that  of \edit{$L$} (which is necessarily even).  
 
 Much of the basic theory for $\ZZ$-lattices has an extension to Eisenstein lattices. 
Let $\LL^\vee$
denote $\Ecal$-module of \emph{antilinear} maps $\LL\to \theta\Ecal$ ($\Ecal$ acts by post-composition). An $\Ecal$-linear map 
$\LL\to  \LL^\vee$ is defined by assigning to $v\in\LL$ 
the antilinear form $z\in \LL\mapsto \la v, z\ra$. Note that the real part of the underlying map on abelian groups gives, when multiplied with $\frac{2}{3}$, the evident map 
$L\to L^\vee=\Hom(L, \ZZ)$ (this explains our notation).  The 
\emph{discriminant module} $\disc(\LL)$ of $\LL$ is  defined as the cokernel of the corresponding  $\Ecal$-homomorphism: 
\[
\disc(\LL):=\coker(\LL\hookrightarrow \LL^\vee).
\]
We say that $\LL$ is \emph{unimodular} 
 if this module is trivial (which is equivalent to 
$L$ being unimodular). The hermitian form on $\LL$ extends to one on $\LL^\vee$ if we allow it to take values in the Eisenstein field $\QQ(\omega)$. It then  induces a hermitian pairing on the $\Ecal$-module  $\disc(\LL)$, 
\[
\la\, ,\, \ra_\disc : \disc(\LL)\times \disc(\LL)\to \QQ(\omega)/ \theta\Ecal, 
\]
called the \emph{discriminant form}. When  $\LL$ (or equivalently, $L$) is nondegenerate, this form  enables us to describe the overlattices of $\LL$, i.e.,  finitely generated $\Ecal$-modules $\LL'$ sandwiched as  
$\LL\subset \LL'\subset \edit{\QQ\otimes \LL}$  such that
$\la\, ,\, \ra$ still takes on $\LL'\times\LL'$ its values in $\theta\Ecal$. This is the case if and only $\LL'/\LL$  defines a  submodule 
$I\subset  \disc(\LL)$ that is totally isotropic for $\la\, ,\, \ra_\disc$. Then  
the discriminant form establishes an antiduality between $I$ and $\disc(\LL)/I^\perp$ and 
$I^\perp/I$ inherits a hermitian pairing from 
$\disc(\LL)$ and  can be identified as such with  $\disc(\LL')$.  So then  
\[
|\disc(\LL)|=|\LL'/\LL|^2 .\,  |\disc(\LL') |.
\]

Any $a\in L$ with $a\cdot a >0$  defines an orthogonal  reflection in $\QQ\otimes_\ZZ L$ by the formula 
\[
\sigma_a(x):= x-2\frac{(x\cdot a)}{(a\cdot a)}a.
\]
This reflection preserves $L$ if and only if $2(x\cdot a)/(a\cdot a)\in \ZZ$ for all $x\in L$. The latter condition is obviously fulfilled when 
$a\cdot a=2$, but there are also situations of interest, where $a\cdot a=4$ and  $x\cdot a$ even for all $x\in L$.

The condition $a\cdot a >0$ is of course equivalent to  $\la a, a\ra >0$. When $a$ is regarded as an element of $\LL$, then it  defines unitary pseudo-reflections in $V(\LL)$, among them
\[
s_a(x):=x-2\frac{\la x,a\ra }{\la a , a\ra}a, \quad  t_a(x):=x-\bar\omega \theta \frac{\la x,a\ra }{\la a , a\ra}a, \quad t'_a(x):=x+\omega\frac{\la x,a\ra }{\la a , a\ra}a.
\]
which are of order $2$, $3$ and $6$ respectively  and send $a$ to resp.\  $-a$, $\omega a$, $-\bar\omega a$. So $t_a=t'_a{}^2$ and $s_a=t'_a{}^3$. The pseudo reflection $s_a$, $t_a$, $t'_a$   preserves $\LL$ if and only if resp.\ $x\in \LL\mapsto 2\la x,a\ra$, $\theta \la x,a\ra $, $\la x,a\ra $ takes its values in $\la a , a\ra\Ecal$.  For example, if  $a\in L$ is such that $\la a,a\ra=6$, then the reflection $s_a$ is defined if and only if $a\in \theta\LL^\vee$. This is equivalent to $\theta^{-1}a$ defining an element of order $3$ in the discriminant module $\LL^\vee/\LL$.

If  $a\in L$ is such that $a\cdot a=2$, then $\la a, a\ra=3$ and $t_a$ preserves $L$. When considered as an element of $\LL$, we  shall call such an $a$ a \emph{$3$-vector} and  $t_a\in \U(\LL)$ a \emph{triflection}. 
Note that the six $3$-vectors $\pm a$, $\pm \omega a$ and $ \pm \bar\omega a$ all define the same triflection. If in addition,
$\rho a -a\in 3L^\vee$, then the \emph{hexaflection} $t'_a$ is defined as an element of $\U(\LL)$ .

\begin{conventionsnot}\label{conv:groupnames}
The group $\U(\LL)$ acts on the discriminant module $\LL^\vee/\LL$. We denote the kernel of this action by $\G(\LL)$.  It is clear from the definition that any  triflection in $\LL$ is contained in $\G (\LL)$ and so  the  
subgroup generated by these is a subgroup of $\G(\LL)$ that  is normal in $\U(\LL)$. 
We shall denote this group by $\Gg(\LL)$. So $\Gg(\LL)\subset \G(\LL)\subset \U(\LL)$.
 \end{conventionsnot}

\begin{remark}\label{rem:discriminant_duality}
Suppose $\LL$ is unimodular and $\MM\subset \LL$ is a nondegenerate sublattice. 
Then it is easily checked  that 
the   pointwise stabilizer  of $\MM$ in  $\U(\LL)$ restricts isomorphically to $\G(\MM^\perp)$.

We can say more if $\MM$ is a primitive sublattice.
Then every antilinear form $\MM\to \theta\Ecal$ extends to an antilinear form $\LL\to \theta\Ecal$ 
(choose an $\Ecal$-basis of $\MM$ and extend this basis to $\LL$). 
Since $\LL$ is unimodular, such a form can be realized as  $x\in \LL\mapsto \la a, x\ra$ for some  
$a\in\LL$). Two such $a$ define the same antilinear form on $\MM$ if and only if their difference is in $\MM^\perp$. 
So this gives an identification 
\[
\LL/(\MM +\MM^\perp)\cong  \MM^\vee/\MM
\]
as $\Ecal$-modules. 
This produces an isomorphism $\disc(\MM)\cong \disc(\MM^\perp)$ of the discriminant modules   which changes the sign of their  discriminant forms. This works also in the opposite direction: if $\MM$ and $\NN$ are nondegenerate 
Eisenstein lattices, then any  isomorphism of discriminant modules of $\disc(\MM)\cong \disc(\NN)$ which changes the sign of their  discriminant form determines a unimodular overlattice  $\LL\supset \MM\operp \NN$ (the graph of such an isomorphism is maximally isotropic in $\disc(\MM)\operp \disc(\MM^\perp)$). We also see that the   pointwise stabilizer  of $\MM^\perp$ in  $\U(\LL)$ 
restricts isomorphically to $\G(\MM)$. We have a similar characterization of $\G(\MM^\perp)$.

These properties have  well-known  $\ZZ$-counterparts, but are strictly speaking, not a formal consequence  of them, 
as this must take into account $\Ecal$-module structures. 
\end{remark}

\begin{examples}\label{example:basiccases}
Denote by $\LL_k$ the hermitian $\Ecal$-module with basis 
$r_1,r_2,\dots , r_k$  such that  $\la r_i\, r_i\ra =3$,  $\la r_i,  r_{i+1}\ra =\theta$ and $\la r_i,  r_{i+1}\ra=0$ when 
$|i-j|>1$.  (Our sign conventions are that of \cite{allcock:y555}; Allcock denotes in \cite{allcock:leech} this lattice endowed  with the opposite hermitian form by $\Lambda_ \Ecal^k$.)

For $1\le i\le 4$, these lattices are positive definite: the underlying $\ZZ$-lattices are root lattices 
of type $A_2,D_4,E_6,E_8$ (this is why they are denoted in \cite{allcock:y555} by $A_2^\Ecal, D_4^\Ecal,E_6^\Ecal ,E_8^\Ecal$) and 
hence the group $\U(\LL_i)$ is  finite. The discriminant module $\LL_i^\vee/\LL_i$ may be identified with  weight lattice modulo root lattice of the associated root system, and so this can be identified with a vector space of dimension one over  the finite  fields  $\FF_3$ (for $i=1$), $\FF_4$ (for $i=2$), $\FF_3$ (for $i=3$), whereas for $i=4$ this is trivial. 

The center $C_{\Gg(\LL_i)}$ of $\Gg(\LL_i)$ is a subgroup of $\mu_6$ (the group of units of $\Ecal$ acting on $\LL_i$ as scalars) and according to Table 1  of \cite{allcock:y555}, of order $3,2,3,6$ respectively. Let $\mu(\LL_i)\subset \mu_6$  be  the subgroup such that $\mu_6=\mu(\LL_i)\times C_{\Gg(\LL_i)}$, so $\mu(\LL_i)$ is  of order $2$, $3$, $2$, $1$ respectively. That same table tells us  that $\U(\LL_i)=\mu(\LL_i)\times \Gg(\LL_i)$ and that  the group
$\mu(\LL_i)\cong U(\LL_i)/\Gg(\LL_i)$ acts on the discriminant module $\LL_i^\vee/\LL_i$ though multiplication by scalars, thus  identifying it 
with the group of units of resp.\  $\FF_3$, $\FF_4$, $\FF_3$, the trivial module.

This also implies that for any embedding   $\LL_i\hookrightarrow \LL_4$, the restriction of the  $\Gg(\LL_4)$-stabilizer of $\LL_i$  to $\LL_i$ is the full  unitary group $\U(\LL_i)$. 
 \end{examples}
 
 \begin{lemma}\label{lemma:discrLi}
The discriminant forms of  $\LL_1$, $\LL_2$, $\LL_3$ are isomorphic with 
$x\in\FF_3\mapsto x^2\in \FF_3$, $x\in \FF_4\mapsto xx'\in \FF_2$ (the norm over $\FF_2$), $x\in\FF_3\mapsto -x^2\in\FF_3$ respectively. The lattice $\LL_4$ has trivial  discriminant form.
\end{lemma}
\begin{proof}
The last assertion follows from the fact that the $E_8$-lattice is unimodular.
The generator of $\LL_1^\vee/\LL_1$ is represented by $\theta^{-1}r$ and 
$\la \theta^{-1}r, \theta^{-1}r\ra=1$. Via the identification  
$\LL_1^\vee/\LL_1\cong \Ecal/\theta\Ecal\cong \FF_3$  that takes $\theta^{-1}r$ to $1\in \FF_3$, this becomes the form $x^2$.  Since $\LL_1$ primitively embeds in $\LL_4$ with orthogonal complement $\LL_3$, it follows that the discriminant form of $\LL_3$ is  isomorphic  to the opposite form $-x^2$.

For $\LL_2$ a generator of $\LL_2^\vee/\LL_2$ is represented by $\half(r_1+r_2)$, so that 
$\LL_2^\vee/\LL_2\cong \half\Ecal/\Ecal\cong\FF_2\otimes \Ecal\cong \FF_4$. Since
$\la \half(r_1+r_2), \half(r_1+r_2)\ra=\frac{3}{2}$ (which is of order $2$ in $\half\Ecal/\theta\Ecal$), the hermitian form takes on $1\in\FF_4$ the value $1\in \FF_2$ and so  the form is on $\FF_4$ is as stated (this  is in fact the only nontrivial hermitian form on $\FF_4$).
\end{proof}

Assume $i\le 4$. We then are in the situation discussed in the previous section: it was shown by Orlik-Solomon \cite{os}  that the discriminant of $V(\LL_i)\to \Gg(\LL_i)\bs V(\LL_i)$ is of Coxeter type $A_i$ and so the 
Galois group of the universal reflection cover relative this discriminant  is isomorphic with $\Sfrak_{i+1}$ 
(a uniform proof is due to Couwenberg \cite{couw} and appears in   \cite{chl} as Thm.\  5.1 and Thm.\  5.3).
Concretely, if $H_i\subset \CC^{i+1}$ is the hyperplane defined by $z_0+\cdots +z_i=0$ with its evident 
$\Sfrak_{i+1}$-action, then the formation 
of its orbifold quotient $H_i\to \Sfrak_{i+1}\bs H^i$ is given by the coefficients $\sigma_2, \dots \sigma_{i+1}$ of  the polynomial 
\[
(z-z_0)\cdots (z-z_i)=z^{i+1}+\sigma_2z^{i-1}-\sigma_3z^{i-2}\cdots +(-1)^{i+1}\sigma_{i+1}.
\]
The discriminant in $\Sfrak_{i+1}\bs H_i$ is the set of polynomials with a multiple root. So if $Z_i\subset H^i$ is the locus where at least two coordinates are equal, then we have an isomorphism  
$\Sfrak_{i+1}\bs (H^i, Z_i)\cong \Gg(\LL_i)\bs  (V(\LL_i), Z(\LL_i))$ and 
\[
(H^i, Z_i)\to (V(\LL_i),Z(\LL_i))_{\Gg(\LL_i)}
\]
becomes a universal reflection cover. For us, the cases $i=3$ and $i=4$ will be especially relevant.
 
The group $\Gg(\LL_1)$ equals $\mu_3$. The groups $\Gg(\LL_2)$,  $\Gg(\LL_3)$, 
 $\Gg(\LL_4)$ are numbered  by  Shephard and Todd in their list \cite{st} by resp.\ \textbf{4},   \textbf{25},  \textbf{32}. The polynomial algebra  $\CC[V(\LL_i)]^{\Gg(\LL_i)}$ has for $i=2,3, 4$ homogeneous generators of degree resp.\  $(4,6)$, $(6,9,12)$, 
 $(12, 18, 24, 30)$. The degrees have as greatest common divisor  the order  $C_{\Gg(\LL_i)}$ (as they must have), for 
the standard $\CC^\times$ action on $V(\LL_i)$ given by scalar multiplication induces a $\CC^\times$ action on $\Gg(\LL_i)\bs V(\LL_i)$ which factors through a faithful action  of $\CC^\times/C_{\Gg(\LL_i)}$. If we identify $\CC^\times/C_{\Gg(\LL_i)}$ with $\CC^\times$  via  $t\mapsto t^{\#( C_{\Gg(\LL_i)})}$, then the  resulting (faithful) $\CC^\times$ action  has degrees $(2,3)$,  $(2,3,4)$, $(2,3,4,5)$ respectively.
These are precisely  the degrees of homogeneous generators of $\CC[H_i]^{\Sfrak_{i+1}}$. 

The group of automorphisms of the pair $\Sfrak_{i+1}\bs (H^i, Z_i)$ is just  $\CC^\times$, the action being induced  by scalar multiplication in $H_i$ (see also Theorem 5.3 of \cite{chl}).
So any isomorphism $\Gg(\LL_i)\bs  (V(\LL_i), Z(\LL_i))\cong \Sfrak_{i+1}\bs (H^i, Z_i)$ is unique up to this  $\CC^\times$-action and that via such an isomorphism the $\CC^\times/C_{\Gg(\LL_i)}$-action on $\Gg(\LL_i)\bs V(\LL_i)$ lifts to scalar multiplication on $H^i$. In particular, we get a natural action of $\mu(\LL_i)\times \Sfrak_{i+1}$ on the pair $(H^i,Z_i)$.  
 
In the context of the generalization which we will discuss next,  we may write   $\hat V(\LL_i)$ for 
$H^i$, $\hat Z(\LL_i)$ for $Z_i$, $W(\LL_i)$ for $\Sfrak_{i+1}$ and  $\hat\G(\LL_i)$ for $\mu(\LL_i)\times \Sfrak_{i+1}$.
%\end{remark}
\smallskip

The following lemma makes clear the `building block' nature of these examples. (We recall that here an Eisenstein lattice always has its hermitian form take values in $\theta\Ecal$.)

\begin{lemma}[Allcock, \cite{allcock:y555} Thm.\ 3] \label{lemma:posdef}
Every positive  definite  Eisenstein lattice $\MM$ generated by its  $3$-vectors decomposes canonically into a direct sum of sublattices, each of which is isomorphic to some $\LL_i$ with $i\le 4$.\hfill $\square$
\end{lemma}

For a positive definite Eisenstein  lattice $\MM$, we set
\[
Z(\MM):=\text{union of triflection mirrors in } V(\MM).
\]
Suppose first that $\MM$ is an  orthogonal direct sum of  $k$ copies of $\LL_i$ with $i\le 4$. Then $\Gg(\MM)=\Gg(\LL_i)^k$ and
$\U(\MM)=\U(\LL_i)\wr \Sfrak_k=(\mu(\LL_i)\times \Gg(\LL_i))\wr\Sfrak_k$. We can also write this as a semidirect product of $\Gg(\MM)$ and 
$\mu(\LL_i)\wr \Sfrak_k$. 
Then $\Gg(\MM)\bs (V(\MM), Z(\MM))$ admits a reflection covering 
\[
(\hat V(\MM), \hat D(\MM)):= (H^i,Z_i)^k\to \Gg(\MM)\bs (V(\MM), Z(\MM)),
\]
where we have now an action of $(\CC^\times\times\Sfrak_{i+1})\wr \Sfrak_k$ on $\hat V(\MM)$ which contains $W(\MM):=\Sfrak_{i+1}^k$ as the Galois group.
The quotient group $\CC^\times\wr \Sfrak_k$ acts on $\Gg(\MM)\bs V(\MM)$ through a faithful action of  
$\CC^\times/C_{\Gg(\LL_i)} \wr \Sfrak_k$ and that action lifts to $\hat V(\MM)$. In particular, the group 
\[
\hat\G(\MM):=(\mu(\LL_i)\times \Sfrak_{i+1})\rtimes\Sfrak_k \cong \Sfrak_{i+1}^k\rtimes (\mu(\LL_i)\wr \Sfrak_k)=W(\MM)\rtimes (\U(\MM)/\Gg(\MM))
\]
naturally acts on this cover  and gives rise to  a reflection covering 
\[
(\hat V(\MM), \hat D(\MM))\to \U(\MM)\bs (V(\MM), Z(\MM)).
\]
with Galois group  $\hat\G(\MM)$.

More generally, if $\MM$ is positive definite and generated by its $3$-vectors, then by Lemma \ref{lemma:posdef} $\MM$ decomposes 
as $\MM_1\operp\cdots\operp\MM_4$ with $\MM_i$ an  orthogonal direct sum with $k_i$ copies of $\LL_i$ and $\U(\MM)=\U(\MM_1)\times\cdots \times \U(\MM_4)$,  $\Gg(\MM)=\Gg(\MM_1)\times\cdots \times \Gg(\MM_4)$  and we can apply the preceding discussion to each factor to produce a reflection covering $(\hat V(\MM), \hat D(\MM))\to \U(\MM)\bs (V(\MM), Z(\MM))$ with Galois group $\hat\G (\MM)=\hat\G(\MM_1)\times\cdots\times 
\hat\G(\MM_4)$. This group contains $W(\MM):= W(\MM_1)\times\cdots\times W(\MM_4)$ and has the form $W(\MM)\rtimes (\U(\MM)/\Gg(\MM))$. We record this as follows.

\begin{corollary}\label{cor:isodiscr}
Let $\MM$ be a positive definite Eisenstein lattice generated by its $3$-vectors.
Denote by $Z(\MM)$ the union of its triflection mirrors. Then there exists a universal reflection covering 
\[
(\hat V(\MM), \hat Z(\MM))\to \U(\MM)\bs (V(\MM), Z(\MM))
\]
whose Galois group is a  semidirect product $W(\MM)\rtimes (\U(\MM)/\Gg(\MM))$. \hfill $\square$
\end{corollary}

So for every subgroup $G\subset \U(\MM)$ which contains $\Gg(\MM)$, the above map factors through 
$G\bs (V(\MM), Z(\MM))$  with the first map defining a universal reflection covering with Galois group  the subgroup 
$\Gcal:= W(\MM)\rtimes  (G/\Gg(\MM))\subset W(\MM)\rtimes (\U(\MM)/\Gg(\MM))$. Such a situation will come up if we are given 
an overlattice $\MM'\supset\MM$ and take $G=U(\MM')\cap  \U(\MM)$.

We shall need some other  properties of lattices $\LL_i$ and it seems best to  derive them here.

\begin{lemma}\label{lemma:corankone}
Let $\MM\subset \LL_i$  (with $i=1,2,3,4$) be a sublattice of corank one spanned by $3$-vectors. If $\MM$ is $0$ ($i=1$) or of type $2\LL_1$ ($i=3$), $3\LL_1$ ($i=4$) or $\LL_3$ ($i=4$), then $\MM^\perp$ is spanned by a $3$-vector $r$. Otherwise $\MM$ is of type $\LL_1$ ($i=2$), $\LL_2$ ($i=3$) or 
$\LL_2\operp\LL_1$ ($i=4$) and  $\MM^\perp$ is spanned by a vector $a$ with $\la a,a \ra=6$.
\end{lemma}
\begin{proof}
This is rather straightforward. Let us just mention what   $r$ or $a$ is for the standard inclusions, leaving the rest as an exercise:
For $2\LL_1\subset \LL_3$, take  $r= r_1-\theta r_2-r_3$, for $\LL_3\subset \LL_4$, take  $r=r_1-\theta r_2-2r_3+\theta r_4$ and for  $3\LL_1\subset \LL_4$ (with $3\LL_1$ spanned by $r_1,r_3, r_1-\theta r_2-2r_3+\theta r_4$), take  $r= r_1-\theta r_2-r_3$.
For $\LL_1\subset \LL_2$ take $a=r_1-\theta r_2$, for both
$\LL_2\subset \LL_3$ and  $\LL_2\operp \LL_1\subset \LL_4$ (where the $\LL_1$-summand is spanned by
the $3$-vector $r_1-\theta r_2-2r_3+\theta r_4$), take $a=r_1-\theta r_2-2r_3$. 
\end{proof}

\begin{corollary}\label{cor:corankone}
The pairs $\{\LL_1, \LL_3\}$, $\{2\LL_1, 2\LL_1\}$,  $\{\LL_2, \LL_2\}$,  $\{\LL_1\operp\LL_2, \LL_1(2)\}$ have the property that their orthogonal direct sum 
embeds in $\LL_4$ such that each item is the orthogonal complement of the other (and so each summand is  primitively embedded in $\LL_4$). We also have  an embedding of $3\LL_1$ in $\LL_3$.
\hfill $\square$
\end{corollary}

\begin{remark}\label{rem:nonewmirrors}
The vector $a:=r_1-\theta r_2-2r_3\in \LL_3$ has the property that $\la r_1, a\ra =\la r_2, a\ra =0$, $\la r_3, a\ra=-3$ and $\la a, a\ra=6$ so that the 
unitary reflection $s_a\in \U(\LL_3)$ is defined (it is what we will later call a long root). So $\theta^{-1}a$ lies in $\LL_3^\vee$ and maps to a generator of the discriminant module $\LL_3^\vee/\LL_3\cong \FF_3$. The reflection $s_a$ acts on this  module as inversion and hence represents the nontrivial element of $\U(\LL_3)/\Gg(\LL_3)$. \end{remark}

Here is another  basic example.

\begin{example}[The hyperbolic Eisenstein lattice]\label{example:hyp}
 This lattice, which we shall denote by $\HH$, has two isotropic generators $e,f$ with $\la e,  f\ra=\theta$.
 So this lattice is $\theta$-self-dual. We note that it contains a pair of perpendicular $3$-vectors, namely  $e+\omega f$ and  
 $e+\bar\omega f$. They define a commuting pair of triflections in $\HH$ and hence $\Gg(\HH)$ contains the group of scalars $\mu_3$.
 \end{example}
 
\begin{lemma}\label{lemma:4-period}
The lattice $\LL_5$ is degenerate (its nondegenerate quotient is of type $\LL_4$) and for $k\ge 6$,
$\LL_{k}$ is isomorphic  with $\LL_4\operp\Hb\operp \LL_{k-6}$.
\end{lemma}
\begin{proof}
The kernel of $\LL_5$ is spanned by $e:=r_1-\theta r_2-2r_3+\theta r_4+r_5$.
If in $\LL_6$ we  put  $f:=-\bar\omega e+r_6$, then $e, f$ is a pair of isotropic vectors with $\la e, f\ra=\theta$.
So they span a copy of $\HH$ perpendicular to $\LL_4$.
Since $\LL_4\operp\HH$ is unimodular in the Eisenstein sense, it will be direct summand  of $\LL_{6+k}$.
Its orthogonal complement is spanned by $-e+r_7, r_8, \dots, r_{k}$, which is indeed a copy of $\LL_{k-6}$.
\end{proof}

\begin{conventionsnot}
We use the convention that  if $\lambda\in \ZZ$ and $\LL$ is an Eisenstein  lattice, then
$\LL(\lambda)$ stands for the lattice that has the same underlying 
$\Ecal$-module, but with  the form of $\LL$   multiplied  with $\lambda$.  We use the same convention for $\ZZ$-lattices.   For the hyperbolic lattice $\HH$, this makes even sense if we take $\lambda\in \Ecal$.
\end{conventionsnot}

\begin{example}\label{ex:HHproperty}
It is clear from this definition that  $\HH (\theta)$ is the Eisenstein lattice spanned by two isotropic 
vectors $e,f$ with inner product $-3$. Its discriminant group is $\Ecal/\theta\Ecal\oplus \Ecal/\theta\Ecal$ with a generators represented by  $\eps:=\theta^{-1}e$ and $\varphi:=\theta^{-1}f$. We have $\la \eps, \varphi\ra=-1 \in \QQ(\omega)/\theta\Ecal$ and so via the isomorphism 
$\Ecal/\theta\Ecal\oplus \Ecal/\theta\Ecal\cong \FF_3\oplus \FF_3$, the discriminant form of $\HH (\theta)$ becomes the hyperbolic form $(x,y)\in \FF_3\oplus \FF_3\mapsto -xy$.

For every $a\in \HH(\theta)$, we have $\la a, a\ra\in 6\ZZ$ and so 
unlike $\HH$, the lattice $\HH(\theta)$ does not contain $3$-vectors. On the other hand 
 the exchange of the basis vectors defines an involution of 
$\HH (\theta)$: this involution is defined by the unitary reflection $s_a$ with $a:=e-f$ (which indeed has the property that 
$\la a, a\ra=6$ and $\la x, a\ra\in 3\Ecal$ for all $x\in \HH (\theta)$).

Note that if we replace $f$ by $-f$ we find  that $\HH (\bar\lambda)\cong\HH (\lambda)$.
\end{example}

\begin{lemma}\label{lemma:interestingiso}
There exist the following lattice isomorphisms
\begin{gather*}
\LL_1\operp\LL_3\operp\HH\cong \LL_4\operp\HH(\theta),\\
\LL_1(2)\operp\HH\cong \LL_1\operp\LL_2\operp\LL_1(-1),\\
\LL_3 \operp\HH\cong \LL_4\operp \LL_1(-1).
\end{gather*}
 $\LL_1\operp\LL_3\operp\HH$ is  isomorphic to  $\LL_4\operp\HH(\theta)$.
The discriminant module of these lattices is a $\FF_3$-vector space of dimension one. 
\end{lemma}
\begin{proof}
We verify  the first isomorphism; the other two  are left to the reader.

Let $(r; r_1,r_2, r_3; e,f)$  be the standard basis of $\LL_1\operp\LL_3\operp\HH$. We put 
$r'_3:=r_3+e$ and $r'_4:= f+r$. Then $r_1,r_2, r'_3, r'_4$ are $3$-vectors span a copy $\LL'_4$  of $\LL_4$. 
Since $\LL_4$ is unimodular, this will be direct (orthogonal) summand. 
So it remains to show that the orthogonal complement of $\LL'_4$ is a copy of $\HH(\theta)$. 
One checks that this is spanned by $u:=r+\theta e$ and $v:=r_1-\theta r_2-2r_3 +\theta f$.
\end{proof}

\begin{remark}\label{rem:peculiariso}
The last two isomorphisms of Lemma \ref{lemma:interestingiso} imply that  there exists an isomorphism
\[
f: \LL_2\operp \LL_3 \operp\HH\cong \LL_2\operp \LL_4\operp \LL_1(-1)\cong \LL_1(2)\operp \LL_4\operp \HH
\]
Note that   $\LL_2\operp \LL_3$ is generated by $3$-vectors, but  $\LL_1(2)\operp \LL_4$ and so these two lattices are not isomorphic.
This tells us that $f$ cannot take an  isotropic 
generator of the $\HH$-summand of  its source to an  isotropic  generator of the $\HH$-summand of its target.
Said differently: the isotopic rank one lattices spanned by them  represent different cusps.
 \end{remark}

\begin{definition}\label{def:special}
Let $\LL$  be an Eisenstein lattice.  We call $a\in \LL$ an \emph{$n$-vector} if $n=\la a,a\ra$ (so $n$ must be an integer divisible by $3$).
If $a$ is in addition primitive, then we say  that $a$  (or
the $\Ecal$-submodule  of $\LL$  spanned by such a vector) is \emph{special}, if its inner product with any lattice vector is divisible by $3$ (or equivalently,  that $\theta^{-1}a$ represents an element   of order three in the discriminant module
$\disc(\LL)=\LL^\vee/\LL$ and \emph{standard} otherwise.
\end{definition}

Note that if $a\in \LL$ is special, then $\theta^{-1}a$ represents an element   of order three in the discriminant module
$\disc(\LL)=\LL^\vee/\LL$.

\begin{remarks}\label{rem:special}
Most relevant here are the special vectors with self-product $3$, $0$ and $6$, to which we shall refer  as resp.\ \emph{special $3$-vectors}, \emph{special isotropic vectors} and \emph{long roots}. We make a few observations about these cases.
\begin{description}
\item[special $3$-vectors] It is clear that a $3$-vector $r$  is special if and only if it spans a direct (orthogonal) summand of type $\LL_1$ of $\LL$. In that case
 the hexaflection $t'_r$ is defined. Furthermore, its third power, the reflection  $s_r$ takes the order three element 
 $[\theta^{-1}a]\in \LL^\vee/\LL$ to   $-[\theta^{-1}a]$.
 \item[special isotropic vectors] 
A  basis vector of  $\HH(\theta)$  is special and hence spans  a special  isotropic submodule.
\item[long roots]  If $a\in \LL$  is such that $\la a,a\ra=6$, then it is clear from the definition that $a$ is special if and only if  the unitary reflection $s_a$ preserves $\LL$. When that is the case,   $[\theta^{-1}a]\in\LL^\vee/\LL$  has order $3$ and 
$s_a$ takes it to $-[\theta^{-1}a]$.
\end{description}

\end{remarks}

For example, the element $a=r_1-\theta r_2-2r_3\in \LL_3$ is a long root and hence $s_a$ is defined (and acts  nontrivially on the discriminant module
$ \LL_3^\vee/\LL_3\cong \FF_3$). The lattice  $\LL_1\operp\LL_3$
has besides long roots also special $3$-vectors. Via Lemma \ref{lemma:interestingiso} we see that  $\LL_1\operp\LL_3\operp \HH$ has also special  isotropic sublattices $\II$; these  can also be distinguished by the type of the positive  definite lattice $\II^\perp/\II$: its type is then
 $\LL_1\operp\LL_3$ resp.\ $\LL_4$.
 
By definition, a  $3$-vector that is not special is standard. It is clear that every $3$-vector contained in sublattice of type $\LL_2$ is standard.

\subsection{Ball quotients over the Eisenstein field}
We will be mostly concerned with lattices $L$ of signature $(2n,2)$ that are generated by its $2$-vectors. 
Then $\LL$ has \emph{Lorentzian} signature $(n,1)$ and is generated by its $3$-vectors. 
The group  $\U(\LL)$ is arithmetic in $\U(V(\LL))$ and hence a lattice in that group (in the sense that it is a discrete subgroup of finite covolume). 
The set of $z\in V(\LL)$ with $\la z,z\ra<0$ make up an open $\CC^\times$-invariant subset $V_-(\LL)\subset V(\LL)$. This defines an open 
ball $\BB(V)\subset \PP(V)$ such that $V_-(\LL)\to \BB(\LL)$ is a $\CC^\times$-bundle on which $\U(V(\LL))$ acts.
The ball $\BB(V)$  is in fact a symmetric domain for $\U(V(\LL))$ (or $\PU(V(\LL))$). The Baily-Borel theory \cite{BB} then tells us that  
the complex  analytic orbifold  $\U(\LL)\bs \BB(V(\LL))$ naturally lives in the  quasiprojective category.

We will often write  $\BB(\LL)$  for  $\BB(V(\LL))$.

\begin{remark}\label{rem:hodge-interpretation}
We can think of  $\BB(\LL)$ as parametrizing polarized  Hodge structures on $L$ with 
$h^{1,-1}=h^{-1,1}=1$, $h^{0,0}=2n$ (so of weight $0$, but we are of course free to make a 
Tate twist  and thus  change the weight by an even integer) endowed with a $\mu_3$-action.  First note
that we have a natural embedding $j: V(\LL)\subset \CC\otimes L$ which up to a scalar $\frac{3}{2}$ preserves hermitian forms
if we take on  $\CC\otimes L$ the hermitian  extension of the given  pairing on $L$.
If $z\in \BB$ defines the  negative line $\ell\subset V(\LL)$, then  we assign to $z$  the Hodge structure $\Hl_z^{\pt,\pt}$ on $L$ by 
\[
\Hl_z^{1,-1}:=j(\ell), \quad \Hl_z^{-1,1}:=\overline{j(\ell)} \text{ and  } \Hl_z^{0,0}:=\big(\Hl_z^{1,-1}+ \Hl_z^{-1,1}\big)^\perp
\]
So the plane $\Hl_z^{1,-1}+\Hl_z^{-1,1}$ is negative definite and its orthogonal complement $\Hl_z^{0,0}$  is positive definite. This Hodge structure is polarized by the given pairing on $L$. The orbifold  $\U(\LL)\bs\BB(\LL)$ underlies a Shimura variety
(strictly speaking, we should take the disjoint union of two copies, one for each embedding of $\Ecal$ in $\CC$ and then divide by the $\Orth(L)$-normalizer  of  the order $3$ subgroup $\la\rho\ra\subset \Orth(L)$).
\end{remark}

\subsection{Reflection covers of Eisenstein ball quotients} In  this subsection we assume that $\LL$ is an  Eisenstein lattice of Lorentzian signature. After Lemma \ref{lemma:loc_coxeter} we also assume that $\LL$ is generated by its $3$-vectors. 
\\

 We say that  a hyperball in $\BB(\LL)$ is a \emph{mirror} of $\LL$ if it is the  fixed point set of a pseudoreflection in $\U(\LL)$. 
The mirrors in  $\U(\LL)$  are locally finite on $\BB(V)$. So their union in $\BB(V)$ (the \emph{mirror arrangement}) is closed and  for every finite index subgroup 
$\G\subset \U(\LL)$ its image in $\G\bs \BB(\LL)$ a hypersurface.  If we restrict to resp.\ the standard 
$3$-vectors, the special $3$-vectors, the  long roots, then we denote the corresponding  
mirror arrangement resp.\  $Z(\LL)^\std$, $Z(\LL)^\spl$, $Z(\LL)^\lr$ and we shall refer to its image  $\BB(\LL)_\G$ as resp.\ the \emph{standard, special, long root mirror divisor} with respect to $\G$.

\begin{lemma}\label{lemma:loc_coxeter}
Let $Z\subset \BB(\LL)$ be the standard arrangement $Z(\LL)^\std$,  the special arrangement  $Z(\LL)^\spl$ or  their union $Z(\LL)^\std\cup Z(\LL)^\spl$. Then the  pair $(\Xcal, D):=\U(\LL)\bs (\BB(\LL),Z)$ is locally of Coxeter type and hence admits a universal reflection covering
$(\hat \Xcal, \hat D)\to (\Xcal,D)$. 
\end{lemma}
\begin{proof}
Let  $p\in \BB(\LL)$ and denote by  $\Gg(p)\subset \U(\LL)_p$ the subgroup  generated by the
triflections whose fixed point set in $\BB(\LL)$ is contained in $Z$ and contains $p$. By Corollary \ref{cor:inheritance},  it suffices to prove that the  discriminant of $\BB(\LL)\to 
\Gg (p)\bs \BB(\LL)$ is of Coxeter type at $p$. We note that $p$  defines a negative definite line $\ell\subset V(\LL)$ and that the  $3$-vectors 
perpendicular to this line span a positive definite sublattice  $\MM\subset \LL$. It is clear that $\Gg(p)$ can be identified with  $\Gg(\MM)$.
By  Lemma \ref{lemma:posdef} the discriminant of $V(\LL) \to \Gg(\MM)\bs V(\LL)$ (and hence of $\BB(\LL)\to \Gg (p)\bs \BB(\LL)$) is of Coxeter type. This settles  the case $Z=Z(\LL)^\std\cup Z(\LL)^\spl$. The other two cases follow from this. 
\end{proof}

In the rest of this subsection we assume that $\LL$ is generated by its $3$-vectors.

Let $\MM\subset\LL$ be a positive definite sublattice  generated by $3$-vectors. So $\MM^\perp$ is Lorentzian and $\BB(\MM^\perp)$ is an intersection of mirrors.
Denote by $\U(\LL)_{[\MM]}$ the  $\U(\LL)$-stabilizer of $\MM$, so that we have a natural map 
\begin{equation}\label{eqn:normalozation}
f: \U(\LL)_{[\MM]}\bs \BB(\MM^\perp)\to \U(\LL)\bs \BB(\LL).
\end{equation}
This map is  the normalization of its image: it is clearly finite and it is also of degree one, for  the $\U(\LL)$-orbit of a  generic point of $\BB(\MM^\perp)$ meets
$\BB(\MM^\perp)$ in its $\U(\LL)_{[\MM]}$-orbit. 

\begin{proposition}\label{prop:mirrortrace}
A proper intersection of a mirror of a $3$-vector in $\LL$ with  $\BB(\MM^\perp)$ is in fact  the mirror of a $3$-vector or a long root in $\MM^\perp$.  

If $\MM$ is spanned by a special $3$-vector, then the normalization map \eqref{eqn:normalozation} is injective. 
\end{proposition}
\begin{proof}
Let $r\in  \LL$ be a $3$-vector whose mirror meets $\BB(\MM^\perp)$ in a hyperball. 
Then $\Ecal r+\MM$ is positive definite and spanned by its $3$-vectors. 
We decompose this lattice  according to 
Lemma \ref{lemma:posdef}. The summand $\MM'\subset \Ecal r+\MM$ which contains $r$  is of type 
$\LL_i$ for some $i\le 4$. We must show that $\MM'\cap \MM^\perp$ contains $3$-vector or a long root.
If $i=1$, then  $r\in \MM^\perp$ and we are done, of course. So assume  $i\ge 2$. Then   $\MM'\cap\MM$ has corank $1$ in $\MM'\cong \LL_i$. It will be a subsum of the decomposition of $\MM$ given by Lemma \ref{lemma:posdef} and hence is spanned by $3$-vectors. 
By Lemma \ref{lemma:corankone}  the rank one lattice $\MM'\cap\MM^\perp$ is then  spanned by a $3$-vector or a long root. 
This proves the first assertion.

For the second assertion, suppose $\MM$ spanned by the special $3$-vector $r$. 
Let $z\in \BB(\MM^\perp)$ and $g\in \U(\LL)$ be such that  $g(z)\in \BB(\MM^\perp)$. 
We must show that there exists a $g'\in \U(\LL)$ which preserves $\MM$ and takes $z$ to $z'$. 
Since $r$ is special, so is $g(r)$ and as $g(z)$ is perpendicular to both $r$ and $g(r)$,  these special vectors  span a positive definite sublattice.  By definition each of them  splits off a perpendicular summand. This can only happen  if  they are  proportional or perpendicular. In the first case, $g$ preserves $\MM$ and so $g':=g$ will do. In  the last case 
$\Ecal r+\Ecal r'$ is an orthogonal direct summand of $\LL$. Hence there is $g_1\in\U(\LL)$ which exchanges $r$ and $r$
and is the identity on $(\Ecal r+\Ecal r')^\perp$. Then $g':=g_1g$ takes $z$ to $g(z)$ and preserves $\MM$.
\end{proof}

\subsection{Allcock's conjecture}\label{subsect:ac}
The case of prime interest here is when $L$ is even unimodular of signature $(26,2)$. Such a lattice is unique up to isomorphism, and for example isomorphic 
to the orthogonal direct sum of three $E_8$ lattices and two standard hyperbolic lattices. According to Allcock (see the proof of Thm.\  5 in \cite{allcock:y555}), the   $\rho\in \Orth(L)$ that satisfy $\rho^2+\rho+1=0$ form a single $\Orth(L)$-conjugacy class  and hence the Eisenstein lattices that have $L$ as underlying integral lattice are mutually isomorphic. We choose $\rho$ such that the resulting Eisenstein lattice $\LL$ is $3\LL_4\operp \HH$. 
\smallskip

\begin{convention}
We shall refer to $\LL$, whose signature is $(13,1)$, as the \emph{Allcock lattice} and 
write  simply  $\G$ for  $\U(\LL)$ and $\BB$ for  $\BB(\LL)$. \emph{In particular, from now on,
$\LL$ stands for this lattice.}
\end{convention}

Allcock and Basak (in some cases jointly) prove that 
\begin{enumerate}
\item[(i)] $\G$ is generated by its  triflections (so $\G=\Gg(\LL)$) and its triflections  make up a single conjugacy class in $\G$ (\cite{basakI}, \cite{allcock:y555}),
\item[(ii)] $\G$ has five orbits in the set of primitive isotropic sublattices $\II\subset \LL$ and  these can be distinguished by the type of  the (positive definite) sublattice of $\II^\perp/\II$ generated by its $3$-vectors. These types are $(0)$ (the underlying $\ZZ$-lattice is then the  Leech lattice), $12\LL_1$, $6\LL_2$, $4\LL_3$ and $3\LL_4$ (\cite{allcock:y555}, Thm.\ 4).  
\item[(iii)] the union $Z$  of the triflection mirrors in $\BB$ considered as a reduced divisor, the zero divisor of a $\G$-automorphic form $\Psi_0$ of weight $4$ (\cite{allcock:leech}, Thm.\ 7.1).
\end{enumerate}

Property (i) implies that $\G$ acts  transitively on the mirrors and so $Z_{\G}$ is irreducible in $\BB_{\G}$. 

Property (ii) shows that whenever  $\II^\perp/\II$  contains a $3$-vector, then the $3$-vectors span an Eisenstein sublattice of maximal rank $12$. (One can show that this is formally consequence of (iii).) 

We denote the ball quotient $\G\bs (\BB, Z)$ by $(\Xcal, D_\Xcal)$.

\begin{corollary}\label{cor:doublecover}
There exists a connected $\mu_2$-covering $\tilde \Xcal\to \Xcal$ whose ramification locus is $D_\Xcal$.
\end{corollary}
\begin{proof}
By property  (iii) above, the divisor $Z$ is defined by a  $\mu_3$-valued automorphic form $\Psi_0$ of weight 4. We may think of $\Psi_0$ as a
 $\G$-invariant function $V^{\A-}/\mu_3\to \CC$ which is homogeneous of degree $-4$ with respect to the $\CC^\times/\mu_3$-action on
 $V^{\A-}$. Since $-1\in \G$, we  rather view $\Psi_0$ as a  $\G$-invariant function $V^{\A-}/\mu_6\to \CC$ which is homogeneous of degree $-2$ with respect to the $\CC^\times/\mu_6$-action. In other words,  if $\chi : \CC^\times =\GL_1(\CC)$ is the tautological character and 
 $\Lcal\to\BB$ denotes the holomorphic line bundle over $\BB$ defined by  $V^{\A-}\otimes _{\chi^{-6}} \CC$, then $\G$ acts on this bundle and 
 $\Psi_0$ can be regarded as a  
$\G$-invariant section of  $\Lcal^{\otimes 2}$ with zero divisor $Z$. The roots of  $\Psi_0$  in $\Lcal$ then define a double cover of $Z$ as asserted.
 \end{proof}

\begin{conjecture}[Allcock \cite{allcock:proposal}]\label{conj:allcock}
The universal reflection Galois group  of the pair $(\Xcal, D_\Xcal)$ is isomorphic with  the Bimonster 
$\Mb\wr\mu_2= (\Mb\times \Mb)\rtimes \mu_2$, where $\Mb$ is the largest sporadic finite simple  group (the Monster).
\end{conjecture}

In view of Corollary \ref{cor:doublecover}, this is equivalent to the  assertion  that the orbifold fundamental group of 
$\tilde \Xcal$  is isomorphic to $\Mb\times\Mb$ (with the covering involution of $\tilde \Xcal$ exchanging the two copies). 

\subsection{Three ball quotients with a modular interpretation}\label{subsect:ballmodular}
In what follows, we will be concerned with certain  Eisenstein sublattices of the lattice $\LL$ 
of the form $\MM^\perp$, where $\MM$ is a positive definite sublattice generated by $3$-vectors.  So by 
Proposition \ref{prop:mirrortrace},  each proper intersection of a mirror of $\LL$ with $\BB(\MM^\perp)$ is the mirror of a pseudoreflection in $\U(\MM^\perp)$.

 Recall that $\G(\MM^\perp)\subset \U (\LL)$ stands for the subgroup of unitary transformations that act trivially on the discriminant lattice of $\MM^\perp$ and can be identified with the subgroup of unitary transformation of $\LL$ that act as the identity on $\MM$. So if      
$\U (\MM^\perp)/\G(\MM^\perp)$ is represented by scalars, then  $\U(\MM^\perp)$ is the image of the $\G$-stabilizer of $\MM^\perp$.

Especially relevant for this paper are:
\begin{enumerate}
\item[$\LL^\ct$:] $=\LL_1\operp 2\LL_4\operp\HH$ and  $\MM=\LL^0_3\subset \LL$, where $\LL^0_3$ is the span of the first three basis vectors  of the first $\LL_4$-copy.
Its discriminant module is the one of the $\LL_1$-summand and hence a one-dimensional $\FF_3$-vector space.
We have  $\U(\LL^\ct)=\{\pm 1\}\times \G(\LL^\ct)$, and  we may represent the nontrivial element of $\U(\LL^\ct)/\G(\LL^\ct)$ by a hexaflection in the $\LL_1$ summand. 
\item[$\LL^\K$:] $=\LL_1\operp\LL_3\operp\LL_4\operp\HH$ and  $\MM=\LL^0_3\operp\LL^0_1\subset \LL$ (the obvious embedding in the first two $\LL_4$-copies). This lattice is perpendicular 
to a standard $3$-vector in $\LL^\ct$. Its discriminant module comes from the (mutually perpendicular) $\LL_1$ and the $\LL_3$ summands and hence is canonically the direct sum of two 
one-dimensional $\FF_3$-vector spaces. As $\LL_1$ and $\LL_3$ are each others orthogonal complement  when  embedded in $\LL_4$, 
they have by Remark \ref{rem:discriminant_duality} opposite discriminant forms. These are distinct  and so we have a canonical isomorphism 
$\U(\LL^\K)/\G(\LL^\K)\cong \mu_2\times\mu_2$. Its  generators are represented  by the reflection in the $\LL_1$-summand resp.\ the reflection in a 
long root in the $\LL_3$-summand.  The latter can also be represented by the minus the identity in that $\LL_3$-summand and hence the product of 
these generators is represented by  $-1\in \U(\LL^\K)$. So if  $\tilde \G_K$ is the subgroup of $\U(\LL^\K)$ generated by $\G_K$ and the reflections in  
the long roots, then  $\PGtilde (\LL^\K)=\PU(\LL^\K)$ and  $[\PU(\LL^\K):\PG^\K]=2$.

Lemma \ref{lemma:interestingiso} shows that $\LL^\K\cong  2\LL_4\operp\HH(\theta)$ and so $\LL^\K$ has  special isotropic vectors.
\item[$\LL^{\DM}$:] $=2\LL_4\operp\HH$  and   $\MM=\LL^0_4\subset \LL$, where $\LL^0_4$ is the first $\LL_4$-copy. This lattice is unimodular and is perpendicular to a special $3$-vector in $\LL^\ct$  in $\LL$. We have $\G(\LL_{10})=\U(\LL_{10})$.  By Lemma \ref{lemma:4-period}, this lattice is isomorphic with  $\LL_{10}$. Deligne-Mostow and Mostow  (\cite{dm}, \cite{mostow}) proved that the ball quotient $\U(\LL_{10})\bs \BB(\LL_{10})$ is isomorphic  with
$\Hcal^{0,\st}_{12}$, the  moduli space of  smooth rational curves endowed with a stable effective degree 12 divisor (so this is a Deligne-Mostow ball quotient). 
\end{enumerate}

The  other ball quotients  also  have a modular interpretation  (recalled below) and our notation refers to this: \emph{ct} for cubic threefold, \emph{DM} for Deligne-Mostow and \emph{K} for Kond\=o. We shall make use of these interpretations  for the construction of the universal reflection covers.  
We also abbreviate 
\begin{gather*}
\G^\ct:=\G(\LL^\ct),\;  \PG^\ct:=\G(\LL^\ct),\; \BB^\ct:=\BB(\LL^\ct),\; \\
 (\Xcal_\ct, D_\ct):=\G(\LL^\ct)\bs (\BB^\ct, Z^\ct), \; 
 Z^{\ct, \std}:=Z(\LL^\ct)^\std,\; \text{etc.} 
\end{gather*}
 and do likewise for the other cases.

\section{Classification  of geometric  sublattices of the Allcock lattice}\label{sect:appC} 
We shall see that moduli spaces of cubic threefolds with a given singularity configuration have the structure of a ball quotient, where the ball is a subball  of the Allcock ball $\BB$ which is defined by sublattices of the Allcock lattice of a particular type.  One of our goals is to classify these sublattices. We will not pursue a complete classification,  but focus  on the ones that are relevant here, namely those that are described by  Definition \ref{def:geomlattice} below.

\begin{definition}\label{def:geomlattice}
We say that a sublattice $\MM\subset \LL$ is \emph{geometric (of type $(p,q)$)} if it is nonzero and isomorphic to $p\LL_3\operp q\LL_4$. Its \emph{outer automorphism group} $\out_{\LL}(\MM)$ is a group of permutations of its factors induced by $\G_{[\MM]}$, the group  of  unitary transformation of $
\LL$ which preserve $\MM$. 

We regard  the collection of all geometric sublattices of $\LL$ as  a partially ordered set for the inclusion relation and denote it  by $\Lscr$.
\end{definition}

The similarity in terminology with notions related of a geometric subgroup of the Bimonster (introduced in Definition \ref{def:geometric subgroup}) is deliberate, for one of our main goals is to establish a correspondence between the two.

The discussion of the Examples \ref{example:basiccases} and Lemma \ref{lemma:discrLi} shows that 
\begin{multline*}
\U(p\LL_3\operp q\LL_4)=\big((\mu_2\times \Gg(\LL_3))\wr \Sfrak_p\big)\times \big(\Gg(\LL_4)\wr \Sfrak_q\big)=\\=
\Gg(p\LL_3 \operp  q\LL_4)\rtimes ((\mu_2\wr\Sfrak_p)\times\Sfrak_q).
\end{multline*}

It is clear that the image of $\G_{[\MM]}$ in $\U(\MM)$ contains $\Gg(\MM)$ and  so this image is  a semidirect product of $\Gg(\MM)$ by a  group  which via an isomorphism $\MM\cong p\LL_3\operp q\LL_4$ is  realized as a subgroup of  $(\mu_2\wr\Sfrak_p)\times\Sfrak_q$. That subgroup is clearly an extension  of  $\out_{\LL}(\MM)$ by a subgroup of $\mu_2^p$. We shall denote that subgroup by $\mu_\MM$ and 
regard  it as a subgroup of the image of $\G_{[\MM]}$ in $\U(\MM)$.

The kernel of this restriction map $\G_{[\MM]}\to \U(\MM)$ is the group $\G(\MM^\perp)$ (which we recall, stands for the group of  unitary transformations of $\MM^\perp$ which act trivially on its discriminant module). 
It is clear that 
\[
\G_{[\MM]}\cong (\mu_\MM\times \Gg (\MM)\times \G (\MM^\perp)). \out_{\LL} (\MM).
\]

\subsection{Primitive embeddings}
We first classify the \emph{primitive} embeddings  of $k_\pt\LL_\pt=k_1\LL_1\operp k_2\LL_2\operp k_3\LL_3\operp k_4\LL_4$  in a unimodular lattice of 
Lorentzian signature. Such  questions have been well studied for integral  lattices, 
notably by Nikulin, and indeed, we build on one his results in an essential way.
Let us state this in a form that we will use.  The following is a consequence of  Nikulin's Thm.\ 1.14.4 in \cite{nikulin}, as we make stronger assumptions.
For a finite abelian group $A$ and a prime $p$, we denote by $\ell_p(A)$ the minimal number of generators of the $p$-primary part $\ZZ_p\otimes A$ of $A$.

\begin{proposition}[Nikulin]\label{prop:nikulin}
Let $L$ be an even unimodular integral lattice and  $M\subset L$ a primitive sublattice. Assume that for 
every prime $p$,  
\begin{equation}\label{eqn:nikulin}
\rk(L/M)>2+\ell_p(M^\vee/M).
\end{equation}
Then every primitive embedding $M\hookrightarrow L$ is obtained by composing $M\subset L$ with an orthogonal transformation. \hfill $\square$
\end{proposition}

We shall refer to the inequalities as \eqref{eqn:nikulin} as the \emph{Nikulin condition}. 

According to the discussion in \ref{example:basiccases} 
the discriminant module of $k_\pt\LL_\pt$ is a sum of two vector spaces over finite fields:
\[
k_\pt\LL_\pt^\vee/k_\pt\LL_\pt\cong \FF_3^{k_1+k_3}\oplus \FF_4^{k_2}.
\]
We also want to identify the discriminant form on the left hand side in terms of the right hand side.
We can of course do this per summand.

If we are given a primitive embedding of $k_\pt\LL_\pt$ in a unimodular Eisenstein lattice, then
the orthogonal complement $\NN$ of this embedding has the property that 
its discriminant module is isomorphic to that of $k_\pt\LL_\pt$. Hence  $\NN$ must have at least 
$\max\{k_1+k_3,k_2\}$ generators, in other words, the corank  of the embedding must be at least 
this number. The following construction comes close.

Recall from Corollary \ref{cor:corankone} that  the pairs $\{\LL_1, \LL_3\}$, $\{2\LL_1, 2\LL_1\}$, 
$\{\LL_2, \LL_2\}$,  $\{\LL_1\operp\LL_2, \LL_1(2)\}$ have the property that their orthogonal direct sum 
embeds in $\LL_4$ with each individual item embedding primitively. This  implies that
$k_\pt\LL_\pt$ embeds primitively in $\psi(k_\pt)\LL_4$, where 
\[
\psi(k_\pt):=\lceil \half(k_1-k_2)_+\rceil +k_2+k_3+k_4
\]
(the first term is zero when $k_1\le k_2$ and $\lceil\half(k_1-k_2)\rceil$ otherwise). To see this,  
note that for  $k_1\le k_2$ we can rewrite $k_1\LL_1\operp k_2\LL_2$ as 
$k_1(\LL_1\operp \LL_2)\operp (k_2-k_1)\LL_2$ and observe that this embeds in $k_2\LL_4$. 
When  $k_1>k_2$, then for  $k_1-k_2=2n$ even resp.\ $k_1-k_2=2n-1$ odd 
(so $n=\lceil \half(k_1-k_2)\rceil$ in either case),  we rewrite this as 
$k_2(\LL_1\operp \LL_2)\operp n(2\LL_1)$ resp.\  $k_2(\LL_1\operp \LL_2)\operp (n-1)(2\LL_1)\operp\LL_1$. 
Both lattices embed in $(n+k_2)\LL_4$.

By inspection, we see that  the corank of the embedding is as follows:
\[
4\psi(k_\pt)-\rk_\Ecal(k_\pt\LL_\pt)=
\begin{cases}
-k_1+2k_2+k_3 & \text{if  $k_1\le k_2$;}\\
k_1+k_3 & \text{if  $k_1-k_2$ is positive and even;}\\
k_1+k_3 +2 & \text{if  $k_1-k_2$ is positive and odd.}\\
\end{cases}
\]
In all cases the right hand side is $\ge \max\{k_1+k_3,k_2\}$, as it should. 

So when $m\ge \psi(k_\pt)$, we have defined a primitive embedding $\MM(k_\pt)\hookrightarrow m\LL_4\operp \HH$.

\begin{lemma}\label{lemma:nikulinforM}
Then the  embedding of the underlying  integral lattices (so of $M(k_\pt)$ in $m E_8\operp 2H$)  satisfies the Nikulin condition. 
\end{lemma}
\begin{proof}
The discriminant group of $M(k_\pt)$ is  isomorphic to   $\FF_3^{k_1+k_3}\oplus  \FF_2^{2k_2}$. So 
$\ell_2(M(k_\pt))=k_1+k_3$,  $\ell_3(M(k_\pt))=2k_2$ and $\ell_p(M(k_\pt))=0$ for all primes $p>3$.
Consider the equality  
\[
-2+ 2\rk_\Ecal (m\LL_4\operp \HH)-2\rk_\Ecal k_\pt\LL_\pt=2+ 2\big (4\psi(k_\pt)-\rk_\Ecal k_\pt\LL_\pt\big).
\]
We found  that the expression in parenthesis in the right hand side is $\ge \max\{k_1+k_3,k_2\}$ and so the Nikulin  condition is satisfied. 
\end{proof}

We will use this  to prove the  corresponding assertion for the Eisenstein lattices.

\begin{lemma}\label{lemma:uniqueembedding}
If $\psi(k_\pt)\le m$, then  every primitive embedding  $j: k_\pt\LL_\pt\hookrightarrow m\LL_4\operp\HH$ is 
 unitary equivalent to the one defined above.
\end{lemma}
\begin{proof}
We proceed with induction on the number of summands $|k_\pt|=k_1+k_2+k_3+k_4$.
We first reduce to the case $k_4=0$. A primitive embedding of $k_4\LL_4$ in  $m\LL_4\operp\HH$  with $k_4\le m$ 
has as orthogonal complement a unimodular lattice of Lorentzian signature. 
By a theorem of  Basak (\cite{basakI}, Lemma 2.6)  any such Eisenstein lattice is isomorphic  with  
$(m-k_4)\LL_4\operp \HH$. This implies that all embeddings of $k_4\LL_4$ in $m\LL_4\operp\HH$ are unitary equivalent. 
So $j$ is equivalent to one which on the $\LL_4$ summands is already in standard form. 
Hence we may assume that $k_4=0$.

We next reduce to the case $k_3=0$. 
Suppose that $k_3>0$  and split off a $\LL_3$-summand: $k_\pt\LL_\pt$ as $\LL_3\operp k'_\pt\LL_\pt$.
For the standard primitive embedding of $k_\pt\LL_\pt$ in $m\LL_4\operp\HH$ described above there exists a $3$-vector  in $m\LL_4\operp\HH$ which  is perpendicular to $k'_\pt\LL_\pt$  and spans with the split off $\LL_3$  a copy of $\LL_4$. 
If we pass to the underlying integral lattices, then by Lemma \ref{lemma:nikulinforM} above, 
$j$ is orthogonally equivalent to the given embedding $k_\pt\LL_\pt\subset m\LL_4\operp\HH$.
Hence we can find a $3$-vector $r$  in $m\LL_4\operp\HH$ which is perpendicular to 
$j(k'_\pt\LL_\pt)$ and such that  $\ZZ r+j(\LL_3)$ is contained in a $E_8$-summand. 
This implies that $\Ecal r+j(\LL_3)$ is positive definite of rank $4$, spanned by its $3$-vectors and indecomposable. By Lemma \ref{lemma:posdef}, this will be a copy of $\LL_4$.  So by a unitary transformation, we can arrange that $j$ is the orthogonal direct sum of a primitive embedding  
$\tilde j: \tilde k'_\pt\LL_\pt\hookrightarrow (m-1)\LL_4\operp\HH$ and the standard embedding $\LL_3\subset \LL_4$.
With induction, we then reduce to the case when $k_3=0$ also.

The same  argument works if both $k_1$ and $k_2$ are positive, since we can split off $\LL_1\operp\LL_2$ from $k_\pt\LL_\pt$ and use the fact 
that the standard embedding of $\LL_1\operp\LL_2$ in $\LL_4$ is perpendicular to a $3$-vector. This allows us to reduce to the case when 
$k_\pt\LL_\pt$  equals $k_1\LL_1$ or $k_2\LL_2$.

If $k_\pt\LL_\pt=k_2\LL_2$ with $k_2>0$,  then, since  $k_2=\psi (k_\pt)\le m$,
the same trick shows that $j$ factors through an embedding of  $(k_1-1)\LL_2\operp\LL_3$ 
and this is a case already taken care of.  So can also assume that $k_2=0$. 

It remains to do the case when $k_\pt\LL_\pt=k_1\LL_1$ with $k_1>0$.  Then 
$\lceil \half k_1\rceil=\psi (k_\pt)\le m$ and another such application shows that  
$j$ factors through an embedding of  
$(k_1-1)\LL_1\operp\LL_2$ and so we are done.
\end{proof}

\begin{proposition}\label{prop:uniqueembeddinghalf}
Let  $\MM\subset \LL$ be a primitively embedded  geometric sublattice isomorphic with  $p\LL_3\operp q\LL_4$. Then $p+q\le 3$
and  the restriction map $\G_{[\MM]}\to \U(\MM)$ is  onto with kernel $\G(\MM^\perp)$. More specifically, 
\[
\G_{[\MM]}\cong (\mu_\MM\times \Gg (\MM)\times \G (\MM^\perp))\rtimes  \out_{\LL} (\MM),
\]
where $\mu_\MM$ is the group generated by antipodal involution $\LL_3$-summands (so $\mu_\MM\cong \mu_2^p$) and  $\out_{\LL} (\MM)$ is the \emph{full} permutation group of the factors of $\MM$ which preserve type (so $\out_{\LL} (\MM)\cong \Sfrak_p\times\Sfrak_q$).

Any primitively embedded  geometric sublattice of $\LL$ isomorphic with $\MM$ is a $\G$-translate of $\MM$. 
\end{proposition}
\begin{proof}
Since $\MM$ is a sublattice of $\LL$, we have (for reasons of signature) that $3p+4q\le 13$. 
So if  $p+q>3$, then  $(p,q)=(3,1)$ or $(4,0)$.  
We claim that a primitive embedding of $4\LL_3$ in $\LL$ does not exist, so that  
\emph{a fortiori} a primitive embedding of   $3\LL_3\operp\LL_4$ in $\LL$ does not exist. 
This is because $\disc(4\LL_3)$ is a $4$-dimensional  $\FF_3$-vector space,  
so that the corank of any primitive  embedding must be at least $4$. But here the corank of such 
an embedding is going to be $14-4\cdot 3=2$. 

It follows from Lemma \ref{lemma:uniqueembedding} that such $\MM$ make up a single equivalence class. In particular, we may assume that $\MM=p\LL_3\operp q\LL_4$ is standard embedded in $\LL$ via its embedding 
in  $(p+q)\LL_4$. Then  the $\G$-stabilizer of $\MM$ contains 
$\Sfrak_p\times\Sfrak_q$ (which is very small, for $p+q\le 3$) as a permutation group of its summands and a lift of $\mu_{\MM}$ appears as a subgroup of 
the copy of $\mu_2^p$ defined by the antipodal  involutions  in  the $\LL_3$-summands.
\end{proof}

\subsection{Nonprimitively embedded geometric sublattices}
In order to classify the nonprimitive embeddings, we need to determine the overlattices of lattices 
$\MM$ of the form $k_\pt\LL_\pt$, i.e., 
embeddings $\MM\hookrightarrow \MM'$ with finite cokernel 
(whose order we call the \emph{degree} of this overlattice). 
Recall from Subsection \ref{subsect:eislattices} that  $ \MM'$ then embeds in $\MM^\vee$ 
and defines an isotropic submodule $I=\MM'/\MM\subset \MM^\vee/\MM=\disc(\MM)$ such that the  
discriminant module  of $\MM'$ can be identified with $I^\perp/I$. 
We say that two overlattices $j':\MM\hookrightarrow  \MM'$ and $j'': \MM\hookrightarrow \MM''$ 
are \emph{equivalent}  if there exists an isomorphism 
$\tau:  \MM'\cong  \MM''$ and a $\sigma\in \U(\MM)$ such that $j''\sigma=\tau j'$.

If we take $\MM=k_\pt\LL_\pt$, then $\MM^\vee/\MM\cong \FF_3^{k_1+k_3}\oplus \FF_4^{k_2}$. This splitting is canonical and  
so an isotropic  submodule $I$ will be  the form $I_3\oplus I_4\subset \FF_3^{k_1+k_3}\oplus \FF_4^{k_2}$. Hence the associated overlattice of $k_\pt\LL_\pt$ can be written accordingly as 
$\MM'_3\operp \MM'_4\operp k_4\LL_4$ with $\MM'_3$ an overlattice of $k_1\LL_1\operp k_3\LL_3$
and $\MM'_4$ one of $k_2\LL_2$. We therefore discuss the overlattices of $k_1\LL_1\operp k_3\LL_3$ and  $k_2\LL_2$ and assume that they are not already realized on a subsum.

We determine the overlattices for small values of $k_i$.

\begin{lemma}\label{lemma:1-3overlattice}
A strict overlattice  of $k_1\LL_1\operp k_3\LL_3 $ with $k_1+k_3\le 2$ only exists when $k_1=k_3=1$  
and is  then a copy of $\LL_4$ (and of degree $3$). 
All such form a single equivalence class. 

A strict overlattice  of $k_1\LL_1\operp k_3\LL_3 $ with $k_1+k_3=3$  
not already realized by a subsum only exists when $(k_1, k_3)$ equals $(3,0)$ 
(it is then a copy of $\LL_3$) or $(0,3)$. This overlattice is of degree $3$,
invariant under the permutation of the summands, 
and makes up a single equivalence class. Its discriminant form is  isomorphic with 
$x\in\FF_3\mapsto -x^2\in \FF_3$ resp.\ $x\in\FF_3\mapsto x^2\in \FF_3$.

If we denote the $\Sfrak_3$-invariant  overlattice of $3\LL_3$  by  $(3\LL_3)'$, 
then $\U((3\LL_3)')$ is isomorphic to $\mu_2\times (\Gg(\LL_3)\wr \Sfrak_3)$ 
(the $\mu_2$-factor acting as $\pm 1$).

\end{lemma}\begin{proof}
From the above, we see that the discriminant form of $k\LL_1\operp l\LL_3$ is
\[
\textstyle (x,y)\in \FF^3\oplus\FF^3\mapsto  \sum_{i=1}^k x_i^2-\sum_{j=1}^{l} y_j^2\in \FF^3
\]
 If  $(k,l)$ equals $(2,0)$ or $(0,2)$ this form does not represent zero. It does so when $k=l=1$: 
 just take $x_1$ and $y_1$ nonzero. The corresponding overlattice is then of degree 3 and a copy of 
 $\LL_4$. The group $\mu_2\wr\Sfrak_2$ permutes all the choices and so all these overlattices are equivalent. 

For $(k,l)$ equal to $(2,1)$ or $(1,2)$ a  nonzero isotropic vector must have one component zero. 
This means that the associated degree 3 overlattice is realized on a proper subsum and of the type just discussed.  

For $(k,l)=(3,0)$ the  nonzero isotropic vectors are obtained by taking all coordinates nonzero. 
Since the group $\mu_2\wr\Sfrak_3\subset \U(3\LL_1)$ permutes them transitively, the  associated 
overlattices make up one equivalence class. 

If we take the isotropic line spanned by $(1,1,1)$ we get the overlattice $\LL_3$ of $3\LL_1$. 
Note that the permutation group  $\Sfrak_3$  of the summands of $3\LL_1$ leaves this overlattice invariant.  

The same argument shows that the strict  overlattices of $3\LL_3$ make up one equivalence class 
(all of degree $3$). The one associated with the span of  $(1,1,1)$ is  an $\Sfrak_3$-invariant 
overlattice $(3\LL_3)'$ of degree 3. 

The rest of the proof is left to the reader.
\end{proof}

So $\U((3\LL_3)')$ is strictly contained in the unitary group of $3\LL_3$, which is  $\U(\LL_3)\wr\Sfrak_3=\big(\mu_2\times \Gg(\LL_3)\big)\wr \Sfrak_3$.

The overlattice $(3\LL_3)'$ is in a sense dual to  the overlattice $\LL_3$ we found for $3\LL_1$. This means that we should be able to primitively  embed  $(3\LL_3)'$ in a unimodular lattice  with orthogonal complement  isomorphic to $\LL_3$. This is indeed the case, but in a sense we can do better than that:

\begin{lemma}\label{lemma:3L3}
The lattice $(3\LL_3)'\operp \LL_1(-1)$  admits a unimodular overlattice. This overlattice is unique up to equivalence and  isomorphic with  $2\LL_4\operp\HH$. 
In particular, any unitary transformation of $(3\LL_3)'$ extends to one of $2\LL_4\operp\HH$.
\end{lemma}
\begin{proof}
The discriminant form of $(3\LL_3)'$ is isomorphic to $x\in \FF^3\mapsto x^2\in \FF_3$ and isomorphic to  the opposite of the one of $\LL_1(-1)$
(which is isomorphic to $x\in \FF^3\mapsto -x^2\in \FF_3$).  The direct sum has two isotropic one-dimensional subspaces, generates by $(1,1)$ and $(1,-1)$. By Remark \ref{rem:discriminant_duality} each of these determines  a unimodular overlattice of $(3\LL_3)'\oplus \LL_1(-1)$. They are obviously equivalent. The theorem of Basak (\cite{basakI}, Lemma 2.6) implies that  such a  unimodular 
lattice is isomorphic to  $2\LL_4\operp\HH$.
\end{proof}

\begin{lemma}\label{lemma:uniqueembedding3}
The  primitive sublattices in $\LL$ of type $(3\LL_3)'$ make up single equivalence class. 
\end{lemma}
\begin{proof}
It suffices to show that the orthogonal complement of such embedding contains a $(-3)$-vector.  
Lemma \ref{lemma:3L3} gives a primitive embedding of  $(3\LL_3)'$  in $\LL$ with orthogonal complement isomorphic to 
$\LL_1(-1)\operp \LL_4\operp\HH$. We must show that  every   primitive embedding  
$j: (3\LL_3)'\hookrightarrow \LL$ is $\G$ equivalent to this one. 
Nikulin's condition  is satisfied for the underlying integral lattice and so by Proposition \ref{prop:nikulin} the orthogonal complement of $j$ contains a vector $s'$ with $\la s',s'\ra=-3$ and is such that the orthogonal complement of $j((3\LL_3)')+\Ecal s'$ in $\LL$ is unimodular and generated by $3$-vectors.  This orthogonal complement has signature $(4,0)$ and is therefore a copy of $\LL_4$. So the primitive hull of 
 $j((3\LL_3)')+\Ecal s'$ is unimodular of signature $(9,1)$. 
 
This reduces the problem to showing that any two unimodular overlattices of $(3\LL_3)'\operp \LL_1(-1)$   are equivalent.
This follows from the fact that the discriminant form is given by $(x_1,x_2)\in \FF_3^2\mapsto x_1^2-x_2^2$: we see that it has only
two isotropic lines, spanned by $(1,1)$ and $(1,-1)$, which are exchanged by the reflection  of $(3\LL_3)'\operp \LL_1(-1)$ in the second summand.
\end{proof}

We next consider the overlattices of $4\LL_3$. 

\begin{lemma}\label{lemma:uniqueembedding1}
Every overlattice of $4\LL_3$ that is not already realized on a 
subsum is unimodular of degree 9 and has the property that the primitive hull of any subsum of type 
$3\LL_3$ in $(4\LL_3)'$   is an overlattice of degree $3$.  These overlattices make up a single equivalence class.
If we denote one such by $(4\LL_3)'$, then 
$\U((4\LL_3)')$  appears as  an extension of $\Sfrak_4$ by  $\mu_2\times \Gg(\LL_3)^4$, were $\mu_2$ acts by scalar multiplication and $\Sfrak_4$ is the  permutation group of the summands (involving certain signs). 

If we identify $(3\LL_3)'$ with the primitive hull of $3\LL_3\operp 0$ in $(4\LL_3)'$, 
then the $\U(4\LL_3)'$-stabilizer of that sublattice induces its full unitary group, realized here as $\mu_2\times (\Gg(\LL_3)\wr \Sfrak_3)$ 
(the $\mu_2$-factor acting as $\pm 1$).

\end{lemma}
\begin{proof} 
The discriminant form or $4\LL_3$ is $(y_1,y_2,y_3,y_4)\in \FF_3^4\mapsto -\sum_{i=1}^4 y_i^2\in \FF_3$ 
An isotropic  plane $I$ is given by  the span of $(1,1,1,0)$ and $(0,1,-1,1)$. 
In that case $I^\perp=I$ and hence the overlattice is unimodular. It is not hard to verify that $\mu_2^4\subset \U(4\LL_3)$ 
permutes these isotropic planes transitively with  kernel the diagonally embedded $\mu_2$ (if we divide out by this kernel, the action becomes simply transitive).  Each of the four coordinates planes $\FF_3^4$ meets $I$ in  a  line $\ell_k$ and 
therefore the primitive hull of the associated subsum of $4\LL_3$ in this overlattice is of degree $3$ and therefore isomorphic to $(3\LL_3)'$. 

It also follows that any  permutation $\sigma\in \Sfrak_4$ lifts to a unitary transformation  of $4\LL_3$  which permutes the summands accordingly and at the same time preserves $I$. Such a lift will then preserve $(4\LL_3)'$. If $\sigma$ is the  identity element, then this lift  preserves each line in $I$. It also acts linearly in $I$ and so must act in $I$ as $\pm 1$. It follows that this lift lies in $\mu_2\times \Gg(\LL_3)^4$.

The proof of the last paragraph is  left to the reader.
\end{proof}

\begin{remark}\label{lemma:4L3}
Since $(4\LL_3)'$ is unimodular, it has a  Niemeier lattice as its underlying integral lattice. Allcock \cite{allcock:y555} shows that  it appears as a direct summand in $\LL$ with orthogonal complement a copy of $\HH$ (\footnote{I was informed by a  referee of an earlier version that the last clause of Lemma \ref{lemma:uniqueembedding1} corrects an error in Thm.\ 4 of \cite{allcock:y555}.}).  In particular, we have a primitive embedding  $(4\LL_3)'\hookrightarrow\LL$ and all such sublattices lie in the same $\G$-orbit.
Via Lemma \ref{lemma:uniqueembedding1} this then also 
yields a primitive embedding $(3\LL_3)'\hookrightarrow\LL$. 
\end{remark}

\begin{corollary}\label{cor:uniqueembedding2}
Nonprimitive embeddings $p\LL_3\operp q\LL_4\hookrightarrow \LL$ exist if and only if $(p,q)$ is one of  $(3,0)$, $(3,1)$ or $(4,0)$. These cases make up a single $\G$-orbit if we also fix the degree of the primitive hull of the image. This degree  can be 3 or 9:  those of degree 3  are represented by primitive embeddings of $(3\LL_3)'$,  $(3\LL_3)'\operp\LL_3$  and $(3\LL_3)'\operp\LL_4$ in $\LL$ and 
those of degree 9  are defined by a primitive embedding of  $(4\LL_3)'$ in $\LL$. 
\end{corollary}
\begin{proof}
Since $\LL_4$ is unimodular, it  follows from Lemma \ref{lemma:uniqueembedding1} 
that the primitive hull of a  nonprimitive embedding of $p\LL_3\operp q\LL_4$ 
gives a primitive embedding of one of the lattices listed: $(3\LL_3)'$, $(3\LL_3)'\operp\LL_4$, 
$(3\LL_3)'\operp\LL_3$ or  $(4\LL_3)'$. So the lemma comes down to the statement that all such 
sublattices are $\G$-equivalent.
We already observed  this for $(4\LL_3)'$
and Lemma \ref{lemma:uniqueembedding3} shows that this also holds for  $(3\LL_3)'$. 
In this last case, the orthogonal complement is $\LL_4\operp\LL_1(-1)$. 
The standard embeddings of $\LL_3$  and $\LL_4$ in the first summand yields primitive embeddings 
of  $(3\LL_3)'\operp\LL_3$  and $(3\LL_3)'\operp\LL_4$ in $\LL$ with orthogonal complement 
isomorphic to $\LL_1\operp\LL_1(-1)$ resp.\ $\LL_1(-1)$. An application of Proposition \ref{prop:nikulin} to the underlying integral lattices gives the 
desired result.
\end{proof}

We will need (much) later also the following somewhat special case. 

\begin{proposition}\label{prop:lowdegreeclass}
Consider the  geometric lattices  $\MM\subset \LL$  of rank $\le 6$.
Then  $\MM$  is  primitively embedded unless $\MM$ has  a direct summand
of the type $\LL_3\operp \LL_1$, $3\LL_1$ or $4\LL_1$ whose primitive hull is of type
$\LL_4$, $\LL_3$ or $\LL_4$ respectively. Each of these cases make up a single $\U(\LL)$-orbit. 
\end{proposition}
\begin{proof}
This follows in a straightforward manner from Lemma \ref{lemma:uniqueembedding} and the preceding discussion.
\end{proof}

If we combine Proposition \ref{prop:uniqueembeddinghalf} with Corollary \ref{cor:uniqueembedding2}, we find:

\begin{corollary}\label{cor:geometricemb}
The $\G$-orbits of geometric  sublattices   $\MM\subset \LL$  are represented by the following 13 types (the decomposition is that of its primitive hull)
\begin{gather*}
\LL_3,  2\LL_3, 3\LL_3 , (3\LL_3)', (3\LL_3)' \operp\LL_3, (4\LL_3)' \\
\LL_4, \LL_3\operp \LL_4,  2\LL_3\operp \LL_4, (3\LL_3)'\operp \LL_4, 2\LL_4, \LL_3\operp 2\LL_4, 3\LL_4. 
\end{gather*}
For every geometric lattice $\MM$ we have an exact sequence
\[
\G_{[\MM]}=(\Gg(\MM)\times \G(\MM^\perp).  (\mu_\MM\rtimes \out_{\LL}(\MM)) 
\]
Here  $\mu_\MM\rtimes \out_{\LL}(\MM)$ acts faithfully on $\MM$ as  the group generated by the permutations up to sign of  the bundled $\LL_3$-summands (i.e., copies of $\LL_3$, $(3\LL_3)'$, $(4\LL_3)'$ and $\LL_4$).  Furthermore, the quotient 
$\G_{[\MM]}/(\Gg(\MM)\rtimes \out_{\LL}(\MM))$ (a not necessarily split extension $\G(\MM^\perp).\mu_\MM$) acts faithfully on $\MM^\perp$.
\hfill $\square$
\end{corollary}
\begin{proof}
Everything  indeed  follows from  Proposition \ref{prop:uniqueembeddinghalf} and  Corollary \ref{cor:uniqueembedding2}, except the last assertion. This follows from the fact that $\mu_\MM$ acts faithfully on the discriminant module of the primitive hull of $\MM$ and hence by Remark \ref{rem:discriminant_duality} also faithfully on the discriminant module of $\MM^\perp$.
\end{proof}

Since the outer automorphisms are transitive on the bundled $\LL_3$-summands, Corollary \ref{cor:geometricemb} also implies:

\begin{corollary}\label{cor:coarse-ct}
The $\G_\ct$-orbits of the geometric lattices isomorphic with $p\LL_3$ which strictly  contain a fixed $\LL^0_3$  are represented by the following six types
\begin{gather*}
\LL^0_3\operp\LL_3, \LL^0_3\operp 2\LL_3,  (\LL^0_3\operp2\LL_3)',  \LL^0_3\operp(3\LL_3)', (\LL^0_3+2\LL_3)'\operp \LL_3, (\LL^0_3\operp 3\LL_3)',\\
\LL^0_3\operp \LL_4, \LL^0_3\operp\LL_3\operp \LL_4, (\LL^0_3\operp 2\LL_3)'\operp \LL_4, \LL^0_3\operp 2\LL_4.
\end{gather*}
\hfill $\square$
\end{corollary}

Corollary \ref{cor:geometricemb} gives us the   $\G$-orbits of the  geometric sublattices of $\LL$. 
There are thirteen conjugacy classes of geometric subgroups of the Bimonster as well. The  way we have labeled them gives  an explicit  link with the Bimonster, for this defines a  bijection $\Phi$
between the conjugacy classes of the geometric subgroups of the Bimonster of type $(p,q)$ and the equivalence classes if geometric sublattices of the Allcock lattice of the same type.  We denote the simplicial complex defined by the partial ordered set $\Lscr$ of geometric sublattice of $\LL$ by $\Sigma(\Lscr)$. Mere inspection shows:
%We shall lift this to the level of complexes. The proof  invokes  Proposition \ref{prop:matching}.

\begin{corollary}[Bimonster matching]\label{cor:bimonstercorrespondence}
The map $\Phi$ defines  isomorphism of simplicial complexes $\Phi: \G\bs \Sigma(\Lscr)\cong  \Gcal\bs\Sigma(\Cscr)$. It preserves type and  bundling  (in particular takes the vertex  represented by the   $\LL_3$-sublattices resp.\  $\LL_4$-sublattices of the Allcock lattice  to the vertex represented by $\Cscr(A_3)$ resp.\  $\Cscr(A_4)$) and gives (therefore) conjugate (relative) outer automorphism groups. \hfill $\square$
\end{corollary}

By the last assertion  we mean that if  the geometric lattice $\MM\subset \LL$ and $W\in \Cscr$  correspond under $\Phi$, then
the set of the summands of $\MM$ (of type $\LL_3$ or $\LL_4$) with its $\out_{\LL}(\MM)$-action is isomorphic with 
the set of factors of $W$ (of type $A_3$ or $A_4$)  with its $\out_{\Gcal}(W)$-action.  We 
do however   \emph{not} assert here  that this isomorphism is canonical.

We need to refine  this correspondence in a manner which involves the groups and 
ball quotients themselves. This involves among other things that  we give  for  any  geometric sublattice $\MM$  the ball quotient 
$\G_{[\MM]}\bs\BB(\MM)$ a modular interpretation in terms of cubic threefolds or genus 5 curves.  These two are related by the fact that both have (intermediate) jacobians of dimension 5.  (This also explains why we called such $\MM$ geometric.) By imposing a  level two structures on them, we will be able to explain to some extent why it is that  we get the same bundling of summands resp.\ factors.
This however necessitates a considerably more detailed study of these moduli spaces  than can be found in the literature. The  next few sections are devoted to this.

\section{A regular  neighborhood of the moduli space of cubic threefolds}\label{sect:neighborhood}

We write $(\Xcal, D)$ for the Allcock ball quotient $\G\bs (\BB, Z)$. 
\\

We recall that $\LL^\ct=(\LL^0_3)^\perp$, that $\G_{[\LL^0_3]}\to \U(\LL^\ct)$ is onto with kernel $\Gg(\LL^0_3)$ and that 
the restriction map yields an isomorphism $\G_{\LL^0_3}\cong \G^\ct$. 
The identification $\LL^\ct=(\LL^0_3)^\perp$ makes $\BB^\ct$  a subball of $\BB$ of codimension 3. It is an intersection of $3$-vector mirrors of $\BB$ and hence $\BB^\ct\subset Z$.
Our focus will be on the union of its $\G$-translates
\[
\KK:=\cup_{g\in \G}\,  g\BB^\ct\subset \BB,   
\]
which is also contained in  $Z$ and thus defines  a closed subset of $D$ that  is the image of the natural map $\Xcal_\ct\to \Xcal$.
We shall therefore refer to the $K:=\G\bs\KK$ in $\Xcal$ as the \emph{cubic threefold locus} in $\Xcal$.  We state our main theorem: 

\begin{theorem}\label{thm:main}
There exists  regular neighborhood of the cubic threefold locus  in the Allcock ball quotient $\Xcal$  which admits a universal reflection cover whose Galois group is isomorphic with  the Bimonster $\Gcal$.
\end{theorem}

The proof will rest on  a modular interpretation of $\Xcal_\ct$, which we will use to prove (in later sections)

\begin{theorem}\label{thm:ctcover}
The pair $(\Xcal_\ct, D_\ct)$ admits a universal reflection cover 
\[
(\hat \Xcal_\ct, \hat D_\ct)\to (\Xcal_\ct, D_\ct)
\]
with Galois group $\Orth(\Tcal)$.
\end{theorem}

The plan is to derive  Theorem \ref{thm:main} from Theorem \ref{thm:ctcover}  as follows. The  regular neighborhood $U_K$ comes with a covering by open subsets  indexed by $\G\bs\Lscr$ whose nerve is $\G\bs\Sigma(\Lscr)$. Via  the Bimonster matching we identify this indexing set  with $\Gcal\bs\Cscr$ and its nerve with $\G\bs\Sigma(\Cscr)$.   
With the help of Theorem \ref{thm:ctcover}  we  produce an open  regular neighborhood $U_K$ of $K$ in $\Xcal$ and a reflection cover  $(\hat U_K, D(\hat U_K))\to (U_K, D(U_K))$ with Galois group $\Gcal$.  The construction endows 
$\hat U_K$ with covering by open subsets, but now indexed by $\Cscr$ and with nerve $\Sigma(\Cscr)$. We then prove that the hypotheses of  Corollary \ref{cor:simply connected cover2} are satisfied, so  that  $\hat U_K$ is simply connected and we can conclude that  $(\hat U_K, D(\hat U_K))\to U_K, D(U_K))$ is universal. 

In order to make this work, we need a somewhat more precise version of Theorem \ref{thm:ctcover}, which we will state at the appropriate moment.

\subsection{A stratification of the ball quotient} 
Lemma \ref{lemma:uniqueembedding}  tells us among other things that every sublattice of $\LL$ isomorphic to $\LL_3$ is a $\G$-translate of $\LL_3^0$
and so  in the union defining ${\KK}$, we may also  take the collection of such sublattices as index set.  Since this union is locally finite, it comes with an obvious $\G$-invariant  stratification. These strata will be open subsets of subballs that are intersections of $\G$-translates of $\BB^\ct$. In order to describe  these strata, we must know how 
$\BB^\ct$ intersects its $\G$-translates. This  means that we  consider intersections of $\G$-translates 
of $\LL^\ct$ with positive definite  orthogonal complement. Equivalently,  we consider positive definite  sublattices $\MM$ of $\LL $ that are generated by its sublattices of type $\LL_3$. Lemma \ref{lemma:posdef} tells us that such an $\MM$ is isomorphic to $p\LL_3\operp q\LL_4$ for certain nonnegative integers $p,q$ with $p+q\ge 1$. In other words, $\MM$ is a  geometric  sublattice of $\LL$. It is also clear that every geometric  sublattice so arises.

For  $z\in \BB$,  denote by $\ell_z\subset V$  the associated complex line and by $\MM_z$ the  span of the $\LL_3$-sublattices of $\ell ^\perp\cap \LL$.
Given a geometric lattice $\MM$, we define $\BB(\MM^\perp)'$  as the set  $z\in \BB$ with $\MM_z=\MM$.
It is clear that $\BB(\MM^\perp)'$ is open-dense  in $\BB(\MM^\perp)$ and that the   $\BB(\MM^\perp)'$  define the natural partition of ${\KK}$ into complex submanifolds. Indeed, $\BB(\MM^\perp)'=\BB(\MM^\perp)\ssm \cup_{\MM'} \BB(\MM'{}^\perp)$, where $\MM'$ runs over the geometric lattices that strictly contain $\MM$. 
The following is clear from the definition. 

\begin{lemma}\label{lemma:orbitintersection}
The collection $\{\BB(\MM^\perp)'\}_{\MM\in \Lscr}$ defines a $\G$-invariant decomposition of ${\KK}$.
In particular, if $g\in \G$ is such that  $g\BB(\MM^\perp)'\cap\BB(\MM^\perp)'$ is nonempty, then $g$ stabilizes $\MM$. 

A triflection mirror in $\BB$  which meets $\BB(\MM^\perp)'$ is defined by a triflection  in $\MM$ or  by a triflection in $\MM^\perp$. 

If $z\in \BB(\MM^\perp)'$, then the $3$-vectors perpendicular to $\MM\oplus \ell_z$ span a sublattice of type $k_1\LL_1\oplus k_2\LL_2$.
\hfill $\square$
\end{lemma}

The $\G$-stabilizer of $\MM$, $\G_{[\MM]}$,  acts on $\MM^\perp$ with finite kernel
and (hence) properly on $\BB(\MM^\perp)$. We put
\[
\Xcal_{\MM^\perp}:=\G_{[\MM]}\bs \BB(\MM^\perp), \quad \Xcal'_{\MM^\perp}:=\G_{[\MM]}\bs \BB(\MM^\perp)'\subset \Xcal_{\MM^\perp}.
\]
So this is a  ball quotient  and  a  (Zariski) open subset of it. 

\begin{corollary}\label{cor:coarse_emb}
The natural map $\Xcal_{\MM^\perp}\to \Xcal$ is finite and is an embedding when restricted to the orbifold $\Xcal'_{\MM^\perp}$.
The loci  $\Xcal'_{\MM^\perp}$ define a stratification (called the \emph{coarse stratification}) of $K=\G\bs {\KK}$ into connected orbifolds whose members are indexed by the $\G\bs\Lscr$ or equivalently (in view of the Bimonster matching \ref{cor:bimonstercorrespondence}), by  $\Gcal\bs\Cscr$. \hfill $\square$
\end{corollary}

Let $W_o\in \Cscr(A_3)$.  So its $\Gcal$-centralizer $W_o^\perp$ is isomorphic with $\Orth(\Tcal)$. Denote by $\Cscr_{\ge  W_o}$ the partially ordered set of  geometric subgroups of $\Gcal_o$ which contain $W_o$.
Every geometric subgroup $W\in \Cscr_{\ge  W_o}$   has 
\begin{enumerate}
\item [(a)] either $W_o$ as a direct factor: $W=W_o\times W'$ with $W'\subset W_o^\perp$,  
\item [(b)] or has a $\Sfrak_5$-factor $\tilde W_o>W_o$ and $W=\tilde W_o \times W'$ with  $W'\subset \tilde W_o^\perp$. 
\end{enumerate}
 In both cases, $W'$ is a geometric subgroup of $W_o^\perp\cong \Orth(\Tcal)$ or trivial, but in the second  case we have something extra, namely the subgroup  $\tilde W_o^\perp\subset W_o^\perp$ which  is a geometrically embedded $\Sfrak_{12}$ in $ W_o^\perp$  which contains $W'$. 

The following theorem is entirely about cubic threefold moduli. The role played by  $\Gcal$ is here only for indexing purposes and via  the Bimonster matching and could be eliminated, but it is the context and  the way we will use this theorem, that makes this a natural way to state it in this manner. This is the  refinement of Theorem \ref{thm:ctcover} we alluded to;  we will prove it in Section \ref{sect:coarsestrata}.

A \emph{coarse stratum of $\hat \Xcal_\ct$} is by definition a connected component of the preimage of a coarse stratum of  $\Xcal_\ct$.

\begin{theorem}\label{thm:match0}
Let $W_o\subset \Gcal$ be isomorphic to $\Sfrak_4$ and choose an isomorphism $W_o^\perp=C_{\Gcal}(W_o)\cong\Orth(\Tcal)$ to make  $W_o^\perp$ act on  $\hat \Xcal_\ct$. Let $W\in \Cscr_{\ge  W_o}$.
\begin{enumerate}
\item [(a)] If $W=W_o\times W'$ (so that $W'\le W_o^\perp$) then the fixed point set of $W'$ in $\hat \Xcal_\ct$ is irreducible and  contains a coarse stratum, denoted $\hat \Xcal\la W\ra$,   as an open-dense subset.
\item [(b)] If  $W=\tilde W_o\times W'$ with $\tilde W_o>W_o$ a $\Sfrak_5$-factor, then $\tilde W_o^\perp$ determines a connected component of the preimage of $\Xcal^{\DM}$ in  $\hat \Xcal_\ct$ and 
its intersection with the fixed point set of $W'$ in  $\hat \Xcal_\ct$  is irreducible and contains a coarse stratum, denoted  $\hat \Xcal\la W\ra$,   as an open-dense subset.
\end{enumerate}
The coarse strata $\hat \Xcal\la W\ra$ thus obtained satisfy the following properties:
\begin{enumerate}
\item [(i)] $\hat \Xcal\la W\ra $ comes with an action of $N_\Gcal(W)/W\cong W'{}^\perp.\out_\Gcal(W)$  such that its subgroup $N_{\Orth(\Tcal)}(W')/W'\cong W'{}^\perp.\out_{\Tcal}(W')$ is its  group of covering transformations (the orbit space is a coarse stratum of $\Mcal_\ct$).
\item [(ii)]  If the coarse stratum of $\Mcal_\ct$ associated with $W'$  is definable by the geometric 
sublattice $\LL^0_3 \subset \MM\subset \LL$, then $\Xcal'_{\MM^\perp}$ is canonically isomorphic with
the $N_\Gcal(W)/W$-orbit space of $\hat \Xcal\la W\ra $ (which we will denote by $\Xcal\la W\ra $).
  %$\out_{\LL}(\MM)$ acts on that stratum as the subgroup  $\out_{\Gcal}(W)$.
\item [(iii)]  If 
$W_1\ge W_o$ and  $W_2\ge W_o$ are in the same geometric conjugacy class of $\Gcal$, then  there is a canonical isomorphism $\Xcal\la W_2\ra \cong \Xcal\la W_2\ra $ which lifts to  an isomorphism 
 $\hat \Xcal\la W_1\ra \cong \hat \Xcal\la W_2\ra $.
 \end{enumerate}
\end{theorem}

\subsection{Regular neighborhoods}\label{sect:nbhds}
We want to extend the embedding of $\Xcal'_{\MM^\perp}\subset \Xcal$  in a  $\G_{[\MM]}$-equivariant manner to a regular neighborhood. 
We begin with making   a few general observations regarding  
tubular neighborhoods of subballs of complex balls. 

Suppose $V$ is a finite dimensional complex 
vector space with a nondegenerate hermitian form of Lorentzian signature and $V_0\subset V$ a 
positive definite subspace.
Then the orthogonal projection of $V$ onto ${V_0}^\perp$ determines the geodesic retraction 
$\pi_{V_0}: \BB(V)\to \BB({V_0}^\perp)$. If  a point of $\BB(V_0^\perp)$ is represented by the negative definite  line 
$\ell\subset {V_0}^\perp$, then 
the fiber of $\pi_{V_0}$ over $[\ell]\in \BB({V_0}^\perp)$ is naturally identified with the set of  
$\varphi\in \Hom_\CC(\ell,{V_0})$ for which 
$\la \varphi (z), \varphi(z)\ra<-\la z,z\ra$ for $z\in \ell$ (associate with $\varphi$ its graph). 
This is the unit ball for the natural inner product on $\Hom_\CC(\ell,{V_0})$. 
Another way to state this is by saying that $\BB(V)$ is naturally realized inside the 
${V_0}$-bundle over $\BB({V_0}^\perp)$ defined by $E_{\BB({V_0}^\perp)}\times_{\CC^\times} {V_0}$, 
where $E_{\BB({V_0}^\perp)}$ is the total space of the  line bundle $\Ocal_{\BB({V_0}^\perp)}(1)$.
This retraction is natural in the sense that if ${V_0}'\subset {V_0}$ is a subspace, then $\pi_{V_0}=\pi_{{V_0}|\BB({V_0}'{}^\perp)}\pi_{{V_0}'}$.  

Every  neighborhood of  $\BB({V_0}^\perp)$ in $\BB(V)$ contains   a neighborhood  whose  fibers  are 
open balls (of varying radius). Such a neighborhood is completely given by a  smooth function on $\BB({V_0}^\perp)$ taking values in the positive reals.
\emph{In this paper, tubular  neighborhoods of (open subsets of) subballs are always understood to be of this form: endowed with the 
geodesic retraction  map and  given a smooth radius function.}
\smallskip

%Later we shall want more freedom and allow the fibers to be  more generally geodesically convex  (in these balls). 
 After these generalities we return to  geometric sublattices of $\LL$. Put  
\begin{gather*}
 V(\MM):=\CC\otimes_\Ecal \MM, \\
 Z(\MM):=\text{the triflection mirror divisor of $U(\MM)$}\subset V(\MM).
\end{gather*}
We further denote by $Z(\MM^\perp)'$  the  triflection mirror divisor of $U(\MM^\perp)$ in $\BB(\MM^\perp)'$. If $\Ucal$ is an open subset of $\BB$, we may abbreviate $Z\cap \Ucal$ by $Z(\Ucal)$.

\begin{lemma}\label{lemma:regularnbhd1}
For every   geometric sublattice   $\MM\subset \LL$,  there exists a tubular neighborhood $\Ycal_\MM$ of $\BB(\MM^\perp)'$ in $\BB $ (in the above sense) such that 
\begin{enumerate}
\item[(i)] if $g\in \G$, then   $g\Ycal_\MM=\Ycal_{g\MM}$ (hence $\Ycal_\MM$ is  $\G_{[\MM]}$-invariant),
\item[(ii)] the pair $(\Ycal_\MM, Z(\Ycal_\MM))$ is (for the geodesic retraction) smoothly locally trivial over the pair 
$(\BB(\MM^\perp)',Z(\MM^\perp)')$  in the sense made precise below,   
\item[(iii)] if for some other geometric lattice   $\NN\subset \LL$,  the intersection $\Ycal_{\NN}\cap\Ycal_{\MM}$ is nonempty, then after perhaps exchanging the two, we have $\MM \supset \NN$,  and then the  geodesic retraction 
\[
(\Ycal_{\MM}\cap\Ycal _{\NN}, Z(\Ycal_{\MM}\cap\Ycal _{\NN}))\to \BB(\MM^\perp)'
\]
is locally trivial in the same sense as in  (ii).
\end{enumerate}
\end{lemma}

By local triviality in  property (ii),  we mean that  each $z\in \BB(\MM^\perp)'$ has a neighborhood  $U_z$ in $\BB(\MM^\perp)'$ invariant under $\G_{[\MM],z}$
such that there exists a $\G_{[\MM],z}$-equivariant  diffeomorphism of pairs:
\begin{center}
\begin{tikzcd}[column sep=small]
(\Ycal_\MM, Z(\Ycal_\MM))\big|_{U_z} \arrow[rr, "\cong"', "C^{\infty}"]&  &(V(\MM), Z(\MM))\times (U_z, U_z\cap Z(\MM^\perp)').
\end{tikzcd}
\end{center}
which commutes with the projection onto $U_z$. 
The meaning of local triviality in  property (iii) is similar.

\begin{proof}[Proof of Lemma \ref{lemma:regularnbhd1}]
We apply  the preceding discussion to  the situation where $V=V(\LL )$ and  $W=V(\MM)$.  
If we choose $\Ycal_\MM$ such that the points of $\Ycal_\MM$ are closer to $\BB(\MM)'$ than to 
$\BB(\NN)'$ for any geometric sublattice $\NN\not=\MM$, then not only (ii) holds, but by 
Lemma  \ref{lemma:orbitintersection}    (i)  is also satisfied. We can easily arrange that 
 (iii) is satisfied as well.  
 \end{proof}
 
We put  
\begin{equation}\label{eqn:UM}
\Ucal_\MM:=\cup_{\NN\subset \MM} \Ycal_\NN
\end{equation}
This is a  $\G_{[\MM]}$-invariant tubular neighborhood of $\BB(\MM^\perp)$ in $\BB$ (and comes a geodesic retraction on $\BB(\MM^\perp)$), but  its  intersection with  a $\G$-orbit in $\BB$  may consist of several 
$\G_{[\MM]}$-orbits. We put 
\begin{equation}\label{eqn:YM}
Y_{\MM}:=\G_{[\MM]}\bs \Ycal_\MM, \quad U_{\MM}:=\G_{[\MM]}\bs \Ucal_\MM.
\end{equation}

It follows from Lemma \ref{lemma:regularnbhd1}-(i) that we have  a geodesic retractions  $Y_{\MM}\to \Xcal'_\MM$ and $U_{\MM}\to \Xcal_\MM$ and that the natural map $Y_{\MM}\to \Xcal$
is an open embedding (thus realizing this  geodesic retraction inside $\Xcal$). We shall  identify  $Y_{\MM}$ with its image  in $\Xcal$. This gives an  open cover of $U_K$ whose nerve is $\Gcal\bs\Sigma(\Cscr)$.  This is no longer the case for the maps $U_\MM\to \Xcal$, for  there will be self-intersections. In fact, the image of the map $U_{\LL_0}\to \Xcal$ is all of  $U_K$. The preimage of $Y_\MM$ in $U_{\LL_0}$ is then a finite number  of copies of $U_{\MM}$: one for each $\out_{\LL}(\MM)$-orbits in the collection of $\LL_3$-summands of $\MM$.

\subsection{Reflection covers of regular neighborhoods}\label{subsect:nbhd}

We denote the  preimage of $D$ in $U_{\MM}$ by $D(U_{\MM})$.  This is also the $ \G_{[\MM]}$-orbit space of $Z(\Ucal_\MM)$.  
Let $Z_{\tfiber}(\Ucal_\MM)$  resp.\ $Z_{\tbase}(\Ucal_\MM)$ stand for the union of mirrors in $(\Ucal_\MM)$ defined by triflections in 
$\MM$ resp.\ $\MM^\perp$ and put $Z'(\Ucal_\MM):= Z_{\tfiber}(\Ucal_\MM)\cup Z_{\tbase}(\Ucal_\MM)$.
In general $Z'(\Ucal_\MM)$ is strictly contained in $Z(\Ucal_\MM)$, but Lemma \ref{lemma:orbitintersection} shows that 
$Z(\Ycal_\MM)= Z'(\Ucal_\MM)\cap \Ycal_\MM)$.  
We denote the images of these loci in  $U(\MM$ by replacing `Z' by `D', so we have 
$D'(U_\MM)= D_{\tfiber}(U_\MM)\cup D_{\tbase}(\U_\MM)$ and $D(Y_\MM)= D'(U_\MM)\cap Y_\MM)$.  
\medskip 

We abbreviate $\Ucal_\ct:=\Ucal_{\LL_3^0}$, $U_\ct$:=$U_{\LL^0_3}$, $V_\ct:=V(\LL_3^0)$  etc. (we had already agreed that  $\Xcal_\ct=\Xcal_{\LL^0_3}$).   

Recall that $V_\ct$ is such that the orbit space $\Gg(\LL_3^0)\bs V_\ct$ comes with 
a reflection covering $\hat V_\ct \to \Gg(\LL_3^0)\bs V_\ct$ with Galois group $\Sfrak_4$.
We shall now identify  $\Sfrak_4$ with a geometric subgroup $W_o$  of $\Gcal$ (so a member of $\Cscr(A_3)$.  Then   its $\Gcal$-centralizer $C_\Gcal(W_o)$ is isomorphic with $\Orth^-_{10}(\FF_2)$. Since  $(\LL_3^0)^\perp=\LL^\ct$, Theorem \ref{thm:ctcover} gives us a universal reflection cover $(\hat \Xcal_\ct,D(\hat \Xcal_\ct))\to ( \Xcal_\ct,D(\Xcal_\ct))$ with covering group 
$\Orth(\Tcal)$. We identify this group with $C_\Gcal(W_o)$. 
We can merge  the two constructions:

\begin{proposition}\label{prop:basicreflectioncovers1}
The universal reflection cover obtained in  Theorem \ref{thm:ctcover} lifts to a commutative diagram a universal  reflection covers 
\[
\begin{tikzcd}
(\hat U_\ct,D'(\hat U_\ct)) \arrow[d, "W_o\times C_\Gcal(W_o)"]\arrow{r} & (\hat \Xcal_\ct,D'(\hat \Xcal_\ct)) \arrow[d, "C_\Gcal(W_o)"] \\
(U_\ct,D(U_\ct)) \arrow{r}&  ( \Xcal_\ct,D(\Xcal_\ct)) \\
\end{tikzcd}
\]
\end{proposition}
\begin{proof}
Recall  that $\Ucal_\ct$ is a disk bundle  in the total space of the vector bundle 
$E_\ct:=E_{\BB(\LL_\ct)}\times_{\CC^\times} V_\ct\to \BB(\LL_\ct) $. This is acted on by $\G_{[\LL_\ct]}$, 
which by  Corollary \ref{cor:geometricemb} has the structure 
\[
\G_{[\LL_\ct]}=\mu_2\times \Gg(\LL_3^0)\times \G(\LL_\ct^\perp)\cong \Gg(\LL_3^0)\times \U(\LL_\ct^\perp).
\]
Since the central factor $\mu_2$ acts  by scalars on $\LL$ and hence trivially on both $\BB$ and $E_{\BB(\LL_\ct)}$, 
its action is through $\Gg(\LL_3^0)\times \G(\LL_\ct)\cong \Gg(\LL_3^0)\times \PU(\LL_\ct)$.

We first construct  a reflection cover in the fiber direction. For this we  transfer the $W_o$-cover $\hat V_\ct\to \Gg(\LL_3^0)\bs V_\ct$
 to the vector  bundle $E_\ct$  and get 
a reflection cover
\[
E_{\BB(\LL_\ct)}\times_{\CC^\times} \hat V(\MM)\to \BB(\LL_\ct), 
\]
still with Galois group $W_o$. That  action   commutes with the  one of $\G(\LL_\ct)$. 
We divide out by the action of the latter
and thus find a bundle $\hat E_{\Xcal_\ct}\to \Xcal_\ct$ with fiber $\hat V_\ct$ and fiber preserving  action of $W_o$. 
Its  pull-back along the universal reflection cover $\hat \Xcal_\ct\to \Xcal_\ct$ with Galois  group $C_\Gcal(W_o)$ gives a 
vector bundle $\hat E_{\hat \Xcal_\ct}\to \hat \Xcal_\ct$ with $W_o\times W^\perp_\ct$-action. The  construction makes it clear that 
the $W_o\times W^\perp_\ct$-orbit space of $\hat E_{\hat \Xcal_\ct}$ is naturally identified with the $\G_{[\LL_3^0]}$-orbit space of $E_\ct$ and the resulting 
map is a reflection cover. Since  $\hat \Xcal_\ct$ is simply connected, so is the vector bundle $\hat E_{\hat \Xcal_\ct}$ over it. 

The disk subbundle $\Ucal_\ct$ of $E_\ct$ determines a disk subbundle  $\hat U_\ct$ of $\hat E_{\hat \Xcal_\ct}$. This gives  a reflection cover 
$(\hat U_\ct,\hat D'_\ct)\to (U_\ct, D'_\ct)$ with a simply connected total space. Since there is ramification over $D'_\ct$, this reflection cover is universal.
\end{proof}

We express Proposition \ref{prop:basicreflectioncovers1} in notation related to the Bimonster and write 
\[
(\hat U\la W_o\ra, D'(\hat U\la W_o\ra))\xrightarrow{N_{\Gcal}(W_o)} (U\la W_o\ra, D'(U\la W_o\ra)), 
\]
for $(\hat U_\ct,D'(\hat U_\ct)) \to (U_\ct,D(U_\ct))$,  where we note that $N_{\Gcal}(W_o)=W_o\times C_\Gcal(W_o)$.

The space $Y_\MM$ space only depends on the $\G$-orbit of $\MM$ in $\LL$, which by the Bimonster matching corresponds to a geometric conjugacy class in $\Cscr$. We are therefore justified to denote it by
$Y\la W\ra$ instead,  where $W$ represents that conjugacy class.

\begin{proposition}\label{prop:basicreflectioncovers}
For  every $W\in \Cscr$, there is a (connected)  reflection cover 
\[
 (\hat Y\la W\ra ,D(\hat Y\la W\ra )) \to  (Y\la W\ra,D(Y\la W\ra))
\]
with  Galois group $N_{\Gcal}(W)$. It is unique up to a covering transformation. 

If  $W\in \Cscr_{\ge  W_o}$, then 
$\hat Y\la W\ra$ comes with an open embedding  in  $\hat U\la W_o\ra$ which is equivariant for $N_{\Gcal}(W_o\le W)=N_{\Gcal}(W_o)\cap N_{\Gcal}(W)$. This embedding has the property that if  $W, W'\in \Cscr_{\ge  W_o}$, then 
$\hat Y\la W\ra$ and $\hat Y\la W'\ra$ meet in $\hat U\la W_o\ra$ if and only if $W$ and $W'$ satisfy 
an inclusion relation.  
\end{proposition}
\begin{proof}
The proof is similar to the one of  Proposition \ref{prop:basicreflectioncovers1}. 
Let $\MM\subset \LL$ be a geometric  sublattice which matches $W$ under the Bimonster matching (so that in particular,  $W\cong W(\MM)$). 

Recall  that $\Ycal_\MM$ is a disk bundle  in the total space of the vector bundle 
$E_{\BB(\MM^\perp)'}\times_{\CC^\times} V(\MM)$. The  loci $Z_{\tbase}(\Ycal_\MM)$  and $Z_{\tfiber}(\Ycal_\MM)$ are the intersections with
$E_{Z(\MM^\perp)'}\times_{\CC^\times} V(\MM)$ and $E_{\BB(\MM^\perp)'}\times_{\CC^\times} Z(\MM)$, where the meaning of the notation implicitly introduced here should be evident. The inclusion is such that the universal reflection cover is going to be the restriction of  the one associated with the
$\G_{[\MM]}$-orbit space of this bundle

According to  Corollary \ref{cor:geometricemb},  the group $\G_{[\MM]}$ has the structure 
\[
\G_{[\MM]}= \big(\Gg(\MM)\times \G(\MM^\perp)\big). (\mu_\MM\rtimes \out_{\LL}(\MM)),
\]
where $\mu_\MM$ is a commutative group generated by involutions which preserve each factor of $\MM$.
The first step amounts to  constructing a reflection cover in the fiber direction which exhibits  all the symmetry of the situation. This is not difficult: the group $\G_{[\MM]}$ contains  $\Gg(\MM)$ 
as a normal subgroup and acts on $\MM$ through $(\mu_\MM\times \Gg(\MM))\rtimes \out_{\LL}(\MM)$.
By Corollary \ref{cor:isodiscr} we can find  a universal reflection cover 
\[
(\hat V(\MM), \hat Z(\MM))\to \Gg(\MM)\bs ( V(\MM), Z(\MM))
\]
with Galois group $W:=W(\MM)$, but for which the $W$-action extends to an action of 
$W. (\mu_\MM\rtimes \out_{\LL}(\MM))$  on the whole reflection cover.
We transfer this to the vector  bundle $E_{\BB(\MM^\perp)'}\times_{\CC^\times} V(\MM)\to \BB(\MM^\perp)'$ and get 
a reflection cover
\[
(E_{\BB(\MM^\perp)'}\times_{\CC^\times} \hat V(\MM), E_{\BB(\MM^\perp)'}\times_{\CC^\times} \hat Z(\MM)) \to \BB(\MM^\perp)', 
\]
still with Galois group $W$, but for which the $W$-action now extends to one of $\big(W(\MM)\times \G(\MM^\perp)\big). (\mu_\MM\rtimes \out_{\LL}(\MM))$
on the whole bundle. The fiberwise transform of $\Ycal_\MM$ under the finite maps 
\[
\hat V(\MM)\xrightarrow{W}\Gg(\MM)\bs V(\MM)\xleftarrow{\Gg(\MM)} V(\MM)
\]
 yields a subbundle $\Ycal_{\MM,\tfiber}\subset E_{\BB(\MM^\perp)'}\times_{\CC^\times} \hat V(\MM)$ over $\BB(\MM^\perp)'$, 
which is  invariant under  this group. The
$\G(\MM^\perp).\mu_\MM$-orbit space of $\Ycal_{\MM,\tfiber}$ yields  a holomorphic bundle pair
\begin{equation}\label{eqn:fibercover}
(Y_{\MM, \tfiber}, Z_\tfiber(\hat Y_{\MM, \tfiber})) \to \Xcal'_{\MM^\perp}.
\end{equation}
%It  comes with an action of $W\rtimes\out_{\LL}(\MM)$, which acts on the base $\Xcal'_{\MM^\perp}$ through $\out_{\LL}(\MM)$.  
By Theorem \ref{thm:match0}-(ii), we  
can identify $ \Xcal'_{\MM^\perp}$ with $\Xcal\la W\ra$.
%, then by Theorem \ref{thm:match0}-(ii) the $\out_{\LL}(\MM)$-action corresponds to an action  of 
%$\out_{\Gcal}(W)$. We make that identification and
Let   $\hat Y\la W\ra$  be the pull back  of 
$Y_{\MM, \tfiber}\to \Xcal'_{\MM^\perp}\cong \Xcal\la W\ra$
along the $N_\Gcal(W)/W$-cover  $\hat \Xcal\la W\ra\to \Xcal\la W\ra$. Then $\hat Y\la W\ra$  comes with a faithful  action of  $N_{\Gcal}(W)$ whose orbit space is $Y\la W\ra\cong Y_\MM$.

In case $\MM\supset \LL^0_3$, then $\Ycal_\MM$ is an open subset of $\Ucal_\ct$ and such subsets cover  $\Ucal_\ct$.  If $\MM'\supset  \LL^0_3$ and  $\Ycal_\MM\cap \Ycal_{\MM'}\not=\emptyset$, then by Lemma 
\ref{lemma:regularnbhd1}-iii,  we have  $\MM\subset \MM'$  or $\MM\supset \MM'$. 
The above construction realizes $Y\la W\ra$ as an open subset of $U_\ct$ and thus the other two assertions follow.
\end{proof}

Since $Y\la W\ra$ only depends on the conjugacy class of $W$, we can (and will) make our choices more carefully, by choosing
for each conjugacy class one representative $W_o$ (which we could take to contain $W_o$, but that is not needed here)
and then form 
\[
\Gcal \times _{N_{\Gcal}(W_o)} \hat Y\la W_o\ra
\]
so that is the quotient of $\Gcal \times  \hat Y \la W_o \ra$ by the action of $N_{\Gcal}(W_o)$ defined by $h(g,y):=(gh^{-1}, h(y))$.  Its connected components are naturally indexed by the conjugates $W$ of $W_o$ in $\Gcal$ and we may take  for 
$\hat Y\la W\ra$  the associated connected component.  These new choices differ from the previous ones essentially only 
by a covering transformation. It is however still true that when $W_o <W$, we have an open embedding 
$\hat Y\la W\ra\hookrightarrow \hat U_\ct$ which is equivariant with respect to $N_{\Gcal}(W_o\le W)$.
We can now prove our main theorem.

\begin{proof}[Proof of Theorem \ref{thm:main}] 
We apply the same procedure to $\hat U_\ct \la W_o\ra$: the 
 connected components of  $\Gcal \times _{N_{\Gcal}(W_o)} \hat U_\ct \la W_o\ra$ are indexed by the 
members of $\Cscr(A_3)$. For every $W_o\in \Cscr(A_3)$ we denote by the associated connected component by
$\hat U\la W_o\ra$. It comes with a $N_{\Gcal}(W_o)$-covering $\hat U\la W_o\ra\to U\la W_o\ra$. It has the property that whenever $W_o\le W$ we have an open $N_{\Gcal}(W_o\le W)$-equivariant embedding $\hat Y\la W\ra \hookrightarrow \hat U\la W_o\ra$. 
If $W_o\le W'$ and two such embeddings meet, then $W'\le W$ or $W'\ge W$. So if $\Wcal \subset \Cscr_{\ge  W_o}$ is a subset for which the images of the embeddings $\{\hat Y\la W\ra \hookrightarrow \hat U\la W_o\ra\}_{W\in \Wcal}$ have a point in common, then  $\Wcal$ defines a chain in $\Cscr$ and hence a simplex of $\Sigma(\Cscr)$.
We use these embeddings to make identifications in $\hat U_\ct \la W_o\ra$ and denote the resulting pair  
$(\hat U, D(\hat U))$. By construction it comes with an action of $\Gcal$ and a $\Gcal$-equivariant covering $\{\hat Y\la W\ra\}_{W\in \Cscr}$ whose nerve is $\Sigma(\Cscr)$. The $\Gcal$-orbit space is $(U_K, D(U_K)$ and the resulting map 
$(\hat U, D(\hat U))\to  (U, D(U))$ is a reflection covering.

The hypotheses of Corollary \ref{cor:simply connected cover2} are satisfied and hence $\hat U\la W\ra$ is simply connected. So this reflection  cover is universal.
\end{proof}

The rest of this paper is mainly devoted to the proof of Theorems \ref{thm:ctcover} and \ref{thm:match0}. We obtain some other properties on cubic threefold moduli along the way.

\newpage
\begin{center}
 \textbf{PART II: THE FINE STRUCTURE OF THE MODULI SPACE OF CUBIC THREEFOLDS}
 \end{center}

\section{The moduli space of cubic threefolds as a ten-ball quotient}\label{sect:cubic3fold}
The basic references for this section are Clemens-Griffiths \cite{cg}, 
Collino \cite{collino}, Allcock-Carlson-Toledo \cite{act}, Looijenga-Swierstra \cite{ls}, 
Casalaina-Laza \cite{cl} and Casalaina-Grushevsky-Hulek-Laza \cite{cghl}.

\subsection{The Torelli theorem for cubic threefolds}
We find it convenient  to think of cubic threefolds not as sitting in 
$\PP^4$, but in the projectivization $\PP V$ of a given  
complex vector $V$ of dimension five.  Then the cubic threefolds in $\PP V$ are parametrized by 
the projectivization of the space $\Kcal:=\CC[V]_3$ of cubic forms on $V$. Note that $\Kcal$ is as an $\SL(V)$-representation the dual of $\sym^3V$ and hence irreducible.

Let $X\subset \PP V$ be a reduced cubic threefold.  If $F\in \Kcal$ is a   defining equation,   
then let  $Y\subset \PP(V\oplus \CC)$  be the cubic 
$4$-fold defined by $\tilde F:=F -T^3$. We think of $Y$ as the  $\mu_3$-cover $\pi: Y\to \PP V$ 
which ramifies over $X$ 
with the $\mu_3$-action on $V\oplus \CC$ defined by $\rho(\omega): (v,t)\mapsto (v,\omega t)$. The $\mu_3$-action on $Y$ has $X$ as its fixed point locus.  Clearly, $F$ is determined by $X$ up to a scalar and that  ambiguity is reflected by the fact that $Y$ is only defined up to 
a covering transformation.

Assume that $X$ is smooth, so that $Y$ is also smooth. If  $\eta_Y \in \Hl^2(Y)$ denotes the hyperplane class, then
$\eta_Y^4[Y]=3$, or equivalently,  $\eta^2_Y\in  \Hl^4(Y)$ has self-intersection $3$.
The lattice $\Hl^4(Y)$ is unimodular (by Poincar\'e duality)  and odd (for $\eta^2_Y\cdot \eta^2_Y=3$) and its  signature is $(21,2)$. These properties  determine its isomorphism type: there exists a basis on which the intersection matrix takes a diagonal form with $\pm 1$ on the diagonal. We often identify $\Hl^4(Y)$ with $\Hl_4(Y)$ by means of the intersection pairing.

The \emph{primitive cohomology} $\Hl^4_\circ(Y)$ of $Y$ is  the kernel  of $\eta_Y^2\cup: \Hl^4(Y)\to \Hl^8(Y)\cong \ZZ$. 
This is also  the orthogonal complement of $\eta_Y^2$ in $\Hl^4(Y)$ for the intersection form and  
Poincar\'e  duality  $\Hl^4(Y)\cong \Hl_4(Y)$ turns   this is into the kernel of the natural map $\Hl_4(Y)\to \Hl_4(\PP V)$. It has  signature $(20, 2)$ and  discriminant $-3$ and comes with a distinguished set of elements of self-intersection 2 (the Poincar\'e duals of vanishing cycles) that generate it. In particular,   $\Hl^4_\circ(Y)$ is  even.  These two  properties  $\eta^2_Y\cdot \eta^2_Y=3$ and  $\eta^2_Y$ has even orthogonal complement characterize the pair $(\Hl^4(Y), \eta_Y^2)$ up to isomorphism.

The $\mu_3$-action on  $\Hl^4_\circ(Y)$ has no eigenvalue $1$ and hence this action makes $\Hl^4_\circ(Y)$ a torsion free $\Ecal$-module of rank $11$ (which is then necessarily free). 
The intersection pairing  is $\rho$-invariant, so that  $\Hl^4_\circ(Y)$ comes with a hermitian  form $\la\; ,\; \ra$ of (Lorentzian) signature $(10,1)$.
It was established in  \cite{act} and \cite{ls} that this Eisenstein lattice is isomorphic with $\LL^{\ct}=\LL_1\operp 2\LL_4\operp\HH$.  
We  recall from Subsection \ref{subsect:ballmodular} that the discriminant module  of $\LL^{\ct}$ is cyclic of order $3$ with  $-1\in \U(\LL^\ct)$ inducing its nontrivial automorphism  and that $\U(\LL^{\ct})=\{\pm 1\} \times \G^{\ct}$, where 
$\G^{\ct}\subset \U(\LL^{\ct})$ is the subgroup of unitary transformation that act trivially  on the discriminant module.

The direct summand  of $\LL_1\subset \LL^\ct$ and the inclusion $\Hl^4_\circ(Y)\subset \Hl^4(Y)$ each single out a generator of their respective discriminant modules. We refer to these as the \emph{prefered generators}.

\begin{definition}\label{def:integralctmarking}
An \emph{$\LL^\ct$-marking} of a smooth cubic threefold $X$ is the choice of an isomorphism of Eisenstein lattices $\Hl^4_\circ(Y)\cong \LL^{\ct}$ 
which takes one  prefered generator to the other.
\end{definition}

Such markings exist, for if  $\varphi: \Hl^4_\circ(Y)\xrightarrow{\cong} \LL^{\ct}$ is  any isomorphism of Eisenstein lattices, then precisely 
one of $\varphi$ and $-\varphi$ will be take a prefered generator to a prefered generator. Via such a marking the  group $\G^{\ct}$ is identified with the group of orthogonal transformations of  $\Hl^4(Y)$ that fix $\eta_Y^2$ and centralize the $\mu_3$-action. 
It is clear that two $\LL^\ct$-markings will  lie in the same $\G^{\ct}$-orbit. There is a bit of an issue though, because $X$ determines $Y$ up to the group $\mu_3$ of covering transformations  and that group acts faithfully on $\Hl^4_\circ(Y)$ as a group of scalars and hence the marking should also be given up a scalar in $\mu_3$.  Fortunately this hardly plays a role in what follows, for eventually  it is the projectivized group $\PG^{\ct}=\G^{\ct}/\mu_3\cong \PU(\LL^\ct)$ which matters. 
\\

The polarized Hodge structure on  $\Hl^4_\circ(Y)$ has Hodge numbers $h^{3,1}=h^{1,3}=1$ and $h^{2,2}=20$. 
Following Griffiths \cite{griffiths}, this   structure can be concretely described as follows. 
Fix an isomorphism  $\mu:\wedge^5 V\cong \CC$ and denote by $d\mu$  the associated translation invariant $5$-form on $V$. 
A generator of $\Hl^{3,1}_\circ(Y)=\Hl^{3,1}(Y)$ is then represented 
by the double residue of the  meromorphic  $6$-form $\Omega_{\tilde F}:=\tilde F^{-2}d\mu \wedge dt$ 
on $V\oplus \CC$. This form is homogeneous of degree zero and  so its  
residue at the hyperplane at infinity gives a nowhere zero $5$-form on $\PP (V\oplus\CC)$ 
with a double pole along $Y$.  
The  (Griffiths) residue of the latter on $Y$ can be understood as an element of 
$\Hl^1(Y, \Omega^3_Y)=\Hl^{3,1}_\circ(Y)$. Note that the $\mu_3$-action on $Y$ is such that
$\rho(\omega)$ multiplies this generator with $\omega$. 
This produces an $\SL(V)$-invariant  identification of the line spanned by 
$F^{-2}$ and $\Hl^{3,1}_\circ(Y)$. (We shall here think of $F^{-2}$ as a linear form on  $\sym^2 \Kcal$, 
that is, as an element of $\sym^2 (\sym^3V)$, but  if we wish to kill the $\mu_3$-ambiguity, 
then we should take $F^{-6}\in \sym^6 (\sym^3V)$ instead.)  
The Hodge index theorem implies that the intersection pairing is negative on $\Hl^{3,1}_\circ(Y)$. 
So the images of $\Hl^{3,1}_\circ(Y)$ under the (complexified) $\LL^\ct$-markings make up a  
$\PG^\ct$-orbit in $\BB^\ct$. 

According to Allcock \cite{allcock:3GIT} and Yokoyama  \cite{yokoyama:3GIT}, 
the stable cubic $3$-folds are those that admit singularities of type $A_4$ at worst. 
Let us denote by $\Mcal_\ct^\st$ the corresponding GIT quotient and by  
$\Mcal_\ct\subset  \Mcal_\ct^\st$ the  open subset which parametrizes smooth cubic $3$-folds. 
Allcock-Carlson-Toledo \cite{act} and Looijenga-Swierstra \cite{ls} independently proved:

\begin{theorem}\label{thm:cubicperiods}
The map which assigns to a stable cubic $3$-fold 
$X\subset\PP V$ the Hodge structure of the $\mu_3$-cover of $\PP V$ 
ramified along $X$, identifies  $\Mcal_\ct^\st$ with the complement of the special mirror divisor 
$X_\ct\ssm D^{\ct, \spl}$,  
and the locus of singular stable cubic $3$-folds (over whose generic point we have a 
cubic $3$-fold with an ordinary double point as its singular set) with  its intersection with the 
standard mirror divisor  $D^{\ct, \std}$.
\end{theorem}

\subsection{Stratification of the discriminant}\label{subsect:discrstrat}
 We  make the isomorphism appearing in Theorem \ref{thm:cubicperiods} as one of stratified spaces. In what follows, $\Kcal_\sing$ stands for the locus parametrizing singular cubic $3$-folds, 
$o$ for  some point of $\PP \Kcal\ssm  \PP \Kcal_\sing$ such that  the smooth  cubic $3$-fold  $X_o$ it defines has no nontrivial automorphism.
The associated cubic $4$-fold is denoted $Y_o$ and we fix an isomorphism $H^4_\circ(Y)\cong \LL^\ct$ of Eisenstein lattices. We shall need a bit of singularity theory, for which we shall use \cite{looij:icis} as our general reference.

\begin{definition}\label{def:degeneration}
Suppose $X$  is a stable  cubic threefold. By a \emph{degeneration of $X_o$ into $X$} we shall mean that we are  given 
an arc  $\g: [0,1]\to \PP \Kcal$  from $o$ to the point defining $X$ which maps $[0,1)$ to $ \PP \Kcal\ssm  \PP \Kcal_\sing$,  or rather a connected component in the space of such arcs.
\end{definition} 

With such a degeneration  is  associated  a \emph{vanishing manifold}: a 
compact submanifold with boundary $F_\g\subset X_o$ whose interior is open in $X_o$, given up isotopy such that the 
pull-back along $\g$ gives a family of threefolds over 
$[0,1]$, which  is topologically obtained from $[0,1]\times X_o\to [0,1]$ by contracting in $\{1\}\times X_o$ each connected component of $F_\g$ to a point. Note that this identifies $X_o\ssm F_\g$ with the smooth part $X_\reg$ of $X$. The connected component of $F_\g$ that maps to a given singular point $p\in X$ is called the \emph{Milnor fiber} of $p$. It is known that the reduced homology of  a Milnor fiber is concentrated in the middle dimension
(so $3$ in this case).  
The image of the map $\Hl_3(F_\g)\to \Hl_3(X_o)$ is called the \emph{vanishing homology} of the degeneration.  This map respects the intersection pairing, but since $\Hl_3(F_\g)$ might be 
degenerate, this map need not be injective.

The  situation is similar for the associated family of cubic fourfolds. Then $\g$ determines a $\mu_3$-invariant vanishing submanifold $\tilde F_\g \subset Y_o$ (unique up to  $\mu_3$-equivariant isotopy). The intersection pairing and the  $\mu_3$-action give $\Hl_4(\tilde F_\g)$ the structure of an Eisenstein lattice.  As we will see below this lattice is positive definite, so that $\Hl_4(\tilde F_\g)\to \Hl_4(Y_o)$  will always be an embedding of Eisenstein lattices. We shall also see that its image need not be primitive. Another choice of $\g$ has the effect of changing this  embedding by an element of $\G^\ct$ and so without that choice only the  $\G^\ct$-orbit of that positive definite sublattice is well-defined.

To make this concrete, assume for a moment that the singular set of $X$ consists of a single $A_i$-singularity,  so with a local analytic equation of the type 
\[
z_1^{i+1} +z_2^2 +z_3^2+z_4^2=0, \quad (i=1,2,3,4).
\]
Then the interior of $F_\g$ is diffeomorphic with the affine threefold $X_f$ defined by $f(z_1) +z_2^2+z_3^2+z_4^2=0$, 
where $f$ is any separable monic polynomial of degree $i+1$. If $C_f$ denotes the affine hyperelliptic plane curve  defined by  $f(z_1) +z_2^2=0$, then $\Hl_1(C_f)$ is a symplectic lattice of type $A_i$ and we have a natural identification $\Hl_1(C_f)\cong \Hl_3(X_f)\cong \Hl_3(F_\g)$ as symplectic lattices, see \cite{looij:icis}.
(One can show that in this case the map $\Hl_3(F_\g)\to \Hl_3(X_o)$ is injective, so that its image (the vanishing homology) is a copy of the symplectic $A_i$-lattice.) 
The  cubic fourfold $Y$  associated with $X$ has a singularity with a local analytic equation given by  
\[
z_0^3+z_1^{i+1} +z_2^2 +z_3^2+z_4^2=0,
\]
with $\mu_3$ acting by scalar multiplication on the $z_0$-coordinate.
 For $i=1,2,3,4$, this is a simple singularity of type $A_2,D_4,E_6, E_8$ respectively. 
Singularity theory tells us  that this defines a vanishing lattice in $\Hl_4(Y_o)$ (where $Y_o$ is the smooth 
 cubic fourfold associated with $X_o$) that is $\mu_3$-equivariantly isomorphic to $\Hl_4(Y_f)$, 
 where $Y_f$  is the affine fourfold $Y_f$ defined by  $z_0^3=f(z_1) +z_2^2+z_3^2+z_4^2$. 
 This  vanishing lattice is isomorphic with a root lattice of  the above type (see for example  \cite{looij:icis}). 
 The $\mu_3$-invariants in $\Hl_4(Y_f)$ are reduced to  $\{0\}$ and so the 
 $\mu_3$-action turns this vanishing lattice into an Eisenstein lattice of type  $\LL_i$. 
 Theorem \ref{thm:cubicperiods} implies that all the $3$-vectors in such lattice are standard
(this is of course only a nontrivial statement for $i=1$), for otherwise $X$ 
would map to a point of  $\Xcal^{\ct, \std}$ that lies in the image of special mirror. 
The contrast with the symplectic situation is that for $i\le 4$,  
$\LL_i$ is nondegenerate (it is positive definite), whereas the symplectic lattices of type $A_1$ 
and $A_3$ are degenerate. So if we do not make  a choice for $\g$, then we merely end up with
a $\G^\ct$-orbit of sublattices of $\LL^\ct$ isomorphic with $\LL_i$.

In  general  a stable cubic threefold with a singularity configuration of type 
$k_\pt A_\pt=k_1A_1+k_2A_2+k_3A_3+k_4A_4$ (meaning that it has exactly $k_i$ singularities of type $A_i$) determines a $\G^\ct$-orbit of sublattices of $\LL^\ct$ isomorphic with 
\[
k_\pt\LL_\pt:=k_1\LL_1\operp k_2\LL_2\operp k_3\LL_3\operp k_4\LL_4.
\]

\begin{proposition}\label{prop:stratabijection}
Let $\Scal(k_\pt A_\pt)\subset \Xcal_\ct\ssm D^{\ct, \spl}$ be the locus represented by the points $z\in \BB^\ct\ssm Z^{\ct, \spl}$ with the property that the $3$-vectors perpendicular to $z$ span  sublattice of type $k_\pt\LL_\pt$ (such $3$-vectors are necessarily standard) and  let  $\Mcal_\ct(k_\pt A_\pt)$ stand for the locus in  $\Mcal^\st_\ct$ parametrizing  cubic 
threefolds with a singularity configuration of type  $k_\pt A_\pt$.
Then the  isomorphism of Theorem \ref{thm:cubicperiods} maps $\Mcal_\ct(k_\pt A_\pt)$ isomorphically onto $\Scal(k_\pt A_\pt)$.
\end{proposition}
\begin{proof}
Let $z\in \BB^\ct\ssm Z^{\ct, \spl}$ and assume  that the $3$-vectors perpendicular to $z$ span  sublattice 
$\MM\subset \LL^\ct$ of type $k_\pt\LL_\pt$ and let $X\subset \PP(V)$ represent $z$. 
Choose a $\aut(X)$-invariant slice germ $(S,o)\subset \PP(\CC[V]_3)$  to the $\PSL(V)$-orbit of $[X]$ in
$\PP(\CC[V]_3)$. This defines  (by definition) a  universal  deformation  $f:(\Xscr,X)\to (S,o)$  of $X$.
The associated map germ  $(S,o)\to X_\ct$ is an orbifold chart for  $\Xcal_\ct$: the natural map 
$(S,o)\to \Gg(\MM)\bs \BB^\ct$ is local isomorphism and  $\G^\ct_z/\Gg(\MM)\cong \aut(X)$.
By Theorem \ref{thm:cubicperiods} the discriminant $D(f)$ of $f$ is the preimage of the standard mirror divisor 
$D^{\ct, \std}$. So the  irreducible components of $D(f)$  meet transversally and are naturally 
indexed by the summands of $\MM$. It remains to show that if an  irreducible component of $D(f)$ is associated 
with  a summand of $\MM$ is of type $\LL_i$, then it  comes from  a unique singular point on $X$ of  type  $A_i$
and that all singularities of $X$ are thus accounted for.

Let  $p\in X$ be a singular point of type $A_i$  ($i\le 4$). We can find affine coordinates $z_1,z_2, z_3, z_4$ at $p$
such that $X$ admits an equation of the form $f(z)=z_2^2+z_3^3-z_4^2+g(z_1,z_2, z_3, z_4)$ with $g$ homogeneous of degree $3$. Then 
$f(z)+ a_0+a_1z_1+\cdots a_{i-1}z_1^{i-1}$ defines a deformation of $X$ as a  cubic threefold which 
exhibits a universal deformation of the germ  $(X,p)$. The basic theory of hypersurface singularities \cite{looij:icis} then tells us that  there is a submersive map germ $(S,o)\to \Sfrak_{i+1}\bs (H^i,0)$ such that $D(f)$ is the preimage of $\Sfrak_{i+1}\bs Z_i$ (recall that $H_i\subset \CC^{i+1}$ is the hyperplane defined by $z_0+\cdots +z_i=0$ on which $\Sfrak_{i+1}$ acts in the obvious manner). This implies that $D(f)$ is irreducible. Hence $\MM$ consists of  a single summand. The fact that the discriminant germs of 
$\Gg(\MM)\bs (V(\MM),0)$ and $\Sfrak_{i+1}\bs (H^i,0)$ are stably isomorphic implies that $\MM$ must be of type $\LL_i$.

We thus defined a map  from the singular set of $X$ to the set of summands of $\MM$. 
It is clear that the map is onto. It is  also is injective: if distinct singular points 
$p$ and $q$ of $X$ define the same summand, then  they define the same irreducible component of $D(f)$.
The generic point of this irreducible component parametrizes cubic threefolds with at least two singular points.
This contradicts the fact that the locus of cubic threefolds having at least two singularities is of codimension $\ge 2$.
\end{proof}

\begin{corollary}\label{cor:}
A universal deformation of a stable cubic threefold $X$ defines a versal deformation 
of the multigerm $(X, X_\sg)$ (by which we mean that the singularities can be deformed independently).\hfill $\square$
\end{corollary}

We use  our partial  classification  of geometric sublattices of small rank to deduce that strata of low codimension are irreducible.

\begin{corollary}\label{cor:irredinlowcodimension}
The loci $\Mcal_\ct(k_1A_1+k_2A_2 +k_3A_3)$ with $k_1+2k_3+3k_3\le 3$ are irreducible.
\end{corollary}
\begin{proof}
An irreducible component of $\Mcal_\ct(k_1A_1+k_2A_2 +k_3A_3)$ defines a
$\G_\ct$-orbit of sublattices  $\NN\subset \LL^\ct$ with $\NN\cong k_1\LL_1+k_2\LL_2 +k_3\LL_3$
such that no  $\LL_1$-summand of $\NN$ has the property that its sum with $\LL^0_3$ is contained
in an overlattice isomorphic with $\LL_4$ (otherwise this irreducible component would lie in the special 
divisor).   So by Proposition \ref{prop:lowdegreeclass}, the only way $\LL^0_3\operp \NN\subset \LL$ 
can fail to be primitively embedded is when $\NN\cong 3\LL_1$ and has as  its primitive hull a copy of $\LL_3$.
This proposition also tells us that all the primitively  embedded  $\LL^0_3\operp \NN\subset \LL$ make up a 
single $\U(\LL)$-orbit. This implies that the primitively embedded $\NN\subset \LL^\ct$  make up a single 
$\G^\ct$-orbit, even when $\NN\cong \LL_3$, for the standard embedding of $2\LL_3\subset \LL$ admits 
a unitary transformation which exchanges the factors. This implies the corollary, for the nonprimitively 
embedded copies of $3\LL_1$ in $\LL^\ct$ are already subsumed by the embeddings $\LL_3$, which make up an irreducible component and parametrize the cubic threefolds whose  singular locus is a singleton of type $A_3$.
\end{proof}

The connected components (strata) of $\Mcal_\ct(rA_3+sA_3)$  are in bijective correspondence the $\G^\ct$-orbits of the geometric sublattices of $\LL$ which contain $\LL^0_3$ as a direct summand. These are given by 
Corollary \ref{cor:coarse-ct}. We label them accordingly in Table \ref{table:latticeemb}. We shall need to identify these strata in geometric terms. This will be the subject of Section \ref{sect:coarsestrata}, but we make a beginning below.

\begin{small}
\begin{table}
\begin{centering}
\begin{tabular}{c|c|c|c|c|c}
$\Mcal_\ct(A_3)$ &  $\Mcal_\ct(2A_3)$ & $\Mcal_\ct(2A_3)'$ & $\Mcal_\ct(3A_3)$ & $\Mcal_\ct((2A_3)'+A_3)$ & $\Mcal_\ct(3A_3)'$\\
\hline
$\LL^0_3\operp\LL_3$ & $\LL^0_3\operp 2\LL_3$ & $(\LL^0_3\operp 2\LL_3)'$ & $\LL^0_3\operp(3\LL_3)'$ & 
$\LL_3\operp (\LL_3^0\operp 2\LL_3)' $& $(\LL^0_3\operp\LL_3)'$\\
\end{tabular}
\vskip5mm
\begin{tabular}{c|c|c|c}
$\Mcal_\ct(A_4)$ & $\Mcal_\ct(A_3+A_4)$ & $\Mcal_\ct((2A_3)'+A_4)$ & $\Mcal_\ct(2A_4)$\\
\hline
$\LL^0_3\operp \LL_4$ & $\LL^0_3\operp\LL_3\operp \LL_4$ & $(\LL^0_3\operp 2\LL_3)'\operp \LL_4$ & 
$\LL^0_3\operp 2\LL_4$\\
\end{tabular}
\end{centering}
\caption{The strata of  stable cubic threefolds  with only $A_3$ and $A_4$-singularities and the geometric sublattices representing them.}\label{table:latticeemb}
\end {table}
\end{small}
\medskip

\subsection{Geometric structure of some strata}\label{subsect:geomstructure}
We analyze  what happens when  the cubic threefold  acquires  double point singularities and how 
this leads to an interpretation  of discriminant strata as moduli spaces of curves.
The case when it has only one ordinary double point (an $A_1$-singularity) was already treated in 
some detail by Clemens-Griffiths \cite{cg}
and subsequently revisited by Collino-Murre \cite{cm},   Casalaina-Jensen-Laza \cite{cjl} and 
Casalaina-Grushevky-Hulek-Laza \cite{cghl}.
These authors study such degenerations through the Hodge structure of the cohomology of the 
cubic threefold (or equivalently, their 
intermediate jacobian). We need to do this through the Hodge structures of the 
associated cubic fourfolds with $\mu_3$-action.
\\

The mirror associated with a standard $3$-vector in $\LL^\ct$  is associated with an Eisenstein lattice of type 
$\LL_1\operp\LL_3\operp\LL_4\operp\HH$, 
a lattice that was also used to parametrize genus 4 curves. This is no coincidence,  because these curves show up geometrically. 

Suppose  that the stable cubic threefold $X\subset \PP V$ has at $p\in X$ singularity (so of type $A_i$ with $i\le 4$). 
The projectivized tangent cone $Q_p\subset  \PP(T_p\PP V)$ of $X$ at $p$ is then a quadric which is nonsingular
if $p$ is of type $A_1$, but is a quadric cone otherwise. 
The locus $C_p\subset Q_p$ which parametrizes the lines on $X$ through $p$ 
is the complete  intersection of $Q_p$ with a cubic surface. Since $X$ has multiplicity 2 at $p$, projection away from $p$ defines a birational map from 
$X$ to the projectivized tangent space $\PP T_pV\cong \PP^3$. This map becomes a morphism after 
blowing up $p$ in $\PP V$ and taking the strict transform $\widetilde X$ of $X$:
\[
X\xleftarrow{\pi} \widetilde X\xrightarrow{\pi'} \PP T_pV\cong \PP^3.
\]
Then $\pi^{-1}p=Q_p$. The exceptional set of $\pi'$ is the $\PP^1$-bundle $E_{C_p}$ over $C_p$ 
whose fibers are the lines through $p$ parametrized by $C_p$ and $\pi'|E_{C_p}$ is the retraction onto $C_p$.

To be concrete:  choose an affine coordinate system $(z_1, \dots, z_4)$ for $\PP V$ whose origin is $p$.
Then $X$ is given be an equation  of the form $f_2+f_3$ with $f_i$ homogeneous of degree $i$. The quadric   $Q_p$ is defined by $f_2$ and $C_p$ is the degree $6$ curve on $Q_p$ defined as the complete intersection defined by $(f_2, f_3)$. The birational  map above  is given by the projection 
$(z_1,z_2,z_3,z_4) \mapsto [z_1: z_2 : z_3 :z_4]$. This explains the first part of the following lemma.

\begin{proposition}[Casalaina-Laza, \cite{cl}, Prop.\ (3.1) and Cor.\ 3.2]\label{prop:curveCp}
Assume that  the stable cubic threefold  $X$ has  $p$ as a  singular point and let $A_i$ be its type.
Then $C_p$ is of arithmetic genus 4 and is embedded in  the projective $3$-space $\PP T_pV$ by its dualizing sheaf. 

For $i=1$,  $Q_p$ is smooth. For $i\ge 2$,  $Q_p$ is  the cone over a smooth  conic and then the ruling of $Q_p$ defines a hyperelliptic involution of $C_p$.  

For $i=2$,  $C_p$ does not pass through the vertex 
of that cone, but for  $i=3, 4$ it does:  $C_p$ exhibits there  a node  for $i=3$  and an ordinary  cusp for $i=4$.
\end{proposition}

We expand a little on  the cases $i\ge 2$,  as this will lead us to an inverse construction, which produces  all stable cubic threefolds with a distinguished singular point. Suppose  $Q\subset \PP^3$ is a cone over a smooth conic. Then the blowup of its vertex  $\tilde Q\to Q$ is a  Hirzebruch surface whose  exceptional curve $E$ has self-intersection 
number $-2$  and which comes  canonical retraction of $\tilde Q$ onto $E$. 
The class $e$ of $E$  and the class $f$ of a fiber of $\tilde Q\to E$ make up a basis of $\pic(\tilde Q)=H^2(\tilde Q)$. 
If  $C\subset Q$ is cut out by a cubic surface then it has class $3e+6f$.  If it passes through the vertex and  exhibits there an ordinary double point resp.\ a cusp,  then its strict transform $\tilde C\subset \tilde Q$ has class  $2e+6f$. Since $(2e+6f)\cdot f=2$, the projection 
$\tilde C\to E$ is of  degree 2 and  this endows  $\tilde C$ with a  hyperelliptic involution $\iota$. 
The  curve $\tilde C$  is of  genus $3$, smooth near  $E$, and cuts out on  $E$ in a divisor $\Pcal$  of degree $2$ which is reduced 
when $i=3$ and nonreduced when $i=4$.

\begin{proposition}\label{prop:converseconstruct}% \EL{Gebruiken we dit alleen in \ref{cor:Ai_moduli}?}
Let $W$ be a complex vector space of dimension $4$ and let $C\subset \PP(W)$ be a reduced complete intersection of bidegree $(2,3)$, where we assume that the (unique) quadric $Q$ which contains $C$ has at most one singular point (then called  the vertex). We also assume that $C$  has at worst $A_4$-singularities in the smooth part of $Q$ and  has at the vertex of $Q$ (if it exists) at worst an $A_2$-singularity.

Then $C$ determines (up to projective equivalence) a stable cubic threefold $X_C$ in $\PP V$  with  a distinguished singular point  $p$ and for which the variety of lines on $X$ through  $p$ is identified with $C$: it is obtained by blowing up $C$ in its ambient $\PP(W)$ followed by contraction of  the strict transform of $Q$.  The singularities of  $X_C$ relate to those of $C$ as in Proposition \ref{prop:curveCp}.
\end{proposition}

\begin{remark}\label{rem:canemb}
Since  such a $C$ is embedded in $\PP(W)$ by its dualizing sheaf,  this construction only depends on $C$ as an abstract curve.
\end{remark}

\begin{proof}
We choose an isomorphism of vector spaces $\CC\oplus W\cong V$ and denote by $p\in \PP V$ the point defined by the first summand of $\CC\oplus W$. 
Let $f_2, f_3\in \CC[W]$ be homogeneous of  degree $2$ resp.\ $3$ and define 
$C$  (so that $f_2$ defines $Q$) and consider the map $ w\in W\to (f_2(w)w, -f_3(w)\in W\oplus \CC=V$. It defines a rational section
\[
\varphi: \PP W\dashrightarrow \PP V 
\]
of the projection $pr_p$. It is clear that $\varphi$ blows up at $C$ and maps $Q$ to $p$.
We denote the coordinate of  $\CC$-summand by $w_0$. So the  image  of $\varphi$ is given by   $w_0=-f_3(w)/f_2(w)$ and hence its closure satisfies  the equation 
$f_3(w)+w_0 f_2(w)$. This is an irreducible polynomial and defines therefore an irreducible threefold.
Another choice $f'_3$ of $f_3$  has the form $f_3+\ell f_2$, where $\ell$ is a linear form on $W$. This changes
the equation by replacing $z_0$ by $z_0+\ell(z)$ and so gives the same cubic threefold up to a projective transformation.
We leave it to the reader to check that this defines the inverse of the construction of Proposition 
\ref{prop:curveCp}.
\end{proof}

We can now make Proposition  \ref{prop:curveCp} a bit more precise.

\begin{corollary}\label{cor:Ai_moduli}
We have an identification between the following   loci singular  cubic threefolds and moduli spaces of curves:
\begin{description}
\item [$\Mcal_\ct(A_1)$] $\Mcal_4\ssm (\Hyp_4\cup\Mcal_4^\Theta)$, the moduli space of smooth nonhyperelliptic  genus 4 curves on a smooth quadric,
\item [$\Mcal_\ct(A_2)$]  $\Mcal^\Theta_4\ssm \Hyp_4$, the moduli space of smooth nonhyperelliptic genus 4 curves on a quadric cone,
\item [$\Mcal_\ct(A_3)$]  $\Hyp'_{3,[2]}$, the moduli space of smooth hyperelliptic curves of genus $3$ endowed with a 2-element subset $\Pcal$ which is not invariant under the hyperelliptic involution, 
\item [$\Mcal_\ct(A_4)$] $\Hyp'_{3,1}$, the moduli space of smooth hyperelliptic curves of genus $3$ endowed with a non-Weierstra\ss\ point. 
\end{description}
In particular these loci are irreducible.
\end{corollary}
\begin{proof}
Here in the last two cases, the less evident direction goes as follows. Given a pair $(C,\Pcal)$ consisting of a smooth hyperelliptic curve $C$ of genus $3$ endowed with a degree $2$-divisor $\Pcal\ge 0$ not invariant under the hyperelliptic involution, then the  complete linear system $|\omega_C(\Pcal)|$ of  degree 6 defines a morphism of $C$ to a $\PP^3$ whose image lies on a cone  over a smooth conic
and realizes the case $i=3$ and $i=4$ above.
In particular the  corresponding  loci  in the ball quotient are  connected and are open subsets of ball quotients 
(associated with the lattices of type $\LL_1\operp\LL_{4-i}\operp\LL_4\operp\HH$). 
\end{proof}

\begin{remark}\label{rem:chordallimit}
Smooth hyperelliptic curves of genus $3$ endowed with a hyperelliptic orbit (which may be reduced to a Weierstra\ss\ point)
do not occur in this manner, but arise on the chordal locus, which we will discuss below.
\end{remark}

\begin{remark}\label{rem:2a1}
This argument also shows that $\Mcal_\ct(2A_1)$ is dominated  by the moduli space of  nonhyperelliptic  genus 4 curves on a smooth quadric with a node and is hence irreducible.
\end{remark}

\begin{remark}
Allcock \cite{allcock:3GIT} and Yokoyama \cite{yokoyama:3GIT} proved that the GIT compactification
$\Mcal_\ct^\ss\supset \Mcal_\ct^\st$ adds to $\Mcal_\ct^\st$ two singletons and an affine  rational curve which has one of these singletons as its unique boundary point.
As we shall recall below, the literature  provides us with a  more precise version of  
Theorem \ref{thm:cubicperiods} that also says something about this GIT boundary.
One of the singletons of that boundary is  represented by the \emph{chordal cubic} that we will 
discuss in the next subsection. 
The other singleton and the curve  in the GIT boundary $\Mcal_\ct^\ss\ssm\Mcal_\ct^\st$ 
are both mapped to cusps of the ball quotient. Indeed,  $\G^\ct$ has two orbits in the set of primitive isotropic sublattices  
$\II\subset \LL^\ct$: 

\begin{enumerate}
\item[($\infty_\spl$)]  The curve in the GIT boundary  parametrizes threefolds with two singularities of type $A_5$ 
with one of the cubic threefolds parametrized by it having an additional $A_1$-singularity.
The associated cubic fourfold has then two simple elliptic singularities of type $\tilde E_8$ (and possibly  an additional $A_2$-singularity). This means that   
$\II^\perp/\II$  contains a sublattice of type $2\LL_4\operp \LL_1$. Since $2\LL_4\operp \LL_1$ has  no overlattice of the same rank, we must have equality: 
$\II^\perp/\II\cong \LL_1\operp2\LL_4$. Note that such an $\II$ is realized by the span of an  isotropic basis element of the $\HH$-summand of $\LL_1\operp 2\LL_4\operp\HH$ and that $\II^\perp/\II$ has a special $3$-vector.
\item[($\infty_\std$)]  One of the  singletons in the GIT boundary defines a cubic threefold with three $D_4$-singularities. The associated cubic fourfold has then three simple elliptic singularities of type $\tilde E_6$.  This means that $\II^\perp/\II$ contains a sublattice  of type $3\LL_3$. This too, has no overlattices of the same rank and hence  is all of  $\II^\perp/\II$.  This lattice  contains no special $3$-vectors.
\end{enumerate}
 
The curve $\infty_\spl$ and the singleton $\infty_\std$  in the GIT boundary $\Mcal_\ct^\ss\ssm\Mcal_\ct^\st$ are both mapped to cusps. These cusps are different, even in the Allcock ball quotient and are  (via Property (ii) of Subsection \ref{subsect:ac}) characterized by the property that they are represented by primitive isotropic rank one lattices $\II$ in $\LL$ for which $\II^\perp/\II$ is isomorphic to $3\LL_4$ resp.\  its $3$-vectors span  a sublattice isomorphic with $4\LL_3$.
\end{remark}

\subsection{The chordal divisor}\label{subsect:chordal}
A \emph{chordal cubic} $Y\subset \PP V$ is by definition the union of the secants of a normal rational quartic curve in $\PP V$. Such a cubic  is obtained (and best described) by fixing a 2-dimensional complex vector space $W$  and a linear isomorphism $\varphi: \sym^4W\cong V$.
Such a $\varphi$ identifies the linear system of positive divisors degree $4$ on the projective line $\PP W$ with $\PP V$. The quartic  curve  is given by the positive divisors that are 4 times a point (so the image of the Veronese map $\PP W\to \PP V$). 
One may check that  the associated chordal cubic $Y$ is the closure of the image of the locus of degree 4 divisors in harmonic position (i.e., with cross ratio $\pm1$). 
The isomorphism $\varphi$ induces a homomorphism $\SL(W)\to \SL (V)$  whose image is precisely the 
$\SL (V)$-stabilizer of $X$ (its kernel is the center $\{\pm 1_W\}$). Let $F_0\in \Kcal$ be a defining equation for $X$.
It is shown in \cite{allcock:3GIT} and \cite{yokoyama:3GIT} that $F_0$ is $\SL(V)$-semistable and that the $\SL(V)$-orbit  of  $[F_0]\in \PP \Kcal$ is closed in the semistable locus $\PP \Kcal^\ss$.
The normal bundle of this orbit comes with an $\SL(V)$-action. This restricts to an $\SL(W)$-action on the normal space of this orbit at $[F_0]$. 

\begin{lemma}\label{lemma:transversal}
The normal space to the $\SL(V)$-orbit  of  $[F_0]\in \PP \Kcal$ is as an $\SL(W)$-representation canonically isomorphic to $\CC[W]_{12}$ (and hence irreducible).
In fact, an $\SL(W)$-invariant transversal slice to the $\SL(V)$-orbit of $[F_0]$ is defined by 
$G\in \CC[W]_{12}\mapsto [F_0+\tilde\varphi(G)]\in \PP \Kcal$, where 
\[
\begin{CD}
\tilde \varphi: \CC[W]_{12} @>>> \CC[\CC[W]_4]_3 @<{\CC[\varphi]_3}<{\cong}< \Kcal
\end{CD}
\]
 is the evident $\SL(W)$-equivariant map. 
 \end{lemma}
\begin{proof}
This is straightforward: the codimension of the $\SL(V).[F_0]\in \PP \Kcal$ is readily computed to be $13$. 
Note that $\Kcal$ is as a representation of $\SL (W)$ isomorphic to the dual of $\sym^3(\sym^4 W)$.
The choice of a Cartan subgroup of $\SL(W)$ amounts to a decomposition of $W$ into $1$-dimensional subspaces, $W=W_1\oplus W_{-1}$,  and then it is clear that the highest weight space of  
$\Kcal$
is $\CC[W_{-1}]_{12}$. This proves that $\Kcal$ contains a unique copy $\CC[W]_{12}$.
It remains to observe  that the tangent space of $\SL(V).[F_0]$ at $[F_0]$ is an $\SL(W)$-module of lower dimension, so that $\tilde\varphi$  must be transversal to $\SL(V).[F_0]$ at $[F_0]$.
\end{proof}

The   blowup of  $\SL(V).[F_0]$ in $\PP \Kcal^\ss$, 
\[
\widetilde\PP \Kcal^\ss\to \PP \Kcal^\ss,
\]
has as exceptional divisor the projectivized normal bundle of the $\SL(V)$-orbit  of  $[F_0]$ in $\PP \Kcal^\ss$.  Lemma \ref{lemma:transversal} shows that the fiber over $[F_0]$ is naturally identified with $\PP (\CC[W]_{12})\cong \PP (\sym^{12}W)$. Concretely, the cubic threefold  defined by   $F_0+\tilde\varphi(G)$ has under the map
$\PP\varphi: \PP W\to \PP V$ as preimage the degree 12 divisor $D$ on $\PP W$ defined by $G$. Note that once we specify a divisor $D$ on $\PP W$ whose support has size $\ge 3$, then the  $\SL(V)$-stabilizer of the pair $(Y, D)$ becomes finite.

This blowup  is algebraically obtained as the proj of the (homogeneous) ideal in $I\subset \CC[\Kcal]$ of polynomials that vanish on  $\SL(V).[F_0]$. 
This ideal also defines a  linearization of the $\SL(V)$-action on $ \widetilde\PP \Kcal^\ss$. 
The associated stability property is the same as the old one away from $\SL(V).[F_0]$ and is over $\SL(V).[F_0]$ 
the same as  Hilbert-Mumford stability: $(Y,D)$ is stable if and only if  the multiplicities of $D$ are $\le 5$. 
In other words, the space  of stable  $\SL(V)$-orbits in $\widetilde\PP \Kcal^\ss$, which we shall denote $\widetilde\Mcal_\ct^\st$, 
contains  the space of Hilbert-Mumford stable $\SL(W)$-orbits  of degree 12 divisors   on  $\PP W$, denoted in Example \ref{example:DM} 
by $\Hcal^{0, \st}_{12}$, 
as a divisor whose complement in  $\widetilde\Mcal_\ct^\st$ can be identified with the space of GIT stable threefolds 
$\Mcal_\ct^\st$: $\Mcal_\ct^\st\cong \widetilde\Mcal_\ct^\st\ssm \Hcal^{0, \st}_{12}$.
In the references (\cite{act} Theorem 1.1 and \cite{ls} Theorem 3.1), 
Theorem \ref{thm:cubicperiods} includes the assertion  that the period map extends to an isomorphism of $\widetilde\Mcal_\ct^\st$  onto
$\Xcal_\ct$ that takes $\Hcal^{0, \st}_{12}$  onto the special mirror  divisor $D^{\ct, \spl}$. As we mentioned in Subsection \ref{subsect:ballmodular}, this is the Deligne-Mostow ball quotient that comes from associating to an effective  degree 12 divisor on $\PP^1$ 
the connected $\mu_6$-cover of $\PP^1$ having that divisor as discriminant and looking as the part of its cohomology that transforms according to the tautological character.  We prefer however to associate with that divisor the intermediate $\mu_2$-cover  so that 
$\Hcal^{0, \st}_{12}$ has the interpretation as the moduli  space of hyperelliptic curves of genus $5$ that are stable in this Hilbert-Mumford sense; we therefore
shall write $\Hyp_5^\st$ for this divisor, but will follow the custom to refer to it as the \emph{chordal divisor}. This
has the effect that we regard this divisor as a locus  with global inertia group $\mu_2$ (the hyperelliptic involution becomes a central element of its orbifold fundamental group).  So
\[
\Mcal_\ct^\st\cong \widetilde\Mcal_\ct^\st\ssm \Hyp_5^\st.
\]

\begin{remark}\label{rem:collino}
Degenerations of smooth cubic threefolds into the chordal cubic were already  considered by Collino \cite{collino}. He proved that the associated Hodge structure 
on the third integral cohomology does not degenerate (after a possible base change), but becomes that of a hyperelliptic curve of genus $5$ tensored with $\ZZ (-1)$.
This justifies identifying  the chordal divisor  with $\Hyp_5^\st$ (which is also the approach taken by  Casalaina-Jensen-Laza in \cite{cjl}).
\end{remark}

\begin{remark}\label{rem:modulariso}
The moduli space of cubic $3$-folds with a $A_4$-singularity is obtained as an open subset of the ball quotient  associated with the geometric lattice $\LL^0_3\operp\LL^4$.  Something similar is true for the part of the chordal divisor defining hyperelliptic genus  5 curves with an $A_3$-singularity:  the associated  geometric  lattice is then $\LL^0_4\operp\LL^3$. Since these lattices are $\G$-equivalent, we expect to find an isomorphism in terms of their modular interpretations. 

This is indeed the case: we have seen that in the first case we get the moduli space of hyperelliptic curves  $C$ of genus $3$ endowed with a non-Weierstra\ss\ point $z\in C$. 
On the other hand, if we normalize  a hyperelliptic genus 5 curve  with an $A_3$-singularity, then we get a hyperelliptic genus 3
curve and no information is lost, if we remember the preimage of the singular point. This preimage is a regular orbit  for the hyperelliptic involution. The isomorphism then simply assigns to $(C;z)$ the pair $(C; \{z, \iota(z)\})$. So both parametrize affine hyperelliptic curves 
given by an equation $y^2=x^8+c_2x^6+c_3x^5+\cdots +c_8$. Note that the orbit 
of $(c_2, \dots, c_8)\in \CC^8\ssm \{0\}$ under the $\CC^\times$-action $t.(c_2, \dots, c_8)= (t^2c_2, \dots, t^8c_8)$ is a complete invariant.
(The associated ball quotient is of course one of Deligne-Mostow type; its fractional exponents are $\tfrac{1}{6}$ with multiplicity $8$ and $ \tfrac{2}{3}$ with multiplicity $1$.)

This isomorphism extends to the other coarse strata in their closure: if a $\MM$ and $\MM'$ are geometric sublattices of $\LL$ which  
contain $\LL^0_3$ resp.\   $\LL^0_4$ as a direct summand and are $\G$-equivalent, then and isomorphism between the associated ball quotients is obtained directly in terms of  their modular interpretations and is obtained  as a specialization of the isomorphism described above.
\end{remark}

\begin{remark}\label{rem:gradedalgebra}
As we noticed in the discussion leading up to Theorem \ref{thm:cubicperiods}, we have an $\SL(V)$-equivariant identification of the  line  in $\sym^2(\sym^3 V)$ 
spanned  by $F^{-2}$ with the positive line in $V^\K$ represented by $\Hl^{3,1}_\circ(Y)$.  
Any $\G^\ct$-automorphic form $\Phi$ on $\BB^\ct$ of weight $w\ge 0$ defines via a $\LL^\ct$-marking 
of $Y$,  an element of 
$\Hl^{3,1}_\circ(Y)^{\otimes w}$ and hence, via the above isomorphism,  
an element in the line spanned by $F^{-2w}$, or equivalently, a linear form on the line spanned by $F^{2w}$.  
If we let $F$ vary, we then find that $\Phi$ produces an $\SL(V)$-invariant polynomial $e(\Phi)$ of degree 
$2w$ on  $\Kcal$. 
We thus obtain a $\CC$-algebra homomorphism of graded $\CC$-algebras:
\[
\textstyle e: \sum_{w\ge 0} \Hl^0(\BB^\ct, \Ocal_{\BB^\ct}(-w))^{\G^\ct}t^{2w}\to  \sum_{w\ge 0} I_{w}^{\SL(V)}t^w,
\]
where we note that the right hand side is zero in odd degrees  (because $-1_V\in\SL(V)$ acts on $\Kcal$ as $-1_\Kcal$) and the left hand side is zero in degrees not divisible by $3$ (because $\mu_3\subset \G^\ct$ acts faithfully on $\Ocal_{\BB^\ct}(-1)$  as fiberwise scalar multiplication). The map $e$ defines an integral extension  and  is of degree one if we projectivize (it is probably an isomorphism). We thus  obtain the isomorphism  $X_\ct\cong \widetilde\Mcal_\ct^\st$ by passing to the associated proj.
\end{remark}

\section{A universal reflection cover for stable cubic threefold moduli}\label{sect:level2} 

We show that  a  universal reflection cover  of the moduli space stable cubic threefolds is obtained by imposing a level 2 structure.

\subsection{A universal reflection cover for stable cubic threefold moduli} 
For a smooth cubic threefold $X$, there is no obvious relation between $\Hl_ 3(X)$ and $\Hl_ 4(Y)$. There is one at the prime $3$, but that is not relevant here, as we will be more concerned with the prime 2. So the discussion of $\Hl_ 3(X)$ necessitates a separate story.

The group $\Hl_3(X)$ is free abelian of rank $10$ and the intersection pairing makes this a unimodular symplectic lattice. 
We denote by $\Delta_X\subset \Hl_3(X)$ the set  of  vanishing cycles (the cycles represented by an embedded oriented $3$-sphere in $X$  that `vanish' when we let  $X$ degenerate into a cubic threefold with an $A_1$-singularity in the sense of Subsection \ref{subsect:discrstrat}). We recall from Subsection \ref{subsect:hypersurfacesing} that this set  has also the following properties 
\begin{enumerate}
\item [(i)] $\Delta_X$ generates $\Hl_3(X)$ over $\ZZ$,
\item [(ii)] the symplectic transvections $\{T_\delta\}_{\delta\in \Delta_X}$ generate in $\Sp(\Hl_3(X))$  the monodromy group (which  we shall denote by  $\G_{\Delta_X}$),
\item [(iii)] $\Delta_X$ is a $\G_{\Delta_X}$-orbit.
\end{enumerate}

Not every primitive vector in  $\Hl_3(X)$ is a vanishing cycle and so $\G_{\Delta_X}$ cannot  be the full symplectic group of 
$\Hl_3(X) $. Indeed, $\G_{\Delta_X}$ also fixes a quadratic form $q_X: \Hl_3(X; \FF_2)\to \FF_2$ 
which is characterized by the property that it takes on the image of $\Delta_X$ the value $1$ 
and that its associated alternating bilinear form, 
\[
(a,b)\in \Hl_3(X; \FF_2)\times \Hl_3(X; \FF_2)\mapsto q_X(a+b)-q_X(a)-q_X(b)\in \FF_2,
\] 
is the mod 2 reduction of the intersection pairing. 
These two  properties  formally imply that $q_X$ is invariant under the mod 2 reduction 
of the monodromy group. 
Since the intersection pairing is unimodular, this  quadratic form is nonsingular.  According to Randal-Williams \cite{rw}, Remark 3.8, its Arf invariant is $1$. This property characterizes the isomorphism type of the quadratic form $q_X$.  As mentioned in Subsection \ref{subsect:hypersurfacesing}, the monodromy group is 
not subject to other constraints, in other words,  it is the preimage of  the orthogonal group 
$\Orth(q_{X})\subset \Sp(\Hl_3(X; \FF_2))$ under the reduction map $\Sp(\Hl_3(X))\to \Sp(\Hl_3(X; \FF_2))$. It follows 
that the image of $\Delta_X$ is  all of $q_X^{-1}(1)$.

Recall that  $(\Tcal, q)$ is a nonsingular  quadratic space over $\FF_2$  of dimension 10 and Arf invariant $1$ and comes with  an embedding of the permutation group $\Sfrak_E$  in $\Orth(\Tcal)$. 

\begin{definition}\label{def:Tmarking}
A \emph{$\Tcal$-marking} of a smooth cubic threefold $X$ 
is the choice of an isomorphism $(\Hl_3(X, \FF_2), q_X)\cong (\Tcal, q )$ of quadratic 
$\FF_2$-vector spaces. An \emph{isomorphism} between two such is an isomorphism 
between the threefolds in question that takes the marking of one to  the marking of the other.
\end{definition} 

It is clear from the preceding that every smooth cubic threefold $X$ admits a
$\Tcal$-marking and that the group  $\Orth(\Tcal)$ permutes its  $\Tcal$-markings simply transitively. 
The smooth $\Tcal$-marked cubic threefolds are parametrized by a Deligne-Mumford stack that we denote $\underline\Mcal_\ct[\Tcal]$.
The underlying orbifold $\Mcal_\ct[\Tcal]$ is connected, for  the monodromy representation on $\Hl_3(X, \FF_2)$ is the full orthogonal group of $q_X$.

Let $X$ be a stable cubic threefold whose singular set is of type $k_\pt A_\pt$. So for 
$i=1,2,3,4$ and $j=1, \dots  k_i$ we have a singular point $x_{i,j}\in X$ of type $A_i$ and these are pairwise distinct. 
Let $\g$ be an arc  defining a degeneration of  $X_o$  into $X$
(see Subsection \ref{subsect:discrstrat}). We observed that  
this determines on   $X_o$ a vanishing manifold $F_\g$ (each connected component is a Milnor fiber 
of a $3$-dimensional $A_i$-singularity). The reduced homology of such a connected component is 
 concentrated in degree $3$ and is spanned by a set of vanishing cycles: this yields a subset
of $\Delta_{X_o}(X, \g)\subset \Delta_{X_o}$. It dependence on $\g$ is clearly through the homotopy class of $\g$
(in an obvious sense). The homotopy class of any other $\g'$ is obtained by precomposing $\g$ with a loop based at $[X_o]$ which parametrizes smooth cubic threefolds. Such a loop induces a symplectic transformation of $\Hl_3(X_o)$ which is orthogonal when reduced mod $2$ (and any such transformation is thus obtained). In particular,  the orbit of  
$\Delta_{X_o}(X, \g)$ under this monodromy group  is well-defined.

\begin{lemma}\label{lemma:stratalattices1}
The map $\Hl_3(X_o)\to \Hl_3(X)$ defined by $\g$  is onto  and its kernel is the vanishing homology defined by $\g$.  This identifies 
$\Hl_3(X)$ with the  group of covariants of $\Hl_3(X_o)$ with respect to the action of the (local monodromy)
group $\G_{\Delta_{X_o}(X, \g)}$ generated by the  $T_\delta$ with $\delta\in \Delta_{X_o}(X, \g)$.

Dually, $\g$ defines a primitive  embedding  $\Hl^3(X)\hookrightarrow \Hl^3(X_o)$ whose cokernel embeds in 
$H^3(F)$ and whose image is the group of $\G_{\Delta_{X_o}(X)}$-invariants.

The same results hold for homology with $\FF_2$-coefficients.
\end{lemma}
\begin{proof}
 Consider the following piece of the exact homology sequence of the pair $(X_o,F)$:
 \[
\cdots\to \Hl_3(F)\to  \Hl_3(X_o)\to  \Hl_3(X_o,F)\to  \Hl_2(F)\to\cdots
\] 
The first assertion then follows  from the vanishing of  $\Hl_2(F)$ and the identification $\Hl_3(X_o,F)\cong \Hl_3(X, X_\sing)\cong \Hl_3(X)$. The chacterization of $\Hl_3(X)$ as a group of covariants then simply follows from
the fact that for every $\delta\in \Delta_{X_o}(X)$, the symplectic transvection $T_\delta$ is nontrivial (there always exists a vanishing cycle $\delta'$ with $\delta\cdot\delta'=1$).

The proof of the cohomological version is similar and the last assertion follows from this.
\end{proof}

Let us express  the preceding  in terms  of markings. The normalization of $ \underline{\Mcal}_\ct^\st$ in the cover 
$\underline{\Mcal}_\ct[\Tcal]\to \underline{\Mcal}_\ct$  yields a ramified 
$\Orth(\Tcal)$-cover
\[
\underline{\Mcal}_\ct^\st[\Tcal]\to \underline{\Mcal}_\ct^\st.
\]
As agreed earlier, we denote the  underlying orbifolds simply by omitting the underlining, so that we have an orbifold  morphism  $\Mcal_\ct^\st[\Tcal]\to \Mcal^\st_\ct$.  The preceding shows that given 
a  connected component  $S$ of  $\Mcal_\ct(k_\pt A_\pt)$, then each connected component $\hat S$ of its preimage in  
$\Mcal_\ct^\st[\Tcal]\to \Mcal_\ct^\st$ determines a subset of $\Delta(\FF_2)$ and that the subsets thus obtained are transitively permuted by $\Orth(\Tcal)$. We thus find:
 
 \begin{corollary}\label{cor:stratalattices2}
Let $\hat S$ be a  connected component  of $\Mcal_\ct(k_\pt A_\pt)[\Tcal]$. Then $\hat S$ parametrizes 
stable cubic threefolds $X$ with a singularity configuration of type $k_\pt A_\pt$ endowed with a \emph{residual $\Tcal$-marking} in the following sense:  with $\hat S$ is associated a subspace  
$\Tcal(\hat S)\subset \Tcal$  spanned by a root system of type $k_\pt A_\pt$ such that  $\Hl_3(X; \FF_2)\cong \Hl^3(X_\reg; \FF_2)$  is identified with $\Tcal/\Tcal(\hat S)\cong \Tcal(\hat S)^\perp$. 
The map $\hat S\mapsto \Tcal(\hat S)$  reverses the incidence  relation: if $\hat S$ lies in the closure of $\hat S'$, then 
$\Tcal(\hat S)\supset \Tcal(\hat S')$.

When $k_4=0$ and $k_1+2k_2+3k_3\le 3$, then this defines a bijection  between  the set of connected components of $\Mcal_\ct(k_\pt A_\pt)[\Tcal]$ and the subspaces of $\Tcal$ of  type $k_\pt A_\pt$. 
\end{corollary}
\begin{proof}
The  first paragraph merely records our earlier observations. The las paragraph is a consequence of the fact that in this range  the loci $\Mcal_\ct(k_\pt A_\pt)$ are irreducible by Corollary \ref{cor:irredinlowcodimension}.
\end{proof}

Denote by $\Mcal_\ct^\st(k_\pt A_\pt)$ the closure of $\Mcal_\ct(k_\pt A_\pt)$ in $\Mcal_\ct^\st$.
So a connected component  of  $\Mcal_\ct(k_\pt A_\pt)$ determines an irreducible component of $\Mcal_\ct^\st(k_\pt A_\pt)$ and vice versa. We can therefore also express the preceding results in terms of such irreducible components. For example, 
by Corollary \ref{cor:Ai_moduli} the locus $\Mcal^\st_{\ct}(A_i)$ is for $i=1,2,3,4$) 
% and $\Mcal^\st_{\ct}(2A_1)$ 
irreducible and so  the  irreducible components of their preimages in
$\Mcal^\st_{\ct}[\Tcal]$ are in  bijective correspondence with subspaces of $\Tcal$ of type $A_i$ in the sense of Definition \ref{def:geomsubgroup1}.

\subsection{Proof of Theorem \ref{thm:ctcover}}
We prove this theorem, which is  about a 10-dimensional ball quotient, via the  modular interpretation of that  ball quotient
(we do not know how to obtain such a result without this). It then takes the following form.

\begin{theorem}\label{thm:ct1}
The $\Orth(\Tcal)$-covering $\Mcal^\st_{\ct}[\Tcal]\to \Mcal^\st_{\ct}$ is a  universal reflection covering relative to  $\Mcal^\st_{\ct}(A_1)$. 
\end{theorem}

A key step of the proof  is:

\begin{lemma}\label{lemma:sconnectivity}
The preimage $\Mcal^\st_{\ct}(A_1)[\Tcal]$ of $\Mcal^\st_{\ct}(A_1)$ in $\Mcal^\st_{\ct}[\Tcal]$ is simply connected.
\end{lemma}
\begin{proof}
We regard the collection of  irreducible components of $\Mcal^\st_{\ct}(A_1)[\Tcal]$  as a closed covering of that space.  
Such an  irreducible component defines a vanishing cycle  $\delta\in\Tcal$. These vanishing cycles make up a $\Orth(\Tcal)$-orbit (this is the set $q^{-1}(1)$) and so we index these irreducible components accordingly: 
\[
\Mcal^\st_{\ct}(A_1)[\Tcal]=\cup_{\delta\in q^{-1}(1)}\Mcal^\st_{\ct}[\delta].
\]
The nerve of this covering is a simplicial complex $\Sigma(\Tcal)$ with vertex set $q^{-1}(1)$. We want to apply 
the generalized van Kampen theorem \ref{lemma:generalizedvK} to this covering.  
Each $\Mcal^\st_{\ct}[\delta]$ is isomorphic with $\widehat\Hyp_4$ and hence simply connected by Proposition \ref{prop:hypkernel}. If  $\delta, \delta'$ are distinct in $q^{-1}(1)$, then by  Corollary  \ref{cor:stratalattices2}, 
$\Mcal^\st_{\ct}[\delta]\cap \Mcal^\st_{\ct}[\delta']$ is either an irreducible component of 
$\Mcal^\st_{\ct}(2A_1)[\Tcal]$ (when $\delta\cdot \delta'=0$) or of $ \Mcal^\st_{\ct}(A_2)[\Tcal]$
(when $\delta\cdot \delta'=1$). 

If we have three pairwise distinct $\delta, \delta', \delta''$ elements of  $q^{-1}(1)$, then 
$\Mcal^\st_{\ct}[\delta]\cap \Mcal^\st_{\ct}[\delta']\cap \Mcal^\st_{\ct}[\delta'']$ is by Corollary  \ref{cor:stratalattices2} it   given by an embedding of $3A_1$, $A_2+A_1$ or $A_3$ in $\Tcal$ and hence even connected.
This shows among other things  that the nerve of this covering contains the $2$-skeleton of the  full simplex spanned by $q^{-1}(1)$ and hence is simply connected.
It then follows from Lemma \ref{lemma:generalizedvK} that $\Mcal^\st_{\ct}(A_1)[\Tcal]$ is simply connected.
 \end{proof}

The proof of Theorem \ref{thm:ct1} will use a Lefschetz hyperplane theorem for quasiprojective varieties. It says that 
if $X$ is an irreducible  complex quasiprojective variety and $Y\subset X$ is a hyperplane section such that $X\ssm Y$ is smooth, then 
up to homotopy, $X$ is obtained from $Y$ by attaching  cells of dimension $\ge \dim_\CC X$.  So when $\dim_\CC X\ge 3$, the variety 
$Y$ is connected and $Y\subset X$ induces an isomorphism on fundamental groups   (see for instance Hamm-L\^e \cite{hle}, Thm.\ 1.1.1). 

\begin{proof}[Proof of Theorem \ref{thm:ct1}]
The discriminant  hypersurface $\Mcal^\st_{\ct}(A_1)$ of  $\Mcal^\st_{\ct}$ supports an ample divisor.
This implies that $\Mcal^\st_{\ct}(A_1)[\Tcal]$ supports an ample divisor in $\Mcal^\st_{\ct}$. The complement
of $\Mcal^\st_{\ct}(A_1)[\Tcal]$ in $\Mcal^\st_{\ct}[\Tcal]$ is $\Mcal_{\ct}[\Tcal]$ and this is a smooth orbifold.

We want to apply the above   Lefschetz hyperplane theorem to the case when $X=\Mcal^\st_{\ct}[\Tcal]$ and 
$Y=\Mcal^\st_{\ct}(A_1)[\Tcal]$. We know  that $\Mcal^\st_{\ct}(A_1)$ represents an ample divisor, and hence so 
does $\Mcal^\st_{\ct}(A_1)[\Tcal]$, 
but the cited result does not directly apply here, as the complement $\Mcal_{\ct}[\Tcal]$ might be  an
 orbifold and so need not be smooth. We can however pass to a smooth setting by considering the moduli space of cubic  threefolds with a 
 principal level $6$ structure that lifts the $\Tcal$-marking (this amounts to imposing an additional level $3$ structure). 
  The  orbifold covering $\Mcal_{\ct}[6]\to \Mcal_{\ct}$ factors through $\Mcal_{\ct}[\Tcal]$ and it is well-known that 
  this cover  has a smooth total space (see for example Serre \cite{serre}). We take its normalization over 
  $\Mcal^\st_{\ct}$:  $\Mcal_{\ct}^\st[6]\to \Mcal_{\ct}^\st$.
It then follows that the inclusion $\Mcal_{\ct}^\st(A_1)[6]\subset 
\Mcal_{\ct}^\st[6]$ induces an isomorphism on fundamental groups. This implies the same for  
$\Mcal_{\ct}^\st(A_1)[\Tcal]\subset  \Mcal_{\ct}^\st[\Tcal]$.
Hence by  Lemma \ref{lemma:sconnectivity},  $\Mcal^\st_{\ct}[\Tcal]$ is simply connected. 

It remains to see that this is a reflection covering relative to  $\Mcal^\st_{\ct}(A_1)$.
It is clear that the ramification of $\Mcal^\st_{\ct}[\Tcal]\to \Mcal^\st_{\ct}$ is over 
$\Mcal^\st_{\ct}(A_1)$. The monodromy of this ramification is of order two (for a symplectic transvection in 
$\Tcal$ has that property)
and so we have a reflection covering relative to  $\Mcal^\st_{\ct}(A_1)$ indeed.
\end{proof}

\begin{corollary}\label{cor:ct2}
The ball quotient $(\Xcal_\ct, D_\ct)$  has a universal reflection cover  
\[
(\hat \Xcal_\ct, \hat D_\ct) \to (\Xcal_\ct, D_\ct).
\]
with Galois group $\Orth(\Tcal)$ and is obtained by  imposing  a $\Tcal$-level structure  on the objects it parametrizes. 
It  does not  ramify over the special mirror divisor $D^{\ct,\spl}$. 

Every connected component of the  normalized pull-back of this cover along the map
$\Xcal_{\DM}\to \Xcal_\ct$ (which itself defines the   normalization of the special mirror divisor)  is a reflection cover 
$(\hat \Xcal_{\DM}, \hat D_\DM)\to (\Xcal_{\DM}, D_\DM)$ has as Galois group a  reflection subgroup of $\Orth(\Tcal)$ of type $A_{11}$.
 \end{corollary}
\begin{proof}
Theorems \ref{thm:ct1} and \ref{thm:cubicperiods}  give us a universal reflection cover of 
the pair $(\Xcal_\ct\ssm D_{\ct, \spl}, D_\ct\ssm D_{\ct,\spl})$ with Galois group $\Orth(\Tcal)$. The normalization over $\Xcal_\ct$ yields an extension
as an  $\Orth(\Tcal)$-cover that we denote  $\hat \Xcal_\ct \to \Xcal_\ct$.

For the proof of the assertions regarding the chordal locus we use  Collino's analysis in \cite{collino}. He proves that for a general one parameter degeneration $\Xscr/\Delta$ over a complex disk  as in Lemma \ref{lemma:transversal} (so $X_t$ is the zero set of $F_0+\tilde\varphi(G)]$, where $G\in \CC[W]_{12}$), a base change $s=t^2$ trivializes the monodromy:
\[
\begin{CD}
\widehat \Xscr @>>> \Xscr\\
@V{\widehat f}VV @V{f}VV \\
\widehat \Delta @>{2:1}>> \Delta
\end{CD}
\]
Collino shows that the local system $R^3\widehat f_*\ZZ|_{\widehat \Delta\ssm \{o\}}$ is trivial
and hence extends canonically across the central point $o$. Let us denote the stalk over $o$ 
by $\Hl_o$. By construction, the antipodal map of $\widehat \Delta$ lifts to an involution  of this extension and hence acts in $\Hl_o$. Collino canonically identifies $\Hl_o$ with $\Hl^1(C)$, where $C$ is  the hyperelliptic curve of genus $5$ defined by $G$ with the involution being induced by the hyperelliptic involution. (This is even true, if we regard $\Hl_o$ as a limiting Hodge structure, but then of course a Tate twist is necessary.)  In particular, the monodromy of $R^3f_*\ZZ_\Xscr|_{\Delta\ssm \{o\}}$ is minus the identity.
It follows that a full level $2$-structure on $R^3f_*\ZZ_\Xscr|_{\Delta\ssm \{o\}}$ determines one on
$\Hl^1(C)$ (which in this case amounts to an identification of its monodromy group  with a reflection subgroup of $\Orth(\Tcal)$ of type $A_{11}$). Hence a connected component as in the statement of the theorem can be identified with the covering of  the moduli space of Hilbert-Mumford stable hyperelliptic curves of genus five defined by a full level 2 structure. %By Proposition \ref{prop:hypkernel} this gives us a universal reflection cover.
\end{proof}

\begin{remark}\label{rem:}
The normalization of $\Mcal_{\ct}^\st(A_1)$  is the  9-dimensional ball quotient of which an open subset parametrizes genus 4 curves. This ball quotient was investigated by K\=ondo \cite{kondo:g=4}. Via Lemma \ref{lemma:sconnectivity}
we obtain a  reflection cover of this ball quotient  whose generic point has the  modular interpretation
 as a genus 4 curve with a full level 2 structure on its Jacobian (so with  a copy of $\Sp_8(\FF_2)$ as covering group). It follows from Lemma \ref{lemma:sconnectivity} that this reflection cover is universal.  
\end{remark}

\section{Identification of coarse strata}\label{sect:coarsestrata}

The section is devoted to the proof of Theorem \ref{thm:match0} and will thus  complete the proof of our main Theorem \ref{thm:main}.

\subsection{Hodge structures and curves associated with a  singular stable cubic threefold}
The Hodge structure defined by a singular stable cubic threefold  is that  of  a curve. We check this in a special case.
It shows that  although $C_p$ may depend on $p$, its Hodge structure does not.

\begin{proposition}\label{prop:homologyidentification}
Let $X$ be a stable cubic threefold and $p\in X$ a singular point of type $A_3$ or $A_4$.
Then we have a natural isomorphism  $\Hl_3(X)\cong \Hl_1(C_p)(1)$ of polarized Hodge structures. 
\end{proposition}
\begin{proof}
For ease of notation,  suppress the subscript $p$  and write $\PP^3$ for $\PP(T_pV)$.
%Assume that  the situation of Proposition \ref{prop:converseconstruct}, $Q$ is a quadric cone  (so that $X=X_C$ comes with a singular point of type $A_3$ or $A_4$). 

Let $L_C\subset X$ be the union of the lines on $X$ through $p$ parametrized by $C$. Let $\tilde X\to X$ be the  blowup of $p$. 
We identify its exceptional divisor with $Q$. The strict transform $\tilde L_C$ of $L_C$ is a $\PP^1$-bundle. 
It meets $Q$ in $C$ and  hence $C$ can be regarded as a section of this bundle. The projection away from $p$ defines the morphism 
$\tilde X\to \PP(W)\cong \PP^3$ which amounts to the retraction $\tilde L_C\to C$  of the $\PP^1$-bundle onto its section.

The blowup/blowdown  gives isomorphisms
\[
\Hl_\pt (X, L_C)\xleftarrow{\cong} \Hl_\pt (\tilde X, Q\cup \tilde L_C)\xrightarrow{\cong} \Hl_\pt(\PP^3, Q).
\]
We make this substitution in the exact sequence of the pair  $(X,L_C)$ and get the exact sequence 
\[
\Hl_4(X)\to \Hl_4(\PP^3, Q)\to \Hl_3(L_C)\to \Hl_3(X)\to \Hl_3(\PP^3, Q)\to \Hl_2(L_C).
\]

The exact sequence of the pair $(\PP^3, Q)$ and the vanishing of $\Hl_3(\PP^3)=0$ imply
that $\Hl_3(\PP^3, Q)$ is the  kernel of $\Hl_2(Q)\to \Hl_2(\PP^3)$. Since $Q$ is a quadratic cone,  $\Hl_2(Q)$ is generated  by the fundamental class of a rule. This class maps to the generator of $\Hl_2(\PP^3)$ and so  $\Hl_3(\PP^3, Q)=0$.

We also find that  $\Hl_4(\PP^3, Q)$ is the cokernel of $\Hl_4(Q)\to \Hl_4(\PP^3)$ and so is 
of order $2$ and  generated by the class of a plane. Such a plane defines a hyperplane section of $X$ through $p$ which then 
maps with degree one  onto that plane. Hence  the map $\Hl_4(X)\to \Hl_4(X, L_C)\cong \Hl_4(\PP^3, Q)$ is onto. This implies that  $\Hl_3(L_C)\cong \Hl_3(X)$.  The  Thom isomorphism identifies $\Hl_3(L_C)$  with $\Hl_1(C)(1)$. 
\end{proof}

\begin{remark}\label{rem:qudraticiso}
This induces an isomorphism $\Hl_3(X; \FF_2)\cong \Hl_1(C_p; \FF_2)$ of symplectic vector spaces. This 
isomorphism takes vanishing cycles to vanishing cycles and is therefore in fact an isomorphism of quadratic spaces.
\end{remark}

Any plane  contained in $X$ and containing $p$ defines a rule of $Q_p$  that is an irreducible component of $C_p$ and conversely,
any  rule of $Q_p$  that is an irreducible component of $C_p$ determines such a plane. Let $F_p(X)\subset X$ denote the union of the planes in $X$ containing $p$ and let $R_p\subset  C_p$ denote the union of the rule components of $C_p$.  The proof of the following assertion is a minor  modification of the proof of Proposition \ref{prop:homologyidentification} and is therefore omitted.

\begin{proposition}\label{prop:homologyidentification1}
Let $X$ be a stable cubic threefold and $p\in X$ a singular point of type $A_3$ or $A_4$.
Then we have a natural isomorphism  $\Hl_3(X, F_p(X))\cong \Hl_1(C_p,R_p)(1)$ of Hodge structures. \hfill $\square$
\end{proposition}

\subsection{The coarse stratification}\label{subsect:coarse}
For the purposes of this paper (and as is also suggested by Corollary \ref{cor:coarse_emb}), we want to treat $A_1$ and  $A_2$-singularities of a cubic threefold as `not too serious'. This leads to: 

\begin{definition}\label{def:}
A \emph{coarse stratum} of $\Mcal^\st_{\ct}$ of type $(r,q)$ (where $r$ and $q$ are nonnegative integers) is a connected component  of a subvariety of the form
\[
\Mcal^*_\ct(rA_3+qA_4):=\cup_{k_1, k_2} \Mcal_\ct(k_1A_1+k_2A_2+rA_3+qA_4),
\]
A coarse stratum of $\Mcal^\st_{\ct}[\Tcal]$ is a connected component  of the preimage of a coarse stratum of $\Mcal^\st_{\ct}$.
\end{definition} 

The coarse strata clearly  partition $\Mcal_\ct^\st$.   
The difference 
$\Mcal^\st_\ct(rA_3+qA_4)\ssm \Mcal^*_\ct(rA_3+qA_4)$ is  of codimension $\ge 3$ in $\Mcal^\st_\ct(rA_3+qA_4)$ and so 
the inclusion $\Mcal_\ct(rA_3+qA_4)\subset \Mcal^*_\ct(rA_3+qA_4)$ induces  an bijection on connected components.
Proposition \ref{prop:stratabijection}  makes  the  connection with the eponymous notion  in the sense of  Corollary \ref{cor:coarse_emb}.

Corollary \ref{cor:stratalattices2} assigns to a  coarse stratum of $\Mcal^\st_\ct[\Tcal]$ a geometric subgroup $W\subset \Orth(\Tcal)$ such that this stratum is contained in the $W$-fixed point set. We shall later see (and this in fact one of the main assertions of Theorem \ref{thm:match0}) that it is open-dense in
$\Mcal^\st_\ct[\Tcal]^W$, so that $W$ determines the stratum and we are justified in denoting it by $\hat\Xcal\la W\ra$.
\smallskip

We have a similar coarse stratification of the chordal divisor $\Xcal_\DM$.  We identify this  with the GIT-moduli space of stable divisors of degree $12$, or rather,  of the  space $\Hyp_5^\st$ of hyperelliptic curves of genus $5$. This is much more straightforward: 
a stratum is given by a pair $(r,q)$  with  $4r+5q\le 12$ and describes the locus for which we have exactly $r$ 
Weierstra\ss\ points of multiplicity $4$, $q$ Weierstra\ss\ points of multiplicity $5$ and none  of  higher multiplicity.

We recall that  if $C$ is a hyperelliptic curve of genus $5$,
then an  identification of its (12 element) set $W_C$ of  
Weierstra\ss\ points with $E$ (let us call this an $E$-marking), 
determines  an  quadratic isomorphism $H^1(C; \FF_2)\cong \Tcal$ (let us call this a $\Tcal$-marking). 
The monodromy of the smooth hyperelliptic curves of genus 5  realizes the full permutation group of the set of Weierstra\ss\ points and hence  the moduli space of $E$-marked hyperelliptic curves is connected. 
It is a connected component, which  we will denote by $\Hyp_5[E]$,  of the moduli space $\Hyp_5[\Tcal]$ of all $\Tcal$-marked  smooth hyperelliptic curves of genus 5 and the associated group $\Sfrak_E$ can be understood  as its mod  2  monodromy group. We can phrase this as follows: 

\begin{corollary}\label{cor:}
The map which assigns to a connected component of the level $2$ cover $\Hyp_5[\Tcal]$ its mod  2  monodromy group sets up a bijection between these connected components and the reflection subgroups of $\Orth(\Tcal)$ isomorphic with $\Sfrak_{12}$.
\end{corollary}

We extend the cover $\Hyp_5[\Tcal]\to \Hyp_5$ to a (ramified) cover $\Hyp^\st_5[\Tcal]\to \Hyp^\st_5$ in the usual manner.
It has $\Hyp^\st_5[E]\to \Hyp^\st_5$ as a connected component. The following is then clear 

\begin{proposition}\label{prop:coarsestratahyp}
For every subgroup $W'\subset \Sfrak_E$ which is either geometric or trivial,  the fixed point set of $W'$ in
$\Hyp^\st_5[E]$ contains a unique  coarse stratum  as an open-dense subset and this defines a bijection  between  the set of such  subgroups of $\Sfrak_E$  and the coarse strata contained in $\Hyp_5^\st[E]$. The  Galois group for the reflection associated covering of that stratum over its image in $\Hyp_5^\st$  is the centralizer of $W'$ in $\Sfrak_E$. 
\hfill $\square$
\end{proposition}

We can now  prove the corresponding  part of Theorem \ref{thm:match0}.

\begin{proof}[Proof of the assertions regarding part (b) of Theorem \ref{thm:match0}]
Let $W_o\in\Cscr(A_3)$ and identify $W_o^\perp$ with $\Orth(\Tcal)$. Then $\Sfrak_E$  is of the form $\tilde W_o^\perp$, 
where  $\tilde W_o>W_o$ with $\tilde W_o\in \Cscr(A_4)$.  Any $W\in \Cscr_{\ge \tilde W_o}$ is given by a subgroup 
$W'< \tilde W_o^\perp\cong \Sfrak_E$ isomorphic with $\Sfrak_4^p\times \Sfrak_5^r$ for certain $p\ge 0$ and $r\ge 0$. 
This amounts to selecting from  $E$ a number of pairwise disjoint subsets of sizes $4$ ($p$ in number) and $5$ ($r$ in number): 
then $W'$ is the product of their permutation groups, $\out(W')\cong \Sfrak_p\times \Sfrak_r$ is the permutation group of these 
subsets which preserves size  and if $E'\subset E$ is the complement of the union of these subsets, then 
$W^\perp\cong \Sfrak_{E'}$. The stratum  $\Xcal\la W\ra=\tilde \Xcal\la W\ra$ is the 
$\PGL_2(\CC)\times (W'\times W'{}^\perp)\rtimes \out(W')$-orbit space of maps $E\to \PP^1$ that are constant on the given subsets  (giving  a prescribed  number of points of multiplicity $4$ and $5$; all other points have smaller multiplicity) and the cover $\hat \Xcal \la W\ra$ is defined by taking the  $\PGL_2(\CC)\times W'$-orbit space. The Galois group of
$\hat \Xcal \la W\ra\to \Xcal \la W\ra$ is  $W'\times \out(W')\cong \Sfrak_{E'}\times \Sfrak_p\times \Sfrak_r$.
\end{proof}

 \begin{small}
\begin{table}
\begin{centering}
\begin{tabular}{|c|c|c|c|c|}
\hline 
$\MM=\LL_4^0\oplus\NN$  &coarse   stratum of & $W<\Sfrak_{12}$ & centr.\ of& $\out_{\LL}(\MM)/$ \\
$\subset \LL$  &degree 12 divisors on $\PP^1$ &    &  $W$ in $\Sfrak_{12}$ & $\out_{\LL^\DM}(\NN)$\\
\hline\hline
$\LL^0_4 $ & $\Hcal^*_{12}$ & $\{1\}$ & $\Sfrak_{12}$ & $1/1$\\ 
\hline
$\LL^0_4\operp\LL_3 $  & $\Hcal^*_{12}(A_3)$ & $\Sfrak_4$ & $\Sfrak_8$ & $1/1$\\ 
\hline
$\LL^0_4\operp2\LL_3 $ & $\Hcal^*_{12}(2A_3)$ & $\Sfrak^2_4$ & $\Sfrak_4$ & $\Sfrak_2/\Sfrak_2$\\ 
\hline
$\LL^0_4\operp (3\LL_3)'$ & $\Hcal^*_{12}(3A_3)$ &  $\Sfrak_4^3$ & 1 &  $\Sfrak_3/\Sfrak_3$ \\ 
\hline
$\LL^0_4\operp \LL_4$ & $\Hcal^*_{12}(A_4)$ &  $\Sfrak_5$ & $\Sfrak_7$ &  $\Sfrak_2/1$\\ 
\hline
$\LL^0_4\operp \LL_4\operp \LL_3$ & $\Hcal^*_{12}(A_3+A_4)$ & $\Sfrak_4\times \Sfrak_5$ & $\Sfrak_3$ & $\Sfrak_2/1$\\ 
\hline
$\LL^0_4\operp 2\LL_4$ & $\Hcal^*_{12}(2A_4)$ & $\Sfrak_5^2$ & $\Sfrak_2$ & $\Sfrak_3/\Sfrak_2$\\ 
\hline
\end{tabular}
\caption{A  matching for the coarse strata of the chordal locus $\Hyp_5^\st$.}\label{table:refined2}
\end{centering}
\end{table}
\end{small}
\smallskip

The chordal locus is also involved in part (iii) of Theorem  \ref{thm:match0}. The following proposition takes care of that.

\begin{proposition}\label{prop:}
The reflection covers we obtain  for the coarse strata associated to a geometric lattices of the form  
$\MM=\LL^0_3\operp\LL_4\operp \NN$ and $\MM'=\LL^0_4\operp\LL_3\operp \NN$ are isomorphic.
\end{proposition}
\begin{proof}
We check this for  $\NN=0$ only, as the other  cases are specializations of this case.

According to Remark \ref{rem:modulariso}  both moduli spaces parametrize affine hyperelliptic curves of genus 3 with two points at infinity. 
If $X$ is a cubic threefold  with a singular point of type $A_3$, then Proposition  \ref{prop:homologyidentification}  identifies $\Hl_3(X)$ with 
$\Hl_1(C)$, where $C$ is the associated  hyperelliptic  genus 4 curve with a node.  Alexander duality then  gives identifications
\[
\Hl^3(X_\reg)\cong \Hl_3(X)\cong \Hl_1(C)\cong \Hl^1(C_\reg)
\]
and so a full level $2$-structure on $\Hl^3(X_\reg)$ amounts to a full level $2$-structure on $\Hl^3(C_\reg)$. Since $C_\reg$ is an affine hyperelliptic curve (of genus 3), this amounts to a numbering of its Weierstra\ss\ points. In terms of the equation 
$y^2=x^8+c_x^6+\cdots c_8$, this amounts to writing the right hand side as $(x-x_1)\cdots (x-x_8)$.  Hence the  reflection covers we get 
are isomorphic.
\end{proof}

In the rest of this section we therefore focus on  the stable cubic threefolds with  $A_3$-singularities at worst. For a finer  analysis, we need the notion of coplanarity, which we discuss next.

\subsection{Coplanarity and its homological incarnation} 
The following notion will help us to separate strata.

\begin{definition}\label{def:coplanar}
We say that two distinct singular points of a stable cubic threefold are \emph{coplanar} if both are contained in a plane contained in $X$.
\end{definition}

In case both  singular points are of type $A_3$, we shall relate this with  a homological property of the associated cubic fourfold. This will allow us to identify coarse strata of cubic threefolds in terms of geometric properties.

We begin with a general observation regarding $A_3$-singularities. 
Let  $p$ be a smooth or  isolated singular  point of a complex-analytic  variety $X$.  If (for some local-analytic  embedding of a neighborhood of $p$ in $\CC^N$) $B_p$  is  the  preimage of  a closed ball centered  at the image of $p$, then  when the radius of the ball is small enough,  $\Lk_p(X):=\p B_p$ is an oriented submanifold of $X$ whose isotopy class is independent of our choices (it  is called the \emph{link} of $p$ in $X$) and $B_p$ is homeomorphic to the cone over  $\p B_p$.

It is well-known that  if  $p$ is a hypersurface singularity of  type $A_3$ and of odd dimension $2m+1>0$, then $\Hl_{2m+1}(\Lk_p(X))$ is infinite cyclic. It is also known that if $(Y,o)$ is  of type $E_6$ and of even  dimension $2m+2$, then $\Hl_{2m+1}(\Lk_o(Y))$ is cyclic of order $3$.  These links are related in a simple manner as follows (we here take $m=1$, the only case that interests us).
An $A_3$-singularity $(X,0)\subset (\CC^4,0)$ has a local equation of the form  $f(z_1,z_2,z_3,z_4)=z_1^2z_2+z_2^2+z_3z_4$. It contains the plane $\Pi$ defined by $z_2=z_4=0$. Note that  $\Lk_0(\Pi)$ is a copy of $S^3$ and naturally oriented. So
$\Lk_0(\Pi)\cong \ZZ$. Let $(Y, 0)\subset (\CC^5,0)$ be the $4$-dimensional $E_6$ singularity defined by $f(z_1,z_2,z_3,z_4)+z_5^3$. Then in  
\[
\ZZ\cong \Hl_3(\Lk_0(\Pi))\to \Hl_3(\Lk_0(X))\to  \Hl_3(\Lk_0(Y))
\]
the first arrow is an isomorphism and the second is a surjection (so with  kernel $3\ZZ$). This remains true if we replace 
$\Pi$ by any smooth complex surface at $0\in \CC^3$ and contained in $(X,0)$.

\begin{corollary}\label{cor:linkhomology}
Let $X$ be a stable cubic threefold and let $(p,q)$ be coplanar pair of $A_3$ singularities of $X$. Then $\Hl_3(\Lk_p(X))$ and $\Hl_3(\Lk_q(X))$ have the same image in $\Hl_3(X\ssm\{p,q\})$ and if $Y$ is the associated cubic fourfold, then $\Hl_3(\Lk_p(Y))$ and $\Hl_3(\Lk_q(Y))$ have the same image in $\Hl_3(Y\ssm\{p,q\})$.
\end{corollary}
\begin{proof}
The difference of the links of $p$ and $q$ in $\Pi$ clearly bound in $\Pi\ssm \{p,q\}$, and hence also in $X\ssm \{p,q\}$ and $Y\ssm \{p,q\}$.
So these links define the same homology  class in degree 3  in these spaces.
\end{proof}

If in the situation of Corollary \ref{cor:linkhomology}, $(X_o, \g)$ defines a degeneration into $X$, then we have disjoint $\mu_3$-invariant Milnor fibers 
$\tilde F_p\subset Y_o$ and $\tilde F_q\subset Y_o$, whose homology in degree 4 determines  two copies  of $\LL_3$ in $\Hl_4(Y_o)$. 
The next corollary tells us something about the relative position of these two copies:

\begin{corollary}\label{cor:linkhomology2}
In the  situation of Corollary \ref{cor:linkhomology}
the  natural map $\Hl^4(Y_o)\to \Hl^4(\tilde F_p)\oplus \Hl^4(\tilde F_q)$ has cokernel of order $\ge 3$.
Equivalently, the associated embedding $2\LL_3\to \LL^\ct$ is not primitive: it has a cokernel of order $\ge 3$.
\end{corollary}
\begin{proof}
The path $\g$ identifies cohomology of the pair $(Y_o, \tilde F_p\cup \tilde F_q)$ with the homology of $(Y, \{p,q\})$ and so we have an exact sequence 
\[
\Hl^4(Y_o)\to  \Hl^4(\tilde F_p)\oplus \Hl^4(\tilde F_q)\to \Hl^5(Y, \{p,q\})%\to \Hl^5(Y_o)
\]
It suffices to show that the image of the last map has order $3$. For this we note that the  second arrow  factors as 
\begin{multline*}
\Hl^4(\tilde F_p)\oplus \Hl^4(\tilde F_q)\cong \Hl_4(\tilde F_p, \p \tilde F_p)\oplus \Hl_4(\tilde F_q, \p \tilde F_p)\twoheadrightarrow
\Hl_3(\p \tilde F_p)\oplus \Hl_3(\p \tilde F_q)\cong\\ 
\cong \Hl_3(\Lk_p(Y))\oplus \Hl_3(\Lk_q(Y))\to \Hl_3(Y\ssm  \{p,q\})\cong \Hl^5(Y, \{p,q\})
\end{multline*}
Here we used that a Milnor fiber of an $4$-dimensional isolated   hypersurface singularity  is $3$-connected (whence the first surjection) and that its  boundary can be identified with the link of that singularity.
By the previous Corollary \ref{cor:linkhomology} both $\Hl_3(\Lk_p(Y_o))$ and $\Hl_3(\Lk_q(Y_o))$ are of order three and have the same image in $\Hl_3(Y\ssm  \{p,q\})$. 
\end{proof}

We will  find  that  the converse also holds: if $p$ and $q$ are not coplanar, then  the  natural map $\Hl^4(Y_o)\to \Hl^4(\tilde F_p)\oplus \Hl^4(\tilde F_q)$ is onto (implying  that $2\LL_3$ embeds primitively in $H_4(Y_o)$). So coplanarity admits in fact a homological characterization. 
\smallskip

In the rest of this subsection, we assume that $p\in X$ is a singular point of type $A_3$  or $A_4$, so that  $Q_p$ a quadric cone. 
The following lemma shows that if $p\in X$ is a singular point of type $A_3$ or $A_4$, then the  planes contained in $X$ are often `visible' in $C_p$.

\begin{lemma}\label{lemma:planerec}
The planes passing through  $p$ correspond to the rules of $Q_p$ contained in $C_p$.
Moreover,  an irreducible component of $C_p$ that is a plane section of $Q_p$ determines a plane in $X$ that does not pass through $p$, 
 and the other irreducible components correspond to plane sections of $Q_p$ contained in $C_p$  not passing through the vertex of $Q_p$ (and all such irreducible components so arise).
\end{lemma}
\begin{proof}
The bijection regarding planes in $X$ containing $p$ is clear and so we focus on last assertion.

Suppose that $C_p$ has as an irreducible component a plane section $Z$ of $Q_p$ which does not pass through the origin.
We choose affine coordinates $(x,y,z, t)$ in our $\PP V$ such that $p$ is the origin and $X$ is given by an equation of the form
$F(x,y,z,t)=xy+z^2 +f_3(x,y,z,t)$ with $f_3$ homogeneous of degree $3$. Then $C_p\subset \PP^3$ is the complete intersection given
by $xy-z^2=0$ and  $f_3(x,y,z,t)$. After perhaps modifying our coordinates we can arrange that $Z$ is the hyperplane section at infinity in $\PP^3$. This means that
 $f_3$ is in the ideal generated by $xy+z^2$ and $t$, so can be written as $(xy+z^2)g_1+ tg_2$ with $g_i$ homogeneous of degree $i$. Then 
 $F(x,y,z,t)=(xy+z^2)(1+g_1)+ tg_2$ and hence vanishes of the affine plane defined by $1+g_1=0$ and $t=0$.
 \end{proof}

We  now determine   what the irreducible components of  $C_p$ can be like. The blowup of its  vertex of $Q_p$ defines a Hirzebruch surface $\tilde Q_p$ with exceptional section $E$ and  comes canonical retraction of $\tilde Q_p$ onto $E$. The class $e$ of $E$  and the class $f$ of a fiber of $\tilde Q_p\to E$ make up a basis of $\pic(\tilde Q_p)=H^2(\tilde Q_p)$. Then the  strict transform $\tilde C_p\subset \tilde Q_p$ of $C_p$ has class  $2e+6f$. Since $(2e+6f)\cdot f=2$, the projection $\tilde C_p\to E$ is of  degree 2
and  this endows  (as mentioned in Proposition \ref{prop:curveCp}) $\tilde C_p $ with its hyperelliptic involution $\iota$ . The  arithmetic genus of  $\tilde C_p$  is $3$.

\begin{lemma}\label{lemma:classes}
Suppose $\tilde C_p$ is reducible and let $Z$ be an irreducible component  of $\tilde C_p$ whose class in $\pic(\tilde Q_p)$ we denote $[Z]$.
Then we are in one of the following cases 
\begin{description}
\item[{$[Z]=f$}]  $Z$ is a rule (i.e., a fiber of the retraction $\tilde Q_p\to E$).
\item[{$[Z]=e+2f$}] $Z$ is  a plane section of $Q_p$ which does not pass through the vertex (so a smooth conic and hence by Lemma \ref{lemma:planerec} represents a plane in $X$ which does not pass through $p$).
\item[{$[Z]=e+3f$}] $Z$ is a section of the retraction $\tilde Q_p\to E$ which meets $E$ with multiplicity one (such sections can be identified with sections
of $\Ocal_E(3)$ and hence are parametrized by a projective space of dimension 3).
\item[{$[Z]=e+4f$}] $Z$ is a section of the retraction $\tilde Q_p\to E$ which meets $E$ with multiplicity two  (such sections can be identified with sections of $\Ocal_E(4)$ and hence are parametrized by a projective space of dimension 4).
\item[{$[Z]=2e+4f$}] $Z$ is an  intersection of  $Q_p$ with a quadric which does not pass through the vertex; such a $Z$ is of degree $2$ over $E$ (and hence comes with a  hyperelliptic involution $\iota$) and has arithmetic genus $1$.
\item[{$[Z]=2e+5f$}] $Z$ is of degree $2$ over $E$ (and hence comes with a hyperelliptic involution $\iota$ of $Z$), meets $E$ transversally in a single point and its  arithmetic genus is $2$.
\end{description}
\end{lemma}
\begin{proof}
If $[Z]=a e+ bf$, then $(a, b)\in \{0,1,2\}\times \{0, 1,\dots , 6\}$, but $(a, b)\not=(2,6)$. We know  that $0\le [Z]\cdot e= -2a+b$ and if we write 
$\tilde C$ as $Z+\tilde C'$, then $0\le [\tilde C']\cdot E = ((2-a)e+(6-b)f)\cdot e=2+2a-b$. This implies that $[Z]$ must be one of the classes listed. It is straightforward to show that
the  irreducible  curves realizing these classes have  the stated form.
\end{proof}

Assume that $X$ has only $A_3$-singularities with $p\in X$ being one of these. We list all the possibilities for the curve $C_p$.
So $C_p$ has  at the vertex a node and its other singularities are all of type $A_3$.  Proposition \ref{prop:converseconstruct} shows that  
this curve $C_p$ is a complete invariant of  the pair $(X,p)$. The description that we give shows that in each  case such curves are 
parametrized by  an irreducible subvariety of the universal stable cubic threefold (we indicate this in the notation by the subscript 
``$1$'' in $\Mcal_1$). We shall later see how to extract from this a complete invariant of $X$, in other words, find its 
image under the   projection  to $\Mcal_\ct$ (which will then be open-dense in some coarse stratum). 
This will then also justify the chosen notation below.
\smallskip

The normalization $\hat C_p\to C_p$  factors as   $\hat C_p\to \tilde C_p\to C_p$, where  we recall that $\tilde C_p\subset \tilde Q_p$ is the strict transform  of $C_p$.  
Assume that $X$ has only $A_3$-singularities (with $p$ being such a singularity). With the help of  Lemma \ref{lemma:classes}, 
we can now determine what $C_p$ can be like (the notation indicates the number of $A_3$-singularities).
\smallskip

$\Mcal_1(A_3)$: $C_p$ has an ordinary double point at the vertex and  $\tilde C_p$ is  a smooth hyperelliptic curve of genus 3.
The preimage of the node is a 2-element set $\Pcal\subset \tilde C_p$ which is not invariant under the hyperelliptic involution.
The pair $(\tilde C_p; \Pcal)$ is a complete invariant of $(X,p)$, and all such pairs arise. %There is no plane contained in $X$. 

\smallskip
We next do the case of two $A_3$-singularities, so that $C_p$   has just one singular point $q_p$. It must be  invariant under $\iota$.

\smallskip
$\Mcal_1(2A_3)$: \emph{$C_p$ is irreducible.}\\
 Then $C=\hat C_p$ is connected of  genus
$1$  and comes with two 2-element subsets:  the fiber $\Pcal$ over the vertex $v$ and  the fiber $\Qcal$ over $q_p$.
The involution $\iota$ on $C_p$ can be reconstructed from this, because it is characterized by the property that $\Qcal$ is an orbit. 
This makes  the triple $(C; \Pcal, \Qcal)$  a complete invariant of $(X,p)$.  
Its moduli space is irreducible.
If we let the two points of $\Pcal $ resp.\  $\Qcal$ coalesce we get at $p$ resp.\ $q$ an $A_4$-singularity. 
%There is no plane contained in $X$. 

\smallskip
$\Mcal_1(2A_3)'$: \emph{$C_p$ is reducible.}\\ 
Then we find that one component of $C_p$ must be a rule $F\subset \tilde Q_p$  and the other component a smooth genus 2 curve $C$ of class $2e+5f$ with $C\to E$ defining its  hyperelliptic involution $\iota$. 
The  curve $C$ meets $E$ in a singleton $v$ and $F$ at the point $q_p$, where it is tangent to $F$ to produce an $A_3$-singularity.  
This implies that $q_p$ is a Weierstrass point on $C$. The triple $(C; v, q_p)$ is a complete invariant. 
The rule represents is  a plane contained in $X$  passing through $p$.
\smallskip

We now turn to the case $k=3$. So $C_p$ has two $A_3$ singularities, $q_p$ and $ q'_p$ say. 
Since $\tilde C_p$ is a hyperelliptic curve of arithmetic genus $3$ with two $A_3$-singularities, it  must be reducible.
\smallskip

$\Mcal_1(3A_3)$:  \emph{$\tilde C_p$ splits into two irreducible components, each of  class $e+3f$.}
So each  irreducible component is given by a section of $\tilde Q_p\to E$  meeting  $E$ transversally once. 
Since $(e+3f)\cdot (e+3f)=4$, the two components must meet in $q_p, q'_p$ and have there a common tangent 
(giving an $A_3$-singularity). %There is no plane contained in $X$. 

\smallskip 
$\Mcal_1(3A_3)'$: \emph{$\tilde C_p$  has two rules $F$ and $F'$ as an irreducible component, also accounting for an $A_3$-singularity.}\\
Then the residual curve $C$ (of class $2e+4f$) is of genus one and cut out by a quadric not passing through the vertex. It must be tangent to both $F$ and $F'$.
So $C$ is  a smooth genus one curve, comes with an involution $\iota$ with two fixed points $q_p, q_{p'}$ singled out; if we choose $q_p$ to be the origin, then $q_{p'}$ is of order 2 and so the isomorphism type $(Q_p, C_p)$ or equivalently $(X,p)$  is determined by a pair consisting of an elliptic curve plus a point of order $2$.
The two rules determine the two planes contained in $X$ each containing $p$ and  another  singular point. 

\smallskip
There are two  remaining cases. These will turn out to define the same family of cubic threefolds (it is here the choice of $A_3$-singularity $p$  which makes the difference).

\smallskip
$\Mcal_1(2A_3)'+A_3)$-(i): \emph{$\tilde C_p$ has  two irreducible components, one of which is a rule $F$.}\\
The residual component $\tilde C'_p$ has class $2e+5f$, is tangent to the rule $F$ (in a point  $q_p\in F\ssm E$), is of degree $2$ over $E$ and has 
an $A_3$-singular point $q_{p'}$. Its normalization of $C$ is of genus zero and comes with four points: its intersection with $E$ (which we may associate with $p$), its intersection with $F$ (which we may associate with $q$) and the unordered pair over the singularity (which we may associate with $q'$). The rule defines  plane  in $X$ containing $p$ and another singular point.

The roles of $p$ and $q$ are here symmetric, but if we place ourselves at the $A_3$-singular point $q'$, then the situation is different. 
This is in fact  accounted for by the remaining case: 

\smallskip
$\Mcal_1((2A_3)'+A_3)$-(ii): \emph{$\tilde C_p$ has  two irreducible components, one of which is a plane section of $Q_p$ not passing through the vertex.}\\
The class of this plane section is then $e+2f$ and the class of the residual component therefore $e+4f$. Both are given by sections of $\tilde Q_p\to E$,
the former does not meet $E$ and the latter, our $C$,  meets $E$ in two distinct points. The two components meet in two distinct points $q_p$, $q_{p'}$ and meet there with multiplicity $2$, so that their union exhibits at each of these  a $A_3$-singular point. We get again four points on $E$ associated with the three $A_3$-singularities, but now $q$ and $q'$ are associated with singletons and $p$ is associated with the unordered pair. The plane section accounts for  plane contained in $X$ containing the two singular points other than $p$.

\smallskip
\begin{proposition}\label{prop:}
Coplanarity is a transitive property on the set of singular points (all assumed to be of type $A_3$) of $X$.
\end{proposition}
\begin{proof}
In view of our list above, this is something we only need to check for the case $\Mcal_1(3A_3)'$: we must show that $X$ contains a plane which contains the two singular points not equal to $p$. Choose affine coordinates $(x,y,z,w)$ for $\PP(V)$  such that$p$ is the origin 
and a defining equation for $X$ has the form $x^2+yz + f_3(x,y,z,w)$ with $f_3$ homogeneous of degree $3$ and the two singular points other than $p$  are on the lines $x=y=w=0$ and  $x=z=w=0$ (so in the $\PP^3$ with homogeneous coordinates 
$[x:y:z:w]$ the quadric cone $Q_p$ is given by $x^2+yz =0$ and $C_p$ is its intersection with the cubic surface defined by $f_3$).
The hyperplanes contained in $X$ must be given by $x=y=0$ and  $x=z$ and hence  $f_3$ must lie in the intersection of the ideals generated by $(x,y)$ and $(x, z)$. This is ideal generated by $(x, yz)$, which is also the ideal generated by $x^2+yz$ and $x$. If we write $f_3$ accordingly as $g_1(x^2+yz) +x g_2$ with $g_i$ homogeneous of degree $i$, then
$f_3=(1+g_1)(x^2+yz) +x g_2$ and so vanishes on the affine plane $\Pi$ defined by $1+g_1=0$ and $x=0$. 
This plane meets the line $x=y=w=0$ resp.  $x=z=w=0$  where $g_1(0,0,z,0)=-1$ (or in projective coordinates, where 
$g_1(0,0,z,0)+t=0$). These intersection  points are singular points of $X$.
\end{proof}

We can now match these cases with our classification in terms of geometric sublattices:

\begin{proposition}\label{prop:A3matching}
The cases listed above match the $\U(\LL_\ct)$-orbits of geometric sublattices of $\LL$ of the form $\LL_3^0\operp \NN$ with 
$\NN\cong \LL_3^r$ and are given in Table \ref{table:latticeemb} in the manner suggested by the notation.  Table \ref{table:refined}  makes this precise.
\end{proposition}
\begin{proof}
The curve $C_p$ has a node at the vertex of $Q_p$, but all other singular points must be of type $A_3$. 
 We must  match the six cases found above with the six items  listed in Table \ref{table:latticeemb}. The map which assigns to one of our strata an item in the list in Table \ref{table:latticeemb} 
must be onto and hence will be a bijection. According to Corollary \ref{cor:linkhomology2} this  must be such that a coplanar pair of $
A_3$-singularities yields two $\LL_3$-summands  of the associated geometric lattice  which are not primitively embedded. 
This forces what the images of $(2A_3)'$, $(3A_3)'$ and
$(2A_3)'+A_3$ will be. This then determines the images of the remaining cases $A_3$, $2A_3$ and $3A_3$.
\end{proof}

 \begin{small}
 \begin{table}
\begin{centering}
\begin{tabular}{|c|c|c|c|c|}
\hline 
$\MM=\LL_3^0\operp\NN$ &  coarse stratum & type of $W$ & type of $W^\perp$  &  $\out_{\LL^\ct}(\NN)\subset $\\
$\subset \LL$  & in $\Mcal^\st_\ct$ &in $\Gcal$  & in $\Orth(\Tcal)$  & $\out_{\LL}(\MM)$ \\
\hline\hline
 $\LL^0_3 $ & $\Mcal^*_\ct$ & $A_3$  & $\Orth^1_{10}(\FF_2)$ & $\{1\}=\{1\}$\\ 
\hline
 $\LL^0_3 \operp\LL_3$ & $\Mcal^*_\ct(A_3)$ & $2A_3$ & $\FF_2^6\rtimes \Sfrak_8 $& $\{1\}\subset \Sfrak_2$\\
 \hline
 $\LL^0_3 \operp\LL_4$ & $\Mcal^*_\ct(A_4)$ & $A_3+A_4$ & $\Sfrak_8$ & $\{1\}=\{1\}$\\
\hline
 $\LL^0_3 \operp 2\LL_3$  & $\Mcal^*_\ct(2A_3)$ &  $3A_3$ & $\Hfrak(A_2)\rtimes \Sfrak_3$& $\Sfrak_2\subset \Sfrak_3$\\
 \hline
 $\LL^0_3 \operp\LL_3\operp\LL_4$ & $\Mcal^*_\ct(A_3+A_4)$ & $2A_3+A_4$ & $\Sfrak_4$ & $\{1\}\subset \Sfrak_2$\\
\hline
 $\LL^0_3 \operp2\LL_4$ & $\Mcal^*_\ct(2A_4)$ & $A_3+2A_4$ & $\Sfrak_3$ & $\{1\}=\{1\}$\\
 \hline
 $(\LL^0_3 \operp 2\LL_3)'$  & $\Mcal^*_\ct(2A_3)'$ & $(3A_3)'$ & $\FF_2^4\times \Sfrak_5$&$ \Sfrak_2\subset \Sfrak_3$\\
 \hline
 $\LL^0_3\operp(3\LL_3)'$ & $\Mcal^*_\ct(3A_3)'$ & $A_3+(3A_3)'$& $\FF_2^2\rtimes\Sfrak_2$ & $ \Sfrak_3=\Sfrak_3$\\
 \hline
$(\LL^0_3\operp 2\LL_3)'\operp \LL_3$ & $\Mcal^*_\ct((2A_3)'+A_3)$ & $(3A_3)'+A_3$ & $\FF_2$&$\Sfrak_2\subset  \Sfrak_3$\\
 \hline
$(\LL^0_3\operp 2\LL_3)'\operp \LL_4$ & $\Mcal^*_\ct((2A_3)'+A_4)$ & $(3A_3)'+A_4$ & $1$ &$ \Sfrak_2\subset \Sfrak_3$\\
\hline
 $(\LL^0_3\operp 3\LL_3)'$ & $\Mcal^*_\ct(3A_3)$ & $(4A_3)'$ & $\FF^2_2\rtimes \Sfrak_2$&$\Sfrak_3\subset \Sfrak_4$\\
\hline
\end{tabular}
\caption{The matching for the coarse strata in  $\Mcal^\st_\ct[\Tcal]$.}\label{table:refined}
\end{centering}
\end{table}
\end{small}

\subsection{Characteristic  curves associated with singular cubic threefolds}
We here determine  among other things how the curves $C_p$ are related if $p$ runs over the singular points of $X$ of type $A_3$ or $A_4$. This will also explain why they have the same Hodge structure. 
\smallskip

The \emph{weak  normalization}  $\hat C^w_p\to C_p$  replaces each $A_\even $-singularity by a smooth point  and each $A_\odd$-singularity by a normal crossing  (an $A_1$-singularity). So $C^w_p$ is a normal crossing curve and  the  weak normalization is a homeomorphism (and is universal for that property). 
It therefore induces an isomorphism on homology. The irreducible components
of $C^w_p$ which are  a rule of $Q_p$ or a smooth plane section of $Q_p$ are smooth copies of $\PP^1$ which meet the rest of the curve in two points at most.  Let $C^w_p\to \check C_p$ be defined by contracting these components. This is still a normal crossing curve and the contraction map defines an isomorphism on $\Hl_1$. The singular points of $C_p$ of type $A_2$ or $A_4$  define a subset 
$\check P_p$ of the smooth part of $\check C_p$ and the pair $(\check C_p, \check P_p)$ is Deligne-Mumford stable (its automorphism group is finite). Note that this contraction procedure causes singular points of $C_p$ that correspond to coplanar singularities  of $X$ to have the same image.

The goal of this subsection  is to prove:

\begin{proposition}\label{prop:CX}
Let $X$ be a stable cubic surface with at least one  singular point of type $A_3$ or $A_4$. Then 
with  $X$ is canonically associated a Deligne-Mumford  stable pointed curve  $(C_X,P_X)$ with the property that for every singular point $p\in X$ of type $A_3$ or $A_4$ we  have a  canonically identification with $(\check C_p, \check P_p)$.  Furthermore,  $H_3(X)$ is  as a polarized Hodge structure  canonically isomorphic with $H_1(C_X)(1)$. 
\end{proposition}

Suppose $p, q\in X$ are distinct singular points of $X$ of type $A_3$ or $A_4$. 
Since the line $\la p, q\ra$  in $\PP(V)$ has  intersection multiplicity $\ge 2$ with $X$ at $p$ and $q$, it must lie on $X$ and hence  defines a point $q_p\in C_p$ 
(and likewise a point $p_q\in C_p$).  The following implies  that $q_p$ must be a smooth point of $Q_p$. 

\begin{lemma}\label{lemma:smoothpart}
The line in $\PP V$ through $p$ defined by the vertex of this cone meets $X\ssm \{p\}$ in its smooth part.
\end{lemma}
\begin{proof}
Assume that this line contains a singular point $q\not=p$ of $X$. As we just observed, this implies that this line is contained in $X$. Choose an affine coordinate system $(z_1,z_2,z_3,z_4)$
for $\PP V$  such that  $p=(0,0,0,0)$, the vertex of the cone is given by the $z_4$-axis and $q=(0,0,0,1)$. Then $X$ is given by an equation of the form
\[
F(z_1,z_2,z_3, z_4)=f_2(z_1,z_2, z_3) + g_3(z_1,z_2, z_3)+g_2(z_1,z_2, z_3)z_4+ g_1(z_1,z_2, z_3)z^2_4,
\]
with $f_2,g_1,g_2,g_3$ homogeneous of the degree given by the subscript. Since  $(0,0,0,1)$ is a singular point, $g_1=0$ must be identically zero.
But then $X$ is singular along the $z_4$-axis and hence not GIT stable. 
\end{proof}

\begin{lemma}\label{lemma:mult2}
Let $\Pi$ be a plane in $\PP V$ containing the line $\la p,q\ra$. Then the divisor  on $\Pi$ cut out by $X$ contains  $\la p,q\ra$ with multiplicity $>1$
if and only if $\Pi$ defines a line in $\PP T_pV$ that is tangent to $Q_p$  at $q_p$ 
(and a line in $\PP T_qV$ that is tangent to $Q_q$ at $p_q$).
\end{lemma}
\begin{proof}
Choose affine coordinates $(z_1, \dots, z_4)$ for $\PP V$ for which $p$ is the origin,  $\la p,q\ra$ is defined by $z_2=z_3=z_4=0$ and $\Pi$ 
is defined by $z_3=z_4$. Then $X$ is given by an equation of the form $f_2+f_3$ with $f_i$ homogeneous of degree $i$, $Q_p\subset \PP^3$ is the quadric defined by $f_2$ and $q_p=[1:0:0:0]$.  The property that  the divisor  on $\Pi$ cut out by $X$ contains  $\la p,q\ra$ with multiplicity $>1$
amounts to  $f_2(z_1,z_2, 0,0)+ f_3(z_1,z_2, 0,0)$ being divisible by $z^2_2$. It is clear that divisibility property means that $f_2$ is proportional with
$z^2_2$ and $f_3$ is  divisible by $z^2_2$. So 
the line  in  $\PP^3$ defined by $z_3=z_4=0$ is tangent to $Q_p$  at $[1:0:0:0]$.
Conversely, if the line  in  $\PP^3$ defined by $z_3=z_4=0$ is tangent to $Q_p$  at $[1:0:0:0]$, then $f_2(z_1,z_2, 0,0)$ is a multiple of $z_2^2$.
Since $q_p$ is a singular point of $C_p$, $f_3(z_1,z_2, 0,0)$ is then also divisible by $z_2^2$ and hence so is $f_2(z_1,z_2, 0,0)+ f_3(z_1,z_2, 0,0)$.
\end{proof}

With  $p$ and $q$ as above (so distinct singular points of the stable $X$ of type $A_3$ or $A_4$), let  $\Pi$ be a  plane in $\PP V$ containing $\la p, q\ra$. Such a plane defines  a line in $\PP(T_pV)$ through $q_p$ and vice versa. Because of the symmetry this also amounts to specifying a line through  $p_q$ in $\PP(T_qV)$. 

The most degenerate case is when  $\Pi\subset X$. Then $\Pi$ determines a rule of $Q_p$ contained in $C_p$.  For the same reason, $\Pi$ determines (and is determined by) a rule in $Q_q$ contained in $C_q$. There seems however  to be no natural way to identify these two rules.

We therefore assume  that $\Pi\not\subset X$. Then $\Pi\cap X$ is the sum of $\la p, q\ra$ and  a residual conic in $\Pi$. If  this conic is singular, then it is the sum of two lines $\ell+\ell'$ with $p\in \ell$ and  $q\in\ell'$. These correspond to points $z\in C_p$ and $z'\in C_q$.

Conversely,  if  $z\in C_p$ and does not lie on the rule spanned by $q_p$,  then the line $\la z, q_p\ra$ defines a plane $\Pi_z\subset \PP(V)$ through $\la p, q\ra$  which is 
not contained in $X$ and $\Pi_z\cap X$  has as components the line through $p$ defined by the vertex $v_p$, the line $\la p, q\ra $ and 
another line through $q$, which then defines a point  of $C_q$. We denote that point $f^p_q(z)$. 

The proof of Proposition \ref{prop:CX} is now completed by:

\begin{lemma}\label{lemma:pqiso}
Define $C'_p(q)\to C_p$  as the  weak normalization of $C_p$ over the vertex $v_p$ of $Q_p$ and $q_p$ (so that is a homomorphism)
and let $C'_p(q)\to C_p(q)$ contract the preimage of the rule spanned by $q_p$ if that is an irreducible component of $C_p$ and be otherwise the identity. Then the above assignment defines an isomorphism  $f^p_q$ of $C_p(q)$ onto $C_q(p)$ which takes the images of $q_p$ resp.\  
 $v_p$ to the image of $v_q$ resp.\ $p_q$ (we keep on denoting the preimages of points of $C_p$ in  $C'_p(q)$ by the same symbol).
\end{lemma}
\begin{proof}
The construction gives a morphism $f^p_q: C_p\ssm\la v_p,q_p\ra \to C_q$ and shows that $f^p_q$ and $f^q_p$ are each others inverse, wherever  that makes sense.  

Let us first observe that the line  defined by the rule spanned by $q_p$ defines a plane $\Pi$  and this plane  is either contained in $X$ (this is equivalent to the rule spanned by $q_p$ being  contained in $C_p$) or,  by Lemma \ref{lemma:mult2}, $X$ cuts out on $\Pi$ the line $\la p,q\ra$ with multiplicity 2 plus  the line  defined by $v_p$. 

The branches of $C_p$ at $q_p$ have at $q_p$ the same tangent line. This tangent line is a rule if and only if one of the branches is a rule.
If  this tangent line is not  rule, then it defines a plane $\Pi$ and since $X\cap \Pi$ is twice the line  $\la p,q\ra$ plus  the line  defined by 
$v_p$, it follows that   $f^p_q(q_p)$ is then represented by the line $\la p,q\ra$, in other words,  is equal to $p_q$. 
If this tangent line is a rule, then $\Pi$ is contained in $X$. This plane  is in $C_q$ represented by the rule 
spanned by $p_q$. So  the limiting value of $f^p_q(q_p)$ must lie on that rule. This gives therefore  a well defined value in $C_q(p)$.

The discussion for the branches of $C_p$ at the vertex $v_p$ is similar, for then  the limiting line will be the rule spanned by $q_p$ and hence the limit plane is $\Pi$.  If $\Pi$ is not contained in $X$, then we see that $f^p_q(v_p)=p_q$. Otherwise it will be a point of the rule spanned by $p_q$ which is eventually contracted.
\end{proof}

\subsection{A stable model for the curve $C_p$} Let $X$ be a stable cubic threefold and let $p\in X$ be a singular point. The  curve $C_p$ has a topological \emph{stable model}: this is what a `nearby' smooth genus five curve with its Milnor fibers would look like:  a closed connected  topological surface $\Sigma$ of genus $5$ with  a  subsurface with boundary $\Sigma_o\subset \Sigma$  such that
 each connected component of $\Sigma_o$  is the  Milnor fiber of a singularity  that appears on $C_p$: an  $A_i$-singularity giving a subsurface of genus 
 $\lfloor i/2\rfloor$ with one or two boundary components according to whether $i$ is even or odd.  
 The subsurface $\Sigma_o$ may have connected components that are cylinders; we then render the situation  `stable' by shrinking such a cylinders to an  embedded  circle. The union $A$ of these  embedded circles and the other boundary components of $\Sigma_o$ 
is then  such that a connected component of  $\Sigma\ssm A$ represents  either a Milnor fiber or a  connected component of 
 the smooth part of $C_X$. We associate with this pair $(\Sigma,A)$ a \emph{(stable) weighted graph} in the  usual way: vertices  are the connected components of  $\Sigma\ssm \A$ weighted by their genus and two vertices  are connected by as many bonds as the associated connected components  have boundary components in common.  So a vertex that represents  an $A_i$ singularity (we shall call  such a vertex a \emph{virtual vertex}) has weight $\lfloor i/2\rfloor$ and degree $1$ or $2$  according to whether $i$ is even or odd. Among them is the vertex defined by the singularity $p$: we call this  the  \emph{distinguished virtual vertex}.
 Any other vertex (which we shall call a \emph{geometric vertex})  represents a connected component of  the normalization of $C_p$: in the graph we indicate this by its weight. The sum of the vertex weights and the first Betti number of  the graph  must add up to the arithmetic  genus of $C_p$ (which is $5$).  

The weighted graph with its distinguished virtual vertex, thus defined  will be denoted $\G(X,p)$. It is  useful  combinatorial  invariant of the pair  $(X,p)$. If $(X,p)$ is a degeneration of  $(X', p')$, then $\G(X', p')$ is obtained from $\G(X, p)$ by a successive application of bond contractions: a bond contraction is what happens if we remove a connected component of $A$: this has the effect on the graph by identifying the vertices connected by the bond and letting its the weight be such that the arithmetic  genus stays 5.

The classification of the possible types of $C_p$  in the   previous section and Proposition \ref{prop:A3matching} then gives:
\smallskip

\begin{proposition}\label{prop:ka3}
Assume $X$ is a  singular stable  cubic threefold  with only  $A_3$-singularities and let $p$ be one such singularity. 
Then the stable graph $\G(X,p)$ is as below, where 
the labeling is  as in Table \ref{table:latticeemb}.

\begin{center}
\begin{tiny}
$A_3$
\begin{tikzpicture}
\begin{scope}[every node/.style={circle,draw}]
      \node (X) at (-1.5,0) {$A_3$};
    \node (Y) at (0,0) {$3$};
 \end{scope}     
\begin{scope}[={[black]}, every edge/.style={draw=black}]  
 \path  (X) edge[bend right=10] node {} (Y);
    \path  (X) edge[bend left=10] node {} (Y);
   \end{scope}
\end{tikzpicture}
\hskip4mm
$2A_3$
\begin{tikzpicture}
\begin{scope}[every node/.style={circle,draw}]
      \node (X) at (-1.5,0) {$A_3$};
    \node (Y) at (0,0) {$1$};
     \node (Z) at (1.5,0) {$A_3$};
 \end{scope}     
\begin{scope}[={[black]}, every edge/.style={draw=black}]  
 \path  (X) edge[bend right=10] node {} (Y);
    \path  (Z) edge[bend left=10] node {} (Y);
    \path  (Z) edge[bend right=10] node {} (Y);
    \path  (X) edge[bend left=10] node {} (Y);
   \end{scope}
\end{tikzpicture}
\hskip4mm $(2A_3)'$
\begin{tikzpicture}
\begin{scope}[every node/.style={circle,draw}]
     \node (X) at (-0.5,0) {$A_3$};
     \node (Y) at (0.5,0) {$A_3$};
     \node (Z) at (0,1) {$2$};
 \end{scope}     
\begin{scope}[={[black]}, every edge/.style={draw=black}]  
    \path  (X) edge node {} (Z);
    \path  (Y) edge node {} (Z);
    \path  (X) edge node {} (Y);
 \end{scope}
  \end{tikzpicture}
\smallskip

\hskip4mm 
$3A_3$
\begin{tikzpicture}
\begin{scope}[every node/.style={circle,draw}]
    \node (A) at (0,1)  {$0$} ;
    \node (B) at (-1,0)  {$A_3$} ;
    \node (C) at (0,0)  {$A_3$} ;
    \node (D) at (1,-0) {$A_3$};
    \node (E) at (0,-1)  {$0$} ;
 \end{scope}     
\begin{scope}[={[black]}, every edge/.style={draw=black}]  
  \path  (A) edge node {} (B);
    \path  (A) edge node   {}(C);
    \path  (A) edge node   {} (D);
    \path  (E) edge node {} (B);
    \path  (E) edge node {} (C);
    \path  (E) edge node {} (D);
 \end{scope}
  \end{tikzpicture}
$(3A_3)'$
\begin{tikzpicture}
\begin{scope}[every node/.style={circle,draw}]
    \node (A') at (0,0.7)  {$1$} ;
    \node (B') at (-0.7,-0)  {$A_3$} ; 
    \node (D') at (0.7,0) {$A_3$};
    \node (E') at (0,-0.7)   {$\mathbf{A_3}$};
 \end{scope}

\begin{scope}[={[black]},  every edge/.style={draw=black}]
     \path  (A') edge node {} (B');
    \path  (A') edge node   {} (D');
    \path  (E') edge node {} (B');
    \path  (E') edge node {} (D');
\end{scope}
     \end{tikzpicture}   
 $(2A_3)'+A_3$
 \begin{tikzpicture}
\begin{scope}[every node/.style={circle,draw}]
  \node (A'') at (-1,0.5)  {$A_3$} ;
    \node (B'') at (-1,-0.5)  {$A_3$} ;
    \node (C'') at (0,0)  {$0$} ;
    \node (D'') at (1.5,0) {$A_3$}; \end{scope}     
\begin{scope}[={[black]}, every edge/.style={draw=black}]  
    \path  (A'') edge node {} (B'');
    \path  (A'') edge node   {}(C'');
    \path  (B'') edge node   {} (C'');
    \path  (C'') edge[bend right=10] node {} (D'');
    \path  (C'') edge[bend left=10] node {} (D'');
   \end{scope}
  \end{tikzpicture}
   \end{tiny}
   \end{center}
\smallskip
 In case $(3A_3)'$ the distinguished vertex is the one not connected with the geometric vertex (indicated with  a bold font) and the case $(2A_3)'+A_3$ represents the 2 cases $(2A_3)'+A_3$-(i) and  $(2A_3)'+A_3$-(ii): the distinguished vertex can be one of the two on the left or the one on the right. In all other cases, the symmetry group is transitive on the virtual vertices and hence  any virtual  vertex can serve as a distinguished vertex. 

The stable graph of $C_X$ is obtained by erasing the virtual vertices and closing up the adjoining bonds.\hfill $\square$
%Two  $A_3$-singularities of $X$ are coplanar if and only if the associated virtual vertices lie in the same  connected components of the full subgraph spanned by the virtual vertices. 
\end{proposition}

\begin{remarks}\label{rem:}
Any symmetry of these  stable models is realized by a  monodromy  over  the corresponding stratum, but the converse is not true: the three $A_3$-singularities  in the case  $(3A)'$ are transitively permuted by monodromy. Indeed, if we replace $p$ by 
an other singular point $q$, then $C_q$ is in fact isomorphic with 
$C_p$ (this follows from  Proposition \ref{prop:CX}), but $q$ has now taken the position of $p$. 

All the incidence relations among the coarse strata can be  combinatorially explained in terms of these graphs: 
$\G(X',p')$ is obtained from $\G(X,p)$ by (possibly repeated)  contraction  of bonds if and only if $(X,p)$ is a degeneration of $(X',p')$. 
An  $A_3$-singularity  can degenerate into an $A_4$-singularity  if and only if the associated  virtual vertex has a double bond with a geometric vertex.
\end{remarks}

We now can finish the proof of  Theorem \ref{thm:match0}. 

\begin{proof}[Completion of the proof of Theorem \ref{thm:match0}]

Insofar as the theorem involves $A_4$-singularities, this has been taken care of  
in Subsection \ref{subsect:coarse}.  Assume  therefore that $X$ has $r\ge 1$ singular points, all  of type $A_3$. 
We have already  matched the corresponding coarse strata of $\Mcal_\ct^\st$: each of these is given by an orbit of quadratic maps
$rA_3\to \Tcal$, or equivalently, by  a $\Orth(\Tcal)$-conjugacy class of subgroups $W'< \Orth(\Tcal)$ with $W'\cong \Sfrak_4^r$.
We fix  $W_o\in\Cscr(A_3)$ and an isomorphism $C_{\Gcal}(W_o)\cong \Orth(\Tcal)$ and put $W:=W_o\times W'$. Then
$C_{\Gcal}(W)$ is identified with $C_{\Orth(\Tcal)}(W')=W'{}^\perp$.  
A point of the coarse stratum of $\Mcal_\ct^\st[\Tcal]$ over $X$, consists in specifying  a level 2 structure on $H_3(X)$.  To be precise, 
this  consists in specifying  one such $W'$ plus an isomorphism $\Hl_3(X; \FF_2)\cong \Hl_1(C_X; \FF_2)$ onto $\Tcal^{W'}$ as quadratic spaces. 

The stratum thus obtained is denoted   $\hat\Xcal\la W\ra$.  It has  an evident action of  the orthogonal group of $\Tcal^{W'}$.  Corollary \ref{cor:extendedaction} then shows that it comes with an  action of $N_{\Gcal}(W)/W\cong W'{}^\perp.\out_{\Gcal}(W)$  extending the one of  $N_{\Orth(\Tcal)}(W')/W'\cong W'{}^\perp.\out_{\Orth(\Tcal)}(W')$ and that this action is unique for that property.  

This is all that was left to prove; the rest follows from the preceding.
\end{proof}

\end{document}